# The fundamental properties characterizing the structural behaviors of Collatz sequences

**Raouf Rajab**
**National School of Engineers of Gabes- Tunisia**
**raouf.rajab@enig.rnu.tn**

## Abstract

This work represents an in-depth study of the structural behavior of the Collatz sequences. We consider a finite arithmetic progression with a common difference is equal to 2 and the number of terms in the sequence is equal to $2^n$ . After, we consider a $2^n x(n+1)$ matrix ((n+1) columns and $2^n$ rows) such as the first column contains the terms of arithmetic progression and each row of the matrix represent a finite Collatz sequence. Then, each element of the matrix will be replaced by 0 or 1 according to the following rule: the even integer is replaced by 0 and the odd integer is replaced by 1. We obtained a table contains all binary permutation with repetition this property is called the property of structural complementarily. Based on these tables, we can determine any other fundamentals properties characterizing the behavior of Collatz sequences. Thus, we can distinguish two other important properties such as the property of structural cyclical behavior characterized by a structural periodicity and the structural sliding property that can be called also the property of structural divergence.

**Résumé :** Ce travail représente une étude approfondie des comportements structurels des suites de Collatz. On considère une suite arithmétique de raison 2 et de longueur $2^n$ puis on construit une matrice ($(n+1)$ colonnes, $2^n$ lignes) tel que la première colonne contient les termes de la suite arithmétique considérée et chaque ligne du tableau représente une suite finie de Collatz. Ensuite on effectue les conversions suivantes : un entier pair sera remplacé par 0 alors un entier impair sera remplacé par 1. On obtient ainsi un tableau qui contient tous les arrangements avec répétition de deux entiers 0 et 1. Cette propriété est appelée la propriété de la complémentarité structurelle qui se traduit par ces tableaux des arrangements binaires complets caractérisées par une régularité bien déterminée. On se basant sur ces tableaux, on peut déterminer toutes autres propriétés fondamentales caractérisant le comportement structurel des suites de Collatz tel que la propriété de la périodicité structurelle qui décrit le comportement cyclique structurel de ces suites et la propriété du glissement séquentiel structurel qu'on peut l'appeler aussi la propriété de la divergence structurelle. Ces propriétés sont bien détaillées dans cet article.

# 1 Introduction

La conjecture de Syracuse-Collatz est l'un des mystérieux problèmes en mathématique. Le problème est basé sur un principe simple de calcul : on prend un entier naturel non nul quelconque si le nombre est pair on le devise par deux si non on le multiplie par 3 et on l'ajoute 1 si on poursuit le calcul on obtient une suite des nombres que finir par tomber dans un cycle trivial (1, 4, 2). La conjecture annonce que cette propriété est vraie pour tout entier naturel non nul [1]. En réalité, la conjecture de Collatz, absolument pas un problème de convergence ou de divergence mais la divergence ou la convergence ne sont que les conséquences des autres lois qui caractérisent les suites de Collatz. L'objectif de ce travail est la détermination de différentes propriétés fondamentales qui caractérisent ces genres des suites et qui sont en forte relation avec la propriété de convergence de différentes suites de Collatz. On peut montrer ainsi que les suites de Collatz présentent un certain nombre des propriétés très particulières liées à un comportement structurel extrêmement régulier. Par contre, le comportement quantitatif de ces suites est caractérisé par une irrégularité de la distribution des différentes valeurs de leurs termes. L'approche probabiliste conduite à des résultats très encourageants [1]. Comme les variations des termes des différentes suites ne sont pas aléatoires on doit trouver des méthodes non probabilistes basées sur les propriétés générales (de complémentarité ou communes) décrivant les comportements des suites du Collatz qui nous permettent d'expliquer et de prouver la propriété de convergence de toutes les suites vers un cycle bien déterminé.

*Les propriétés fondamentales étudiées et les résultats obtenus*

En réalité, on peut considérer la suite de Collatz comme une suite particulière appartenant à un genre particulier des suites entières qu'on peut l'appeler les suites structuro-entières dont l'expression de l'image d'un entier naturel non nul quelconque dépend de sa forme(ou de sa structure). Le couplage entre la forme de l'entier et l'expression de son image génère un ensemble des propriétés communes et un comportement de groupe exceptionnel pour ce genre des suites.
Pour donner une idée générale sur le contenu du ce travail, les différentes propriétés étudiées dans cet article qui caractérisent le comportement structurel des suites de Collatz sont présentées ci-dessous. Ces propriétés décrivent le comportement du groupe des suites de Collatz. En outre la détermination de différentes propriétés caractérisant les suites du Collatz constitue une étape indispensable pour de résoudre la conjecture de Collatz et d'expliquer profondément cette tendance des suites vers un cycle stable après un certain nombre des transformations successives.
On définie la suite du Collatz notée $T^s$ pour tout entier naturel non nul P comme suit [1]:

$$T^s(P) = \begin{cases} \frac{P}{2} & \text{si } P \equiv 0 \bmod 2 \\ \frac{3P+1}{2} & \text{si } P \equiv 1 \bmod 2 \end{cases} \qquad \text{(e 1.1)}$$

## 1.1 Propriété de la périodicité structurelle ou de la distribution structurelle périodique

C’est une Propriété très intéressante déjà étudiée par Terras. R. Elle traduit la similarité structurelle d’une infinité des suites dont les générateurs possédants une différence entre eux de la forme $2^n j$ avec n la longueur de ces suites isoformes et j un entier naturel non nul quelconque.
L’aspect structurel périodique des suites dites isoformes de Collatz est illustré dans le tableau suivant :

$Sy(P_i, 6)$

$P_{i+1} - P_i = 2^6$

| | | | | | |
|---|---|---|---|---|---|
| 17 | 26 | 13 | 20 | 10 | 5 |
| 81 | 122 | 61 | 92 | 46 | 23 |
| 145 | 218 | 109 | 164 | 82 | 41 |
| 209 | 314 | 157 | 236 | 118 | 59 |
| 273 | 410 | 205 | 308 | 154 | 77 |
| 337 | 506 | 253 | 380 | 190 | 95 |
| 401 | 602 | 301 | 452 | 226 | 113 |

*Figure 1:* Illustration de comportement structurel périodique des suites isoformes de Collatz

Chaque ensemble constitué d’une infinité des suites isoformes dont les premiers termes décrivent une suite arithmétique. Noter que la  raison arithmétique de la suite génératrice (constituée par les éléments de la première colonne du tableau) dépend de la longueur de différentes suites isoformes.

## 1.2 Propriété de la complémentarité structurelle : Aspect structurel parfaitement régulier versus un aspect quantitatif irrégulier

C’est une conséquence de la première propriété. La complémentarité structurelle se traduit par des tableaux d’arrangements complets. En effet l’aspect structurel des suites de Collatz peut être décrit par des tableaux des arrangements binaires complets qui traduits le comportement (structurel) du groupe des suites de Collatz de même longueur n et qui sont générées par une suite arithmétique de raison un entier naturel pair non nul et de longueur finie égale à $2^n$. Dans un tel tableau chaque suite possède sa propre distribution structurelle et les différentes distributions structurelles forment ensemble un tableau d’arrangement binaire complet.
Les conditions d’obtention d’un tel tableau sont simples et on peut les résumer en quatre étapes :
-On considère une suite arithmétique de raison 2 et de premier terme un entier naturel non nul quelconque. La suite est de longueur finie $2^n$ avec n un entier naturel non nul quelconque.
- Pour chaque terme de la suite déjà définie en détermine les n premiers termes de la suite de Collatz.

-La représentation matricielle consiste à construire un tableau de (n+1) colonnes et de $2^n$ lignes tel que la première colonne contient les termes de la suite arithmétique considérée qu'on l'appelle aussi la suite génératrice.
-On fait remplacer chaque terme des suites obtenues par 0 s'il s'agit d'un entier pair ou bien par 1 s'il s'agit d'un entier impair. On peut colorer les cases selon la parité des entiers placées dans les différentes cases au lieu d'utiliser les nombres 0 et 1. Dans les deux cas on obtient tous les arrangements possibles avec répétition de longueur n des deux objets.
Par exemple la conversion de la suite de Collatz générée par 7 et de longueur 6 est comme suit :

$$(7,11,17,26,13,20) \rightarrow (1,1,1,01,0) \qquad \text{(e 1.2)}$$

Lors de ces études sur les suites de Collatz le mathématicien R. Terras est appelé cette suite binaire qui traduit la distribution des entiers impairs et pairs dans une suite donnée le vecteur de parité tout en formalisant quelques propriétés concernant ces vecteurs de parité.

**Exemple 1.1**

L'exemple suivant illustre la construction d'un tableau des arrangements binaires complets à partir d'une matrice de Collatz de dimension bien déterminée selon les règles déjà décrites

| | | | | |
|---|---|---|---|---|
| 1 | 2 | 1 | 2 | 1 |
| 3 | 5 | 8 | 4 | 2 |
| 5 | 8 | 4 | 2 | 1 |
| 7 | 11 | 17 | 26 | 13 |
| 9 | 14 | 7 | 11 | 17 |
| 11 | 17 | 26 | 13 | 20 |
| 13 | 20 | 10 | 5 | 8 |
| 15 | 23 | 35 | 53 | 80 |
| 17 | 26 | 13 | 20 | 10 |
| 19 | 29 | 44 | 22 | 11 |
| 21 | 32 | 16 | 8 | 4 |
| 23 | 35 | 53 | 80 | 40 |
| 25 | 38 | 19 | 29 | 44 |
| 27 | 41 | 62 | 31 | 47 |
| 29 | 44 | 22 | 11 | 17 |
| 31 | 47 | 71 | 107 | 161 |

n colonnes

| | n colonnes | | | |
|---|---|---|---|---|
| | 0 | 1 | 0 | 1 |
| | 1 | 0 | 0 | 0 |
| | 0 | 0 | 0 | 1 |
| | 1 | 1 | 0 | 1 |
| | 0 | 1 | 1 | 1 |
| | 1 | 0 | 1 | 0 |
| | 0 | 0 | 1 | 0 |
| | 1 | 1 | 1 | 0 |
| | 0 | 1 | 0 | 0 |
| | 1 | 0 | 0 | 1 |
| | 0 | 0 | 0 | 0 |
| | 1 | 1 | 0 | 0 |
| | 0 | 1 | 1 | 0 |
| | 1 | 0 | 1 | 1 |
| | 0 | 0 | 1 | 1 |
| | 1 | 1 | 1 | 1 |

$2^n$ Lignes

*Figure 2:* Tableau d'arrangement binaire complet

Noter que la colonne génératrice (première colonne colorée en rouge brique) ne fait pas partie de ce tableau d'arrangements complets mais elle permet de le créer. Le tableau de vecteurs de parité relatives aux différentes suites générées contient tous les arrangements avec répétition de deux entiers 0 et 1. La démonstration de cette propriété fait l'objet de la partie consacrée pour l'étude des chromoformes parfaits de Collatz.

## 1.3 Propriété du glissement séquentiel structurel : une divergence structurelle versus une convergence entière.

Cette propriété peut être appelée aussi la propriété de la divergence structurelle c'est le phénomène opposé de la convergence des termes numériques des suites de Collatz vers un cycle stable. On sait que chaque suite finie de Collatz possède une distribution structurelle bien déterminée (vecteur de parité). Le processus inverse consiste à considéré un arrangement binaire (ou chromatique avec deux couleurs) avec répétition puis on cherche à déterminé une suite de Collatz admettant une structure binaire identique à l'arrangement considéré. On peut distinguer deux cas différents :
-Si la longueur de l'arrangement considéré est finie dans ce cas, ce dernier admet une infinité des suites de Collatz dont leurs structures sont identiques à cet arrangement.
-Si on considère des structures binaires de longueurs infinies, on remarque que certaines structures ne correspondent pas à aucunes distribuions structurelles d'aucunes suites on dit que ce genre des structures (ou arrangements binaires) sont non réalisables ou non convertibles en suites de Collatz. L'exemple suivant illustre cette propriété pour une structure cyclique constituée par la même séquence.

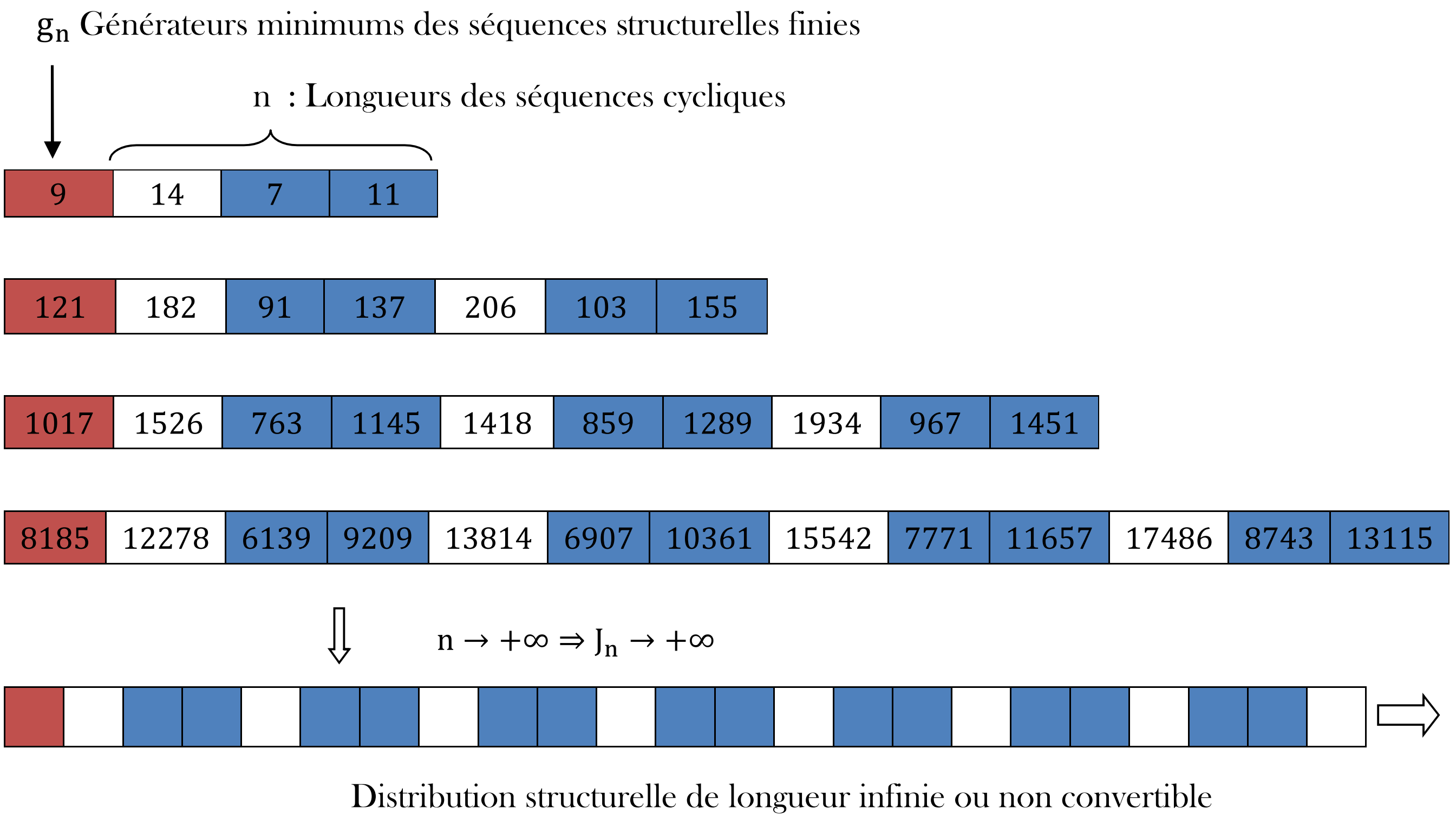

*Figure 3 :* Un exemple illustrant la propriété de la divergence structurelle

On remarque que lorsque les séquences finies tendent vers la structure de longueur infinie, les générateurs des suites tendent aussi vers l'infini.

Donc elle n'existe pas aucune suite de Collatz de longueur infinie tel que sa distribution structurelle chromatique est identique à cette distribution. On montre dans ce travail que la quasi-totalité des distributions chromatiques (ou binaires) de longueurs infinies sont non convertibles en suites de Collatz.

# 2 Préliminaire et notations

**Définition 2.1**

Pour tout entier naturel non nul N, on définit la transformation élémentaire de Collatz notée $T^S$ comme suit:

$$T^S(N) = \begin{cases} \frac{3}{2}N + \frac{1}{2} & \text{si N est impair} \\ \frac{1}{2}N & \text{si N est pair} \end{cases} \quad \text{(e 2.1)}$$

**Définition 2.2**

L'image d'ordre k de N est l'entier obtenu par application k fois successives de la transformation $T^S$ sur N, il est noté $T_k^S(N)$ et il est défini comme suit :

$$\underbrace{T^S\left(T^S\left(\dots\left(T^S(N)\right)\right)\right)}_{\text{k fois}} = T_k^S(N) \quad \text{(e 2.2)}$$

L'indice k dans l'expression $T_k^S(N)$ désigne le nombre des fois d'application de la transformation $T^S$ à l'entier de départ N.
-Noter qu'on ne fait pas distinction entre les deux notations suivantes :

$$T_1^S(N) = T^S(N)$$

**Propriété 2.1**

Pour tous entiers naturels i et j on peut écrire :

$$T_i^S(T_j^S(N)) = T_{i+j}^S(N) \quad \text{(e 2.3)}$$

**Définition 2.3**

On définit les effets transformationnels élémentaires principaux notés $F_1(N)$ relatifs aux deux formes pair et impair comme suit :

$$F_1(N) = \begin{cases} \frac{1}{2} & \text{si N est pair} \\ \frac{3}{2} & \text{si N est impair} \end{cases} \quad \text{(e 2.4)}$$

**Définition 2.4**

On définit les effets transformationnels élémentaires secondaires notés $\varphi_1(N)$ comme suit :

$$\varphi_1(N) = \begin{cases} 0 & \text{si N est pair} \\ \frac{1}{2} & \text{si N est impair} \end{cases} \quad \text{(e 2.5)}$$

**Définition 2.5**

Le couple des effets transformationnels élémentaires est noté $E_1$ est définit comme suit :

$$E_1(N) = \left(F_1(N), \varphi_1(N)\right) \quad \text{(e 2.6)}$$

$E_1(N)$ est appelé aussi vecteur des effets transformationnels élémentaires. On peut l'exprimer comme suit :

$$E_1(N) = \begin{cases} \left(\frac{1}{2}, 0\right) & \text{si N est pair} \\ \left(\frac{3}{2}, \frac{1}{2}\right) & \text{si N est impair} \end{cases} \quad \text{(e 2.7)}$$

**Remarque 2.1** Les effets transformationnels élémentaires principaux ou secondaires sont appelées aussi les coeffcients transformationnels élémentaires de la transformation du Collatz.

**Définition 2.6**

Soit P un entier naturel non nul, pour tout entier naturel n ($n \geq 0$ ) on définit l'indicateur de parité d'ordre n de P qu'on le not $\mathbf{i}_n(P)$ comme suit :

$$\mathbf{i}_n(P) = \begin{cases} 1 \text{ si } T_n^S(P) \text{ est impair} \\ 0 \text{ si non} \end{cases} \quad \text{(e 2.8)}$$

**Exemple 2.1**

On prend le cas de N=27:

$$\mathbf{i}_0(27) = 1$$
$$T_1^S(27) = 41 \implies \mathbf{i}_1(27) = 1$$
$$T_2^S(27) = 62 \implies \mathbf{i}_2(27) = 0$$

**Propriété 2.2**

Soient P un entier naturel non nul, n et m deux entiers naturels quelconques donc l'indicateur de parité d'ordre m de $T_n^S(P)$ correspond à l'indicateur d'ordre (n+m) de P c'est à dire :

$$\mathbf{i}_m\left(T_n^S(P)\right) = \mathbf{i}_{n+m}(P) \quad \text{(e 2.9)}$$

En particulier: $\mathbf{i}_0\left(T_n^S(P)\right) = \mathbf{i}_n(P)$

**Démonstration**

$\mathbf{i}_m\left(T_n^S(P)\right)$ est l'indicateur de parité d'ordre m de $T_n^S(P)$ , il est défini comme suit :

$$\begin{cases} \mathbf{i}_m\left(T_n^S(P)\right) = 1 & \text{si } T_m^S(T_n^S(P)) \text{ est impair} \\ \mathbf{i}_m\left(T_n^S(P)\right) = 0 & \text{si non} \end{cases}$$

Comme $T_m^S(T_n^S(P)) = T_{m+n}^S(P)$ alors on peut écrire :

$$\begin{cases} \mathbf{i}_m\left(T_n^S(P)\right) = 1 & \text{si } T_{m+n}^S(P) \text{ est impair} \\ \mathbf{i}_m\left(T_n^S(P)\right) = 0 & \text{si non} \end{cases}$$

Par définition, on sait que :

$$\begin{cases} \mathbf{i}_{m+n}(P) = 1 & \text{si } T^S_{m+n}(P) \text{ est impair} \\ \mathbf{i}_{m+n}(P) = 0 & \text{si non} \end{cases}$$

On peut conclure que :

$$\mathbf{i}_m\left(T_n^{S1}(P)\right) = \mathbf{i}_{n+m}(P)$$

**Définition 2.7**

Soit P un entier naturel non nul et n un entier naturel quelconque ($n \geq 0$) donc le nombre d'imparité d'ordre n de P qu'on le note $M_n(P)$ est le nombre des entiers impairs contenus dans une suite de la forme : $(P, T_1^S(P), T_2^S(P), \ldots, T_n^S(P))$. Autrement, on peut exprimer $M_n(P)$ en fonction des indicateurs de parité comme suit :

$$M_n(P) = \sum_{k=0}^{k=n} \mathbf{i}_k(P) \qquad \text{(e 2.10)}$$

**Exemple 2.2**

Considérons la suite ci-dessous :

| | | | | | | |
|---|---|---|---|---|---|---|
| 29 | 44 | 22 | 11 | 17 | 26 | 13 |

Cette suite contient 4 entiers impairs : 29, 11,17 et 13 donc on peut déduire que :

$\Longrightarrow M_6(29) = 4$

**Corollaire 2.1**

Soit P un entier naturel non nul alors pour tout entier naturel non nul n, le nombre d'imparité d'ordre n de P vérifie la relation de récurrence suivante :

$$M_n(P) = \sum_{k=0}^{n-1} \mathbf{i}_k(P) + \mathbf{i}_n(P) = M_{n-1}(P) + \mathbf{i}_n(P) \qquad \text{(e 2.11)}$$

En particulier pour k=0 on peut écrire :

$$M_0(N) = \mathbf{i}_0(N)$$

**Proposition 2.1**

On sait que la transformation de Syracuse est définie comme suit pour tout entier naturel non nul N :

$$T^S(N) = \begin{cases} \frac{3}{2}N + \frac{1}{2} & \text{si N est impair} \\ \frac{1}{2}N & \text{si N est pair} \end{cases}$$

Les deux expressions ci –dessus peuvent être remplacées par une expression unique en fonction de l'indicateur de parité comme suit :

$$T_1^S(N) = \frac{3^{\mathbf{i}_0(N)}}{2}N + \frac{1}{2}\mathbf{i}_0(N) \qquad \text{(e 2.12)}$$

Cette expression est appelée l'expression unificatrice de $T^S(N)$

**Proposition 2.2**

Revenons aux expressions des effets transformationnels élémentaires de $T^S$. Les deux expressions de l'effet transformationnel principal élémentaire peuvent être remplacées par la forme unificatrice suivante:

$$F_1(N) = \frac{3^{\mathbf{i}_0(N)}}{2} \tag{e 2.13}$$

Les deux expressions de l'effet transformationnel secondaire élémentaire peuvent être remplacées par la forme unificatrice suivante :

$$\varphi_1(N) = \frac{1}{2}\mathbf{i}_0(N) \tag{e 2.14}$$

**Proposition 2.3**

On se basant sur les deux expressions unificatrices ci-dessus, on peut définir une expression unificatrice pour le vecteur des effets élémentaires. Cette expression unificatrice du vecteur des effets est comme suit:

$$E_1(N) = \left(\frac{3^{\mathbf{i}_0(N)}}{2}, \frac{1}{2}\mathbf{i}_0(N)\right) \tag{e 2.15}$$

**Proposition 2.4**

Soient A et B deux entiers naturels non nuls, l'image de l'entier naturel (A+2B) par $T^S$ est comme suit :

$$T^s(A + 2B) = T^s(A) + 3^{\mathbf{i}_0(A)}\, B \tag{e 2.16}$$

**Démonstration**

La parité de (A+2B) dépend de la parité de A en effet :

-Si A est pair donc A+2B est pair aussi ce qui nous permet d'écrire:

$$T^s(A + 2B) = \frac{A + 2B}{2} = \frac{A}{2} + B = T^S(A) + B \tag{e 2.17}$$

-Si A est impair alors (A+2B) est impair aussi, on écrit alors :

$$T^s(A + 2B) = \frac{3}{2}(A + 2B) + \frac{1}{2} = \left(\frac{3}{2}A + \frac{1}{2}\right) + 3B = T^S(A) + 3B \tag{e 2.18}$$

Comme on a :

$$\begin{cases} 3^{\mathbf{i}_0(A)}B = B & \text{si A est pair} \\ 3^{\mathbf{i}_0(A)}B = 3B & \text{si non} \end{cases} \tag{e 2.19}$$

Donc dans les deux cas, l'image de (A+2B) par $T^s$ peut exprimer comme suit :

$$T^s(A + 2B) = T^s(A) + 3^{\mathbf{i}_0(A)}B$$

**Théorème 2.1**

Soit P un entier naturel non nul. Pour tous entier naturel non nul n ($n \geq 1$), l'expression de $T_n^S(P)$ en fonction de P est comme suit :

$$T_n^S(P) = \frac{3^{M_{n-1}(P)}}{2^n}\, P + \varphi_n(P) \tag{e 2.20}$$

Avec :

-$F_n(P)$ est le coefficient transformationnel cumulatif principal d'ordre n (ou du rang n) de P par la transformation $T^S$ , il est exprimé comme suit :

$$F_n(P) = \frac{3^{M_{n-1}(P)}}{2^n} \qquad \text{(e 2.21)}$$

-$\varphi_n(P)$ est le coefficient transformationnel cumulatif secondaire d'ordre n de P par $T^S$, il est définie comme suit :

$$\begin{cases} \varphi_1(P) = \frac{1}{2}\mathbf{i}_0(P) \\ \varphi_n(P) = \frac{3^{\mathbf{i}_{n-1}(P)}}{2}\varphi_{(n-1)}(P) + \frac{1}{2}\mathbf{i}_{n-1}(P) \end{cases} \qquad \text{(e 2.22)}$$

**Remarque 2.1**

-L'effet tranformationnel cumulatif principal d'ordre n noté $F_n$ traduit tous les effets élémentaires principaux relatifs aux termes du rangs allant du 0 à $(n-1)$.

-L'effet transformationnel cumulatif secondaire d'ordre n noté $\varphi_n$ traduit tous les effets élémentaires relatifs aux termes du rangs allant du 0 à $(n-1)$.

**Démonstration**

On montre par récurrence que cette propriété est vraie pour tout entier naturel non nul n.

L'expression unificatrice de $T_1^S(P)$ s'écrit comme suit:

$$T_1^S(P) = \frac{3^{\mathbf{i}_0(P)}}{2}P + \frac{1}{2}\mathbf{i}_0(P)$$

On sait que $\mathbf{i}_0(P) = M_0(P)$ et $\varphi_1(P) = \frac{1}{2}\mathbf{i}_0(P)$ ce qui nous permet d'écrire :

$$T_1^S(P) = \frac{3^{M_0(P)}}{2}P + \varphi_1(P) \qquad \text{(e 2.23)}$$

Donc la propriété est vraie pour n=1.

On suppose que la propriété est vraie pour tout entier k allant de 2 à n et on montre que la propriété est vraie pour $k = (n+1)$.

On sait que :

$$T_{n+1}^S(P) = T_1^S\left(T_n^S(P)\right) = T_1^S\left(\frac{3^{M_{n-1}(P)}}{2^n}P + \varphi_n(P)\right) \qquad \text{(e 2.24)}$$

On pose :

$$P_n = \frac{3^{M_{n-1}(P)}}{2^n}\, P + \varphi_n(P) \tag{e 2.25}$$

On peut écrire alors :

$$T_1^S\left(\frac{3^{M_{n-1}(P)}}{2^n}\, P + \varphi_n(P)\right) = T_1^S(P_n)$$

$$= \frac{3^{\mathbf{i}_0(P_n)}}{2} P_n + \frac{1}{2}\mathbf{i}_0(P_n) \tag{e 2.26}$$

En utilisant les propriétés d'additivité des indicateurs de parité on peut déduire la relation suivante :

$$\mathbf{i}_0(P_n) = \mathbf{i}_0\left(T_n^S(P)\right) = \mathbf{i}_n(P)$$

Remplaçons $\mathbf{i}_0(P_n)$ par $\mathbf{i}_n(P)$ dans l'équation ci-dessus puis $P_n$ par son expression en fonction de P, on obtient les équations suivantes :

$$T_1^S\left(\frac{3^{M_{n-1}(P)}}{2^n}\, P + \varphi_n(P)\right) = \frac{3^{\mathbf{i}_n(P)}}{2} P_n + \frac{1}{2}\mathbf{i}_n(p)$$

$$= \frac{3^{\mathbf{i}_n(P)}}{2}\left(\frac{3^{M_{n-1}(P)}}{2^n}\, P + \varphi_n(P)\right) + \frac{1}{2}\mathbf{i}_n(P)$$

$$= \frac{3^{(M_{n-1}(P)+\mathbf{i}_n(P))}}{2^{n+1}} P + \frac{3^{\mathbf{i}_n(P)}}{2}\varphi_n(P) + \frac{1}{2}\mathbf{i}_n(P) \tag{e 2.27}$$

On sait que $M_n(P) = M_{n-1}(P) + \mathbf{i}_n(P)$ ce qui nous permet d'écrire :

$$T_1^S\left(\frac{3^{M_{n-1}(P)}}{2^n}\, P + \varphi_n(P)\right) = \frac{3^{M_n(P)}}{2^{n+1}} N + \varphi_{n+1}(P) \tag{e 2.28}$$

Avec:

$$\varphi_{n+1}(N) = \frac{3^{\mathbf{i}_n(P)}}{2}\varphi_n(P) + \frac{1}{2}\mathbf{i}_n(P) \tag{e 2.29}$$

La propriété est vraie pour $k = (n + 1)$ ce qui nous permet de conclure que la propriété est vraie pour tout entier  naturel non nul n.

# 3 Les suites générées de Collatz

**Notation 3.1**

Pour tout entier naturel non nul P, on définit la suite générée de Syracuse de longueur n comme suit :

$$\tilde{S}^{r}(P,n) = \left(T_1^S(P), T_2^S(P), \dots, T_n^S(P)\right) \qquad \text{(e 3.1)}$$

Avec:

-n un entier naturel non nul, il correspond à la longueur de la suite générée.

-P n'appartient pas à la suite $\tilde{S}^{r}(P,n)$, il est appelé le générateur direct de cette suite générée.

**Notation 3.2**

Si la suite générée est de longueur infinie ($n \to +\infty$) donc on la note comme suit :

$$\tilde{S}^{R}(P) = \left(T_1^S(P), T_2^S(P), \dots, T_n^S(P), \dots\right) \qquad \text{(e 3.2)}$$

**Notation 3.3**

On exclut le générateur P de la suite générée puisque les tableaux d'arrangements binaires complets sont constituées uniquement par des suites sous cette forme le générateur est appartient à la suite génératrice permettant de créer mais ceci nous n'empêche pas de définir une suite qui contient le générateur, elle est comme suit :

$$Sy(P,n+1) = \left(P, T_1^S(P), T_2^S(P), \dots, T_n^S(P)\right) \qquad \text{(e 3.3)}$$

Dans quelques cas, on fait recours à l'utilisation des suites $Sy(P,n)$ pour déterminer ou montrer une propriété en relation avec les suites générées par exemple.

On remplace le r par grand R pour designer une suite de longueur infinie.

**Corollaire 3.1**

Le nombre d'imparité de la suite $\tilde{S}^{r}(P,n)$ correspond au nombre des entiers impairs contenus dans la suite $\tilde{S}^{r}(P,n)$, il est noté $\widetilde{M}_n(P)$. Comme l'entier P n'appartient pas à la suite $\tilde{S}^{r}(P,n)$ donc, il est évident que, le nombre d'imparité de cette suite est égal à :

$$\widetilde{M}_n(P) = M_n(P) - \mathbf{i}_0(P) \qquad \text{(e 3.4)}$$

**Exemple 3.1**

-La suite générée de générateur direct 27 et de longueur 5 notée $\tilde{S}^{r}(27,5)$ est comme suit:

$$\tilde{S}^{r}(27,5) = (41,62,31,47,71)$$

La suite peut être représentée sous forme d'une ligne à (n+1) cases. La première case est consacrée pour le générateur qui doit être séparé du reste de la ligne par un trait gras pour indiquer qu'il n'appartient pas à la suite considérée.

| | | | | | |
|---|---|---|---|---|---|
| 27 | 41 | 62 | 31 | 47 | 71 |

La suite générée contient 4 entiers impairs : 41, 31,47 et 71 donc on peut écrire :

$$\widetilde{M}_5(27) = M_5(27) - \mathbf{i}_0(27) = 4$$

Noter que $M_5(27) = 5$

## 3.1 Représentation d'une suite générée sous forme d'une ligne à n cases :

La suite $\tilde{S}^r(P, n)$ peut être représentée sous forme d'une ligne à n cases, on peut ajouter si nécessaire deux autres cases supplémentaires qui contiennent le générateur P et le terme qui suit le dernier terme de la suite c'est à dire le terme $T^S_{n+1}(P)$. Ces deux entiers sont appelés les extermes de la suite générée. L'entier $T^S_{n+1}(P)$ est appelée le preterme de la suite générée.

Les deux cases qui contiennent les extermes doivent être séparées du reste de la ligne par deux traits gras comme montre la figure ci-dessous :

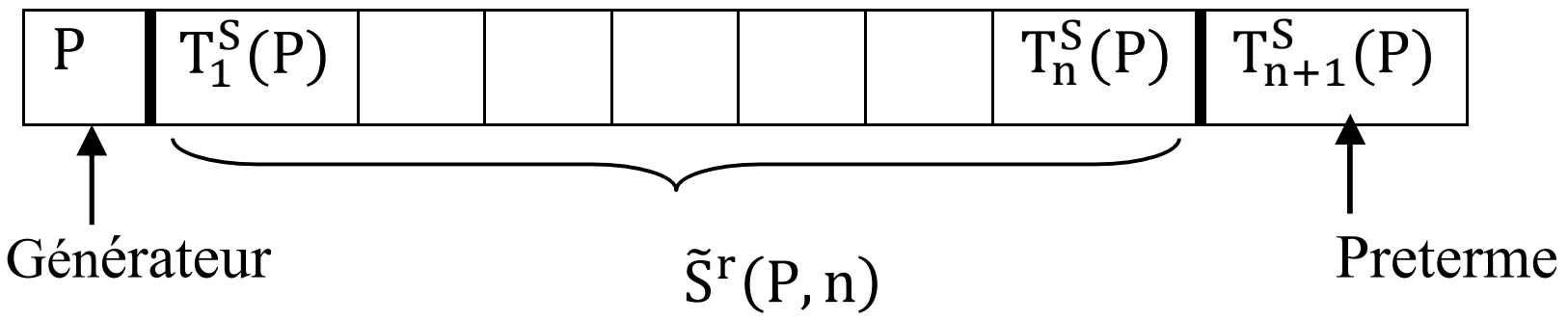

*Figure 4:* La suite générée et leurs extermes

## 3.2 Importance des notions du générateur et du preterme d'une suite générée finie :

Le preterme et le générateur d'une suite générée finie ont une importance particulière, en effet

-Si on veut déterminer si une suite générée est croissante ou décroissante on doit comparer son premier terme et son preterme c'est à dire $T^S_1(P)$ et $T^S_{n+1}(P)$ du fait que le preterme traduit tous les effets transformationnels élémentaires relatifs aux tous les termes de la suite générée.

-Si on cherche à déterminer si deux suites générées sont isoformes ou non (c'est à dire qu'elles possèdent ou non la même distribution des indicateurs des parités) donc on doit comparer leurs générateurs comme va le voir.

## 3.3 Importance de la notion des suites générées de Collatz:

Les suites générées de Syracuse-Collatz ont une importance primordiale dans la détermination des toutes propriétés structurelles de ce genre particulier des suites. L'une de ces propriétés (qu'on va démontrer plus tard) peut être décrite comme suit :

- On fait construire une suite arithmétique de raison 2 et de longueur $2^n$.
- On fait correspondre à chaque terme $P_i$ de la suite arithmétique considérée, la suite générée de générateur $P_i$ et de longueur n :$\tilde{S}^r(P_i, n)$.
- Les suites sont représentées sous forme d'un tableau à n colonnes et de $2^n$ lignes
- Les termes sont convertis en leurs indicateurs de parité, on obtient un tableau des arrangements binaires complets périodiques.

**Exemple 3.2**

La figure ci-dessous représente un simple exemple qui illustre cette propriété :

Suite génératrice ↓

| $P_i$ | $\tilde{S}^r(P_i, n)$ | | |
|---|---|---|---|
| 1 | 2 | 1 | 2 |
| 3 | 5 | 8 | 4 |
| 5 | 8 | 4 | 2 |
| 7 | 11 | 17 | 26 |
| 9 | 14 | 7 | 11 |
| 11 | 17 | 26 | 13 |
| 13 | 20 | 10 | 5 |
| 15 | 23 | 35 | 53 |

Remplacement de $n_{ij}$ par $\mathbf{i}_j(P_i)$ ⇨

| $P_i$ | $\mathbf{i}_1(P_i)$ | $\mathbf{i}_2(P_i)$ | $\mathbf{i}_3(P_2)$ |
|---|---|---|---|
| 1 | 0 | 1 | 0 |
| 3 | 1 | 0 | 0 |
| 5 | 0 | 0 | 0 |
| 7 | 1 | 1 | 0 |
| 9 | 0 | 1 | 1 |
| 11 | 1 | 0 | 1 |
| 13 | 0 | 0 | 1 |
| 15 | 1 | 1 | 1 |

*Figure 5:* Aspect structurel des suites de Syracuse-Collatz

On remarque bien que le tableau contient tous les arrangements binaires (ou chromatique) avec répétition . Cette propriété est générale, en effet quelque soit l'entier naturel n choisi au départ, on obtient toujours des tableaux des arrangements complets.

## 3.4 Détermination des coefficients transformationnels cumulatifs relatifs à une suite générée de Collatz :

On considère la suite générée de Syracuse suivante :

$$\tilde{S}^r(P, n) = (T_1^S(P), T_2^S(P), \dots, T_n^S(P))$$

Le terme $T^S_{n+1}(P)$ (le terme qui suit le dernier terme de la suite générée) traduit tous les effets transformationnels cumulatifs principaux et secondaires relatifs aux différents termes de la suite générée considérée. C'est pour cette raison, pour caractériser une suite générée quelconque de Collatz, on fait comparer son premier terme $T^S_1(P)$ et son preterme $T^S_{n+1}(P)$. Du ce effet, on doit déterminer l'expression de $T^S_{n+1}(P)$ en fonction de $T^S_1(P)$.

**Théorème 3.1**

Soit P un entier naturel non nul. Pour tous entier naturel non nul n $(n \geq 1)$, l'expression de $T^S_{n+1}(P)$ en fonction de $T^S_1(P)$ est comme suit :

$$T^S_n(P) = \frac{3^{M_{n-1}(P)}}{2^n} P + \varphi_n(P) \qquad \text{(e 3.5)}$$

Avec :

-$\tilde{F}_n(P)$ est le coefficient transformationnel cumulatif principal d'ordre n (ou du rang n) de la suite générée $\tilde{S}^r(P, n)$ . Il est exprimé comme suit :

$$\tilde{F}_n(P) = \frac{3^{\tilde{M}_n(P)}}{2^n} \qquad \text{(e 3.6)}$$

-$\tilde{\varphi}_n(P)$ est le coefficient transformationnel cumulatif secondaire d'ordre n de la suite générée $\tilde{S}^r(P, n)$. Il est définie comme suit :

$$\tilde{\varphi}_n(P) = \varphi_{(n+1)}(P) - \frac{3^{\tilde{M}_n(P)}}{2^{n+1}} \mathbf{i}_0(P) \qquad \text{(e 3.7)}$$

**Démonstration**

On sait que :

$$T^S_{n+1}(P) = \frac{3^{M_n(P)}}{2^{n+1}} P + \varphi_{(n+1)}(P)$$

L'expression de P en fonction de $T^S_1(P)$ est comme suit :

$$T^S_1(P) = \frac{3^{\mathbf{i}_0(P)}}{2} P + \frac{1}{2} \mathbf{i}_0(P)$$

Donc l'expression de P en fonction de $T^S_1(P)$ s'écrit comme suit:

$$P = \frac{2T^S_1(P) - \mathbf{i}_0(P)}{3^{\mathbf{i}_0(P)}} \qquad \text{(e 3.8)}$$

On fait remplacer P par son expression en fonction de $T^S_1(P)$, on peut écrire :

$$T_{n+1}^{S}(P)=\frac{3^{M_n(P)}}{2^{n+1}}P+\varphi_{(n+1)}(P)=\frac{3^{(M_n(P)-\mathbf{i}_0(P))}}{2^n}T_1^S(P)+\varphi_{(n+1)}(P)-\frac{3^{(M_n(P)-\mathbf{i}_0(P))}}{2^{n+1}}\mathbf{i}_0(P)$$

Ce qui nous permet d'obtenir l'expression de $T_{n+1}^{S}(P)$ en fonction de $T_1^S(P)$ comme suit :

$$T_{n+1}^{S}(P)=\frac{3^{(M_n(P)-\mathbf{i}_0(P))}}{2^n}T_1^S(P)+\left(\varphi_{(n+1)}(P)-\frac{3^{(M_n(P)-\mathbf{i}_0(P))}}{2^{n+1}}\mathbf{i}_0(P)\right)\quad \text{(e 3.9)}$$

On sait que :

$$M_n(P)-\mathbf{i}_0(P)=\widetilde{M}_n(P)$$

Ce qui implique :

$$T_{n+1}^{S}(P)=\frac{3^{\widetilde{M}_n(P)}}{2^n}T_1^S(P)+\left(\varphi_{(n+1)}(P)-\frac{3^{\widetilde{M}_n(P)}}{2^{n+1}}\mathbf{i}_0(P)\right)\quad \text{(e 3.10)}$$

Les deux coefficients transformationnels globaux (principal et secondaire) de la suite générée $\widetilde{S}^r(P,n)$ sont comme suit :

-Le coefficient transformationnel global principal est comme suit :

$$\widetilde{F}_n(P)=\frac{3^{\widetilde{M}_n(P)}}{2^n}$$

-Le coefficient transformationnel global secondaire est le suivant :

$$\widetilde{\varphi}_n(P)=\varphi_{(n+1)}(P)-\frac{3^{\widetilde{M}_n(P)}}{2^{n+1}}\mathbf{i}_0(P)$$

## 3.5 Classement des suites générées de Syracuse-Collatz

Première classement

A fin de faciliter, l'étude des suites de Collatz, on les classe en deux grandes catégories on fait comparer leurs coefficients transformationnels principaux globaux par rapport à 1 ce qui nous permet de distinguer deux types des suites générées de Collatz selon ce classement :

**Définition 3.1**

Les suites générées de Syracuse de longueur n qui ont des coefficients transformationnels globaux principaux strictement supérieurs à 1 sont appelées les suites de type $A_n$ autrement le coefficient global principal d'une suite générée type A vérifie la condition suivante :

$$\tilde{F}_n(P) = \frac{3^{\tilde{M}_n(P)}}{2^n} > 1 \qquad \text{(e 3.11)}$$

**Définition 3.2**

Les suites générées de Syracuse de longueur n qui ont des coefficients transformationnels principaux globaux strictement inférieurs à 1 sont appelées les suites de type $B_n$. Cette condition se traduit par l'inéquation suivante:

$$\tilde{F}_n(P) = \frac{3^{\tilde{M}_n(P)}}{2^n} < 1 \qquad \text{(e 3.12)}$$

Sous classement des suites générées de type $B_n$:

On fait comparer le preterme et le premier terme d'une suite générée de Syracuse de type $B_n$ et on se basant sur l'équation on peut classer les suites $B_n$ en deux sous catégories:

**Définition 3.3**

C'est sont des suites générées de type $B_n$et de longueur n qui vérifient les deux conditions suivantes :

$$\begin{cases} \frac{3^{\tilde{M}_n(P)}}{2^n} < 1 \\ T^S_{n+1}(P) \le T^S_1(P) \end{cases} \qquad \text{(e 3.13)}$$

Les suites générées de type $B_n$ qui vérifient $T^S_{n+1}(P) = T^S_1(P)$ donc elles sont incluses dans la classe $B^-_n$.

**Définition 3.4**

Elles sont des suites grossissantes de longueur n et de type $B_n$. Elles vérifient les deux conditions suivantes :

$$\begin{cases} \frac{3^{\tilde{M}_n(P)}}{2^n} < 1 \\ T^S_{n+1}(P) > T^S_1(P) \end{cases} \qquad \text{(e 3.14)}$$

Second classement des suites de Collatz

**Définition 3.5**

Les suites générées de Collatz de type $S^+_n$ sont des suites des longueurs n tel que leurs pretermes sont strictement supérieurs à leurs premiers termes autrement si on désigne par $\tilde{S}^r(P, n)$ la suite de Collatz de générateur P et de longueur n et de type $S^+_n$ alors :

$$T^S_{n+1}(P) > T^S_1(P)$$

En réalité cette classe regroupe les deux types des suites les suites de type $A_n$ et les suites de type $B_n^+$

**Définition 3.6**

Une suite générée de Collatz est de type $S_n^-$ est une suite longueurs n tel que son preterme est inférieur à son premier terme autrement si on désigne par $\tilde{S}^r(P,n)$ la suite de générateur P et de longueur n et de type $S_n^-$ alors on elle vérifie :

$$T_{n+1}^S(P) \le T_1^S(P)$$

**Définition 3.7**

On considère la suite générée $\tilde{S}^r(P,n)$ de premier terme P et de longueur n. Le coefficient transformationnel global de cette suite est comme suit :

$$\tilde{F}_n(P) = \frac{3^{\tilde{M}_n(P)}}{2^n}$$

Le coefficient de renversement noté $\alpha_n$ relatif à des suites générées de Collatz de longueur n est défini comme suit :

$$\begin{cases} \dfrac{3^{\tilde{M}_n(P)}}{2^n} < 1 \quad \text{si } 0 \le \tilde{M}_n(P) \le \alpha_n \\ \dfrac{3^{\tilde{M}_n(P)}}{2^n} > 1 \quad \text{si } \alpha_n + 1 \le \tilde{M}_n(P) \le n \end{cases} \qquad \text{(e 3.15)}$$

**Proposition 3.1**

Pour toutes les suites générées de Collatz de longueur n avec n un entier naturel non nul quelconque l'expression du point de renversement $\alpha_n$ est la suivante :

$$\alpha_n = E\left(\frac{Ln(2)}{Ln(3)} n\right) \qquad \text{(e 3.16)}$$

C'est à dire que :

$\alpha_n$ correspond à la partie entière de $\dfrac{Ln(2)}{Ln(3)} n$

**Démonstration**

Pour déterminer l'expression du coefficient d'inversement en fonction de la longueur n des suites générées, on procède comme suit :

-si $0 \leq \mathrm{k} \leq \alpha_{\mathrm{n}}$:

$$\frac{3^{\mathrm{k}}}{2^{\mathrm{n}}} < 1 \Longrightarrow \mathrm{kLn}(3) < nLn(2) \Longrightarrow \mathrm{k} < \frac{\mathrm{Ln}(2)}{\mathrm{Ln}(3)}\mathrm{n}$$

-si $\alpha_{\mathrm{n}} < k \leq n$:

$$\frac{3^{\mathrm{k}}}{2^{\mathrm{n}}} > 1 \Longrightarrow \mathrm{kLn}(3) > nLn(2)$$

Ce qui implique:

$$\mathrm{k} > \frac{\mathrm{Ln}(2)}{\mathrm{Ln}(3)}\mathrm{n}$$

$\alpha_{\mathrm{n}}$ est le plus grand entier naturel inférieur à $\frac{\mathrm{Ln}(2)}{\mathrm{Ln}(3)}\mathrm{n}$ ce qui nous permet de conclure que $\alpha_{\mathrm{n}}$correspond à la partie entière de $\frac{\mathrm{Ln}(2)}{\mathrm{Ln}(3)}\mathrm{n}$ , on peut écrire alors:

Le point de renversement correspond à la partie entière

$$\alpha_{\mathrm{n}} = \mathrm{E}\left(\frac{\mathrm{Ln}(2)}{\mathrm{Ln}(3)}\mathrm{n}\right)$$

**Exemple 3.3**

Prenons le cas des suites générées de longueur 8.

$$\mathrm{r}_8 = 8\mathrm{x}\,\frac{\mathrm{Ln}(2)}{\mathrm{Ln}(3)} = 5{,}04 \Longrightarrow \alpha_8 = \mathrm{E}(\mathrm{r}_8) = 5$$

Le coefficient de renversement d'ordre 8 est le suivant:

$$\alpha_8 = 5$$

Ceci signifie que si une suite générée de longueur 8 possède un nombre des entiers impairs inferieurs ou égal a 5 donc leur coefficient transformationnel global principal est strictement inferieur à 1.

| | | | | | | | | | |
|---|---|---|---|---|---|---|---|---|---|
| 71 | 107 | 161 | 242 | 121 | 182 | 91 | 137 | 206 | 103 |

**Définition 3.8**

Prenons le cas des suites générées de Collatz de longueurs égales à 2n, On définit l'incrément d'ordre 2n comme la différence entre $\alpha_{2\mathrm{n}}$ et n comme suit :

$$\mathrm{e} = \alpha_{2\mathrm{n}} - \mathrm{n} \qquad \text{(e 3.17)}$$

**Lemme 3.1**

Pour des valeurs de n supérieures ou égales à 100 ,un encadrement de l'incrément d'ordre n est comme suit:

$$0.27n > e > 0.25n \qquad \text{(e 3.18)}$$

**Démonstration**

On cherche à déterminer un encadrement de e en fonction de n. Cet encadrement sera utilisé dans la démonstration du théorème

On sait que:

$$\frac{\mathrm{Ln}(2)}{\mathrm{Ln}(3)} = 0.6309\ldots$$

Ceci nous permet d'écrire :

$$2n\,\frac{\mathrm{Ln}(2)}{\mathrm{Ln}(3)} = 1.2618n$$

On peut déduire l'inégalité suivante pour tout entier naturel non nul n:

$$1.27n > 2n\frac{\mathrm{Ln}(2)}{\mathrm{Ln}(3)} > 1.25n$$

On cherche à déterminé un encadrement de l'incrément e pour des valeurs de n suffisamment grandes qui sera utilisé pour la démonstration lorsque n tend vers l'infini donc on peut admettre que lorsque n suffisamment grand (il suffit que $n \geq 100$ par exemple).

$$1.27n > E\left(2n\,\frac{\mathrm{Ln}(2)}{\mathrm{Ln}(3)}\right) > 1.25n$$

Ce qui nous permet d'écrire :

$$1.27n - n > E\left(2n\,\frac{\mathrm{Ln}(2)}{\mathrm{Ln}(3)}\right) - n > 1.25n - n$$

Comme on a :

$$e = E\left(2n\,\frac{\mathrm{Ln}(2)}{\mathrm{Ln}(3)}\right) - n$$

Alors un encadrement de l'incrément e est comme suit:

$$0.27n > e > 0.25n$$

Pour des valeurs de n supérieures ou égales à 100.

**Corollaire 3.1**

Pour tout entier naturel non nul n, les suites générées de type $B_n$ et de longueur n admettant des nombres d'imparité inférieurs ou égaux au coefficient $\alpha_n$ par contre les suites générées de type $A_n$ et de longueur n possèdent un nombre d'imparité strictement supérieur au coefficient d'inversement $\alpha_n$.

## 3.6 Les distributions structurelles linéaires relatives aux suites générées de Collatz

**Définition 3.9**

On considère la suite générée de Collatz de longueur n et de générateur P ci-dessous:

$$\tilde{S}^r(P, n) = (T_1^S(P), T_2^S(P), \dots, T_n^S(P))$$

Par définition, la distribution structurelle binaire linéaire relative à la suite $\tilde{S}^r(P, n)$ est la suite obtenue par remplacement de chaque terme $T_k^S(P)$ de la suite $\tilde{S}^r(P, n)$ par son indicateur de parité $\mathbf{i}_k(P)$ pour tout entier naturel k allant de 1 à n. La distribution structurelle binaire relative à la suite $\tilde{S}^r(P, n)$ est notée comme suit:

$$\tilde{L}^b(P, n) = \left(\mathbf{i}_1(P), \mathbf{i}_2(P), \dots, \mathbf{i}_n(P)\right) \quad \text{(e 3.19)}$$

La suite structurelle binaire d'une suite générée de longueur n correspond donc à un arrangement avec répétition bien déterminé de longueur n de deux entiers 0 et 1. Cet arrangement traduit la distribution des entiers pairs et des entiers impairs dans la suite considérée de Syracuse. Cette distribution structurelle binaire est appelée aussi la conversion structurelle binaire de la suite $\tilde{S}^r(P, n)$.

La figure suivante montre le principe de la conversion structurelle binaire d'une suite générée en un arrangement avec répétition de deux entiers 0 et 1 qui traduit l'aspect structurel de la suite considérée.

| | | | | | | | | | |
|---|---|---|---|---|---|---|---|---|---|
| **P** | $T_1^S(P)$ | $T_2^S(P)$ | $T_3^S(P)$ | | | | | $T_{n-1}^S(P)$ | $T_n^S(P)$ |

⇩

| | | | | | | | | | |
|---|---|---|---|---|---|---|---|---|---|
| $\mathbf{i}_0(P)$ | $\mathbf{i}_1(P)$ | $\mathbf{i}_2(P)$ | $\mathbf{i}_3(P)$ | | | | | | $\mathbf{i}_n(P)$ |

*Figure 6:* conversion d'une suite générée de Syracuse en une suite structurelle binaire

**Exemple 3.4**

Considérons la suite $\tilde{S}^r(27,9)$ représentée comme suit :

| | | | | | | | | | |
|---|---|---|---|---|---|---|---|---|---|
| 27 | 41 | 62 | 31 | 47 | 71 | 107 | 161 | 242 | 121 |

La distribution structurelle binaire linéaire (vecteur de parité) qui correspond à la suite $\tilde{S}^r(27,9)$ est la suivante :

| | | | | | | | | | |
|---|---|---|---|---|---|---|---|---|---|
| 1 | 1 | 0 | 1 | 1 | 1 | 1 | 1 | 0 | 1 |

**Définition 3.10**

Les distributions (ou les arrangements) structurelles chromatiques (de deux couleurs) relative à des suites générées données sont équivalentes aux distributions structurelles binaires non chromatique au lieu d'utiliser les indicateurs de parité 0 et 1 on utilise les indicateurs des parités chromatiques selon les règles suivantes :

- Une case contient un entier impair sera colorée en bleu.
- Une case contient un entier pair sera colorée en une autre couleur (jaune ou autre)

Par conséquent les deux conversions structurelles chromatique et non chromatique reposent sur le même principe et traduisent le même phénomène lié à la distribution des entiers pairs et des entiers impairs dans les suites générées de Syracuse.

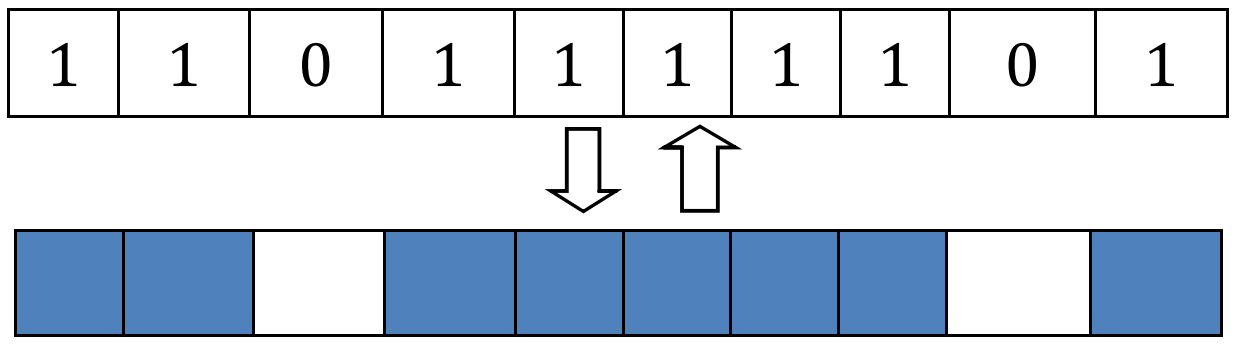

*Figure 7:* Equivalence entre les deux distributions linéaires structurelles chromatique et binaire d'une même suite générée

La ligne obtenue constitue un arrangement avec répétition de deux couleurs. Cette ligne représente la distribution structurelle chromatique linéaire relative à la suite considérée. La distribution structurelle chromatique relative à une suite $\tilde{S}^r(P,n)$ est notée $\tilde{L}^c(P,n)$.

**Exemple 3.5**

On reprend l'exemple de la suite $\tilde{S}^r(27,9)$ , sa distribution chromatique structurelle est la suivante :

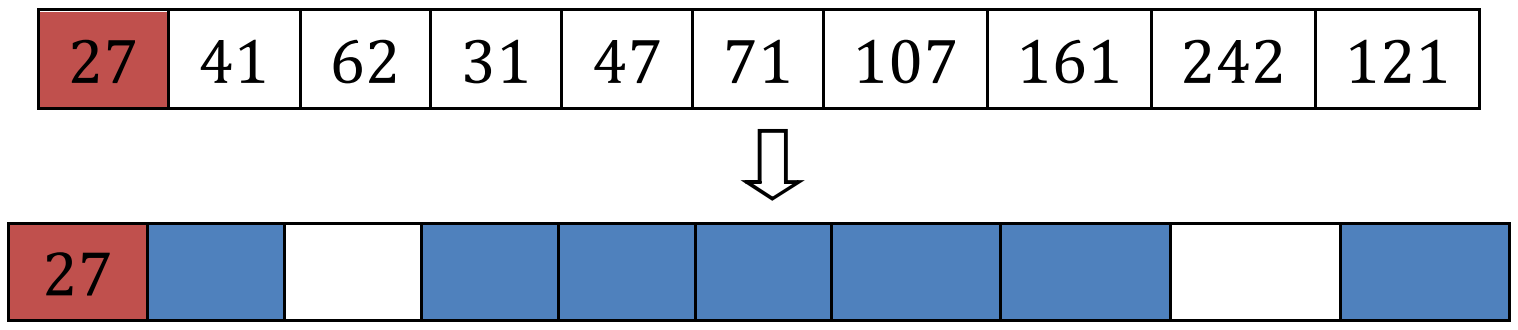

*Figure 8:* Distribution structurelle chromatique de la suite générée $\tilde{S}^r(27,9)$

**Notations 3.4**

A fin de faciliter la compréhension et l'étude de la relation entre les arrangements binaires les suites générées de Collatz, on désigne par :

$\tilde{L}^{s}(P,n)$ (Avec petite $s$ en exposant) l'une de deux distribution structurelles binaire ou bien chromatique

-$\mathbb{V}^b(n)$ l'ensemble constitué par tous les arrangements binaire avec répétition de longueur n de deux éléments 0 et 1. Un arrangement quelconque appartenant à cet ensemble est noté $v_i(n)$.

$$\mathbb{V}^b(n) = \{v_1(n), v_2(n), \dots, v_{2^n}(n)\} \tag{3.20}$$

-L'ensemble constitué par tous les arrangements avec répétition de longueur n de deux couleurs (blanc et bleu) est noté $\mathbb{W}^c(n)$. Ces arrangements sont représentés sous formes des lignes à n cases, certaines cases sont colorées en bleu alors que les autres sont colorées en blanc dans un ordre quelconque. Un élément quelconque de cet ensemble est noté $w_i(n)$.

$$\mathbb{W}^c(n) = \{w_1(n), w_2(n), \dots, w_{2^n}(n)\} \tag{3.21}$$

Puisque le nombre des arrangements possibles avec répétition de taille n de deux éléments égal à $2^n$ donc chaque ensemble $\mathbb{V}^b(n)$ et $\mathbb{W}^c(n)$ contient exactement $2^n$ arrangements avec répétition ce qui nous permet d'écrire :

$$\mathrm{card}(\mathbb{V}^b(n)) = \mathrm{card}(\mathbb{W}^c(n)) = 2^n \tag{3.22}$$

-$\mathbb{S}(n)$ est l'ensemble des suites générées de Collatz de longueur n

**Définition 3.11**

On définie l'application de la conversion structurelle comme suit :

$$T^{conv}: \mathbb{S}(n) \to \mathbb{V}^b(n)$$
$$\tilde{S}^r(P,n) \to V_i(n) \tag{3.23}$$

Comme $\tilde{L}^b(P,n)$ correspond à un arrangement avec répétition de longueur n donc il existe un entier naturel non nul $i \leq 2^n$ tel que :

$$\tilde{L}^b(P,n) = V_i(n) \qquad (3.24)$$

**Définition 3.12**

Deux ou plusieurs suites générées de même longueur n sont dites n-isoformes ou n-isochromatiques si elles possèdent la même conversion (ou distribution) structurelle linéaire. Autrement si on désigne par $\tilde{S}^r(P_1,n)$ et $\tilde{S}^r(P_2,n)$ deux suites générées de Collatz de même longueur n, on dit que ces deux suites sont n-isochromatique si elles vérifient la condition suivante:

$$\tilde{L}^b(P_1,n) = \tilde{L}^b(P_2,n) \qquad (e\ 3.25)$$

Cette condition est équivalente à:

$$\tilde{L}^c(P_1,n) = \tilde{L}^c(P_2,n)) \qquad (e\ 3.26)$$

**Exemple 3.6**

Prenons le cas de deux suites de Collatz $\tilde{S}^r(11,7)$ et $\tilde{S}^r(139,7)$ suivantes :

| | | | | | | |
|---|---|---|---|---|---|---|
| 11 | 17 | 26 | 13 | 20 | 10 | 5 |
| 139 | 209 | 314 | 157 | 236 | 118 | 59 |

Elles s'agissent de deux suites isoformes, elles possèdent la même distribution chromatique linéaire suivante :

| | | | | | | |
|---|---|---|---|---|---|---|
| | | | | | | |

Evidement elles possèdent aussi la même distribution binaire suivante :

| | | | | | | |
|---|---|---|---|---|---|---|
| 1 | 1 | 0 | 1 | 0 | 0 | 1 |

On peut écrire dans ce cas :

$$\tilde{L}^{\mathscr{s}}(11,n) = \tilde{L}^{\mathscr{s}}(139,n)$$

**Théorème 3.2**

Soient P, n et m trois entiers naturels non nuls, alors pour tout entier naturel k tel que $1 \leq k \leq n+1$ on a :

$$T_k^S(P + 2^{n+1}m) = T_k^S(P) + 2^{n+1-k}\, 3^{M_{k-1}(P)}m \qquad \text{(e 3.27)}$$

**Démonstration**

-Pour un entier naturel n bien déterminé, on montre tout abord que cette propriété est vraie pour tout entier naturel k vérifiant $1 \leq k \leq n$.

On vérifie que la propriété est vraie pour k=1. On peut écrire :

$$T_1^S(P + 2^{n+1}m) = T_1^S(P + 2x2^n m)$$
$$= T_1^S(P) + 3^{\mathbf{i}_0(P)}\, 2^n m$$

Comme $\mathbf{i}_0(P) = M_0(P)$

$$T_1^S(P + 2^{n+1}m) = T_1^S(P) + 3^{M_0(P)}\, 2^n m \qquad \text{(e 3.33)}$$

Prenons un entier naturel non nul k quelconque tel que $k \leq n$, on montre que si la propriété est vraie pour k alors elle est nécessairement vraie pour (k+1)

On suppose que:

$$T_k^S(P + 2^{n+1}m) = T_k^S(P) + 2^{n+1-k}\, 3^{M_{k-1}(P)}m$$

Ce qui nous permet d'écrire:

$$T_{k+1}^S(P + 2^{n+1}m) = T_1^S\left(T_k^S(P + 2^{n+1}m)\right)$$
$$= T_1^S\left(T_k^S(P) + 2^{n+1-k}\, 3^{M_{k-1}(P)}m\right) \qquad \text{(e 3.34)}$$

En utilisant la formule (e 2.22) établie pour l'expression unificatrice de l'image d'un entier qui s'écrit sous la forme (A+2B) par $T^S$ on peut écrire alors:

$$T_1^S\left(T_k^S(P) + 2x2^{n-k}\, 3^{M_{k-1}(P)}m\right) = T_{k+1}^S(P) + 3^{\mathbf{i}_0(T_k^S(P))} 2^{n-k}\, 3^{M_{k-1}(P)}m \quad \text{(e 3.35)}$$

On utilise les propriétés de l'indicateur de parité et de nombre d'imparité suivantes :

$$\mathbf{i}_0\left(T_k^S(P)\right) = \mathbf{i}_k(P) \text{ de plus on a } \mathbf{i}_k(P) + M_{k-1}(P) = M_k(P)$$

On peut écrire:

$$T_{k+1}^S(P) + 3^{\mathbf{i}_0\left(T_k^S(P)\right)} 2^{n-k}\, 3^{M_{k-1}(P)}m = T_{k+1}^S(P) + 3^{\mathbf{i}_k(P)} 2^{n-k}\, 3^{M_{k-1}(P)}m$$
$$= T_{k+1}^S(P) + 2^{n-k}\, 3^{M_k(P)}m \qquad \text{(e 3.33)}$$

Ce qui nous permet de déduire que :

$$T_{k+1}^S(P + 2^{n+1}m) = T_{k+1}^S(P) + 2^{n-k}\, 3^{M_k(P)}m$$

On peut conclure que si la propriété est vraie pour un entier k quelconque inférieur ou égal à n alors nécessairement la propriété est vraie pour (k+1)

Comme la propriété est vraie pour k=1 alors elle est aussi vraie pour k=2 puis pour k=3… donc par itération successive jusqu' à k=n, on peut admettre que pour tout entier k tel que $1 \leq k \leq n$

$$T_k^S(P + 2^{n+1}m) = T_k^S(P) + 2^{n-k+1}\, 3^{M_{k-1}(P)}m$$

-Pour k=n+1

On sait que:

$$T_n^S(P + 2^{n+1}m) = T_n^S(P) + 2\, 3^{M_{n-1}(P)}m \qquad \text{(e 3.33)}$$

En utilisant la formule on peut écrire:

$$T_{n+1}^S(P + 2^{n+1}m) = T_1^S\left(T_n^S(P) + 2\, 3^{M_{n-1}(P)}m\right)$$

$$= T_{n+1}^S(P) + 3^{\mathbf{i}_0(T_n^S(P))}\, 3^{M_{n-1}(P)}m \qquad \text{(e 3.33)}$$

On sait que:

$$\mathbf{i}_0\left(T_n^S(P)\right) = \mathbf{i}_n(P) \text{ et } \mathbf{i}_n(P) + M_{n-1}(P) = M_n(P)$$

Ce qui nous permet d'écrire:

$$T_{n+1}^S(P) + 3^{\mathbf{i}_0(T_n^S(P))}\, 3^{M_{n-1}(P)}m = T_{n+1}^S(P) + 3^{M_n(P)}m \qquad \text{(e 3.35)}$$

On peut conclure que la propriété est vraie pour k= n+1 donc pour tout entier naturel k tel que$1 \leq k \leq n+1$, on a:

$$T_k^S(P + 2^{n+1}m) = T_k^S(P) + 2^{n+1-k}\, 3^{M_{k-1}(P)}m$$

**Corollaire 3.2**

Soients P, n et m trois entiers naturels non nuls, alors Pour tout entier k tel que $1 \leq k \leq n$:

$$\mathbf{i}_k(P + 2^{n+1}m) = \mathbf{i}_k(P) \qquad \text{(e 3.29)}$$

Autrement les deux entiers $T_k^S(P + 2^{n+1}m)$ et $T_k^S(P)$ on la même parité pour tout entier k allant de 1 à n.

C'est une conséquence du théorème précédente, en effet comme le terme $2^{n+1-k}\, 3^{M_{k-1}(P)}m$ est un entier pair non nul pour tout entier naturel k inférieur ou égale à n donc on peut conclure que les deux entiers $T_k^S(P + 2^{n+1}m)$ et $T_k^S(P)$ ont la même parité pour tout entier k allant de 1 à n.

**Théorème 3.3**

Soient $P_1$ et $P_2$ deux entiers naturels non nuls. Si $P_1$ et $P_2$ ont la même parité alors les deux suites générées $\tilde{S}^r(P_1, n)$et $\tilde{S}^r(P_2, n)$ sont isochromatiques si et seulement si ils existent deux entiers naturels non nul m et n tel que :

$$P_2 = P_1 + 2^{n+1}m \qquad \text{(e 3.30)}$$

Autrement:

$$\begin{cases} \mathbf{i}_0(P_2) = \mathbf{i}_0(P_1) \\ \tilde{L}^S(P_1, n) = \tilde{L}^S(P_1, n) \end{cases} \Leftrightarrow P_2 = P_1 + 2^{n+1}m \qquad \text{(e 3.31)}$$

**Démonstration**

On suppose que $P_2 = P_1 + 2^{n+1}m$ avec m et n deux entiers naturels non nuls.

D'après le théorème (e 2.22) et le corollaire (e 2.22) pour tout entier naturel k tel que $1 \leq k \leq n$, on a :

$$T_k^S(P + 2^{n+1}m) = T_k^S(P) + 2^{n+1-k}\, 3^{M_{k-1}(P)}m$$

Comme $2^{n+1-k}\, 3^{M_{k-1}(P)}m$ est un entier pair pour tout entier naturel k vérifiant $1 \leq k \leq n$ donc on peut déduire que $T_k^S(P + 2^{n+1}m)$ et $T_k^S(P)$ ont la même parité ceci se traduit par la relation suivante :

$$\mathbf{i}_k(P + 2^{n+1}m) = \mathbf{i}_k(P)$$

Pour tout entier naturel k qui vérifie $1 \leq k \leq n$

Ceci nous permet de conclure que les deux suites $\tilde{S}^r(P, n)$ et $\tilde{S}^r(P + 2^{n+1}m, n)$ ont la même distribution structurelle linéaire de longueur n donc elles sont isochromatiques.

Réciproquement on suppose que les deux suites $\tilde{S}^r(P, n)$ et $\tilde{S}^r(P + 2^{n+1}m, n)$ sont isochromatiques c'est à dire que les deux suites $\tilde{S}^r(P_1, n)$ et $\tilde{S}^r(P_2, n)$ ont la même distribution structurelle linéaire ce qui nous permet de conclure que pour tout entier k tel que $1 \leq k \leq n$:

$T_k^S(P_1)$ et $T_k^S(P_2)$ possèdent la même parité (ou le même indicateur de parité) et elles possèdent le même nombre d'imparité ceci se traduit par les deux relations suivante:

$$\mathbf{i}_k(P_2) = \mathbf{i}_k(P_1) \qquad \text{(e 3.33)}$$

$$M_n(P_1) = M_n(P_2) \qquad \text{(e 3.35)}$$

On montre tout à bord que $\varphi_{n+1}(P_2) = \varphi_{n+1}(P_1)$

On sait que :

$$\begin{cases} \varphi_1(P_1) = \frac{1}{2}\mathbf{i}_0(P_1) \\ \varphi_1(P_2) = \frac{1}{2}\mathbf{i}_0(P_2) \end{cases}$$

Comme $P_1$ et $P_2$ ont la même parité, on peut déduire que :

$$\varphi_1(P_1) = \varphi_1(P_2) \qquad \text{(e 3.44)}$$

La propriété est vraie pour k=1.

En utilisant la relation de récurrence suivante :

$$\varphi_k(P) = \frac{3^{\mathbf{i}_{k-1}(P)}}{2}\varphi_{(k-1)}(P) + \frac{1}{2}\mathbf{i}_{k-1}(P)$$

Donc pour tout entier naturel non nul k on peut déduire l'implication ci dessous :

$$\begin{cases} \varphi_{k-1}(P_1) = \varphi_{k-1}(P_2) \\ \mathbf{i}_{k-1}(P_1) = \mathbf{i}_{k-1}(P_2) \end{cases} \Rightarrow \varphi_k(P_1) = \varphi_k(P_2) \qquad \text{(e 3.33)}$$

Ce qui implique:

$$\varphi_2(P_1) = \varphi_2(P_2) \qquad \text{(e 4.4)}$$

De la même manière on peut démontrer que $\varphi_3(P_1) = \varphi_3(P_2)$ et $\varphi_4(P_1) = \varphi_4(P_2)$ de proche en proche on peut suivre le calcul jusqu'à $k = n + 1$

$$\varphi_n(P_1) = \varphi_n(P_2) \qquad \text{(e 5.55)}$$

Comme on a $\mathbf{i}_n(P_1) = \mathbf{i}_n(P_2)$, puisque on a:

$$\varphi_{n+1}(P_1) = \frac{3^{\mathbf{i}_n(P_1)}}{2}\varphi_n(P_1) + \frac{1}{2}\mathbf{i}_n(P_1)$$

$$\varphi_{n+1}(P_2) = \frac{3^{\mathbf{i}_n(P_2)}}{2}\varphi_n(P_2) + \frac{1}{2}\mathbf{i}_n(P_2)$$

On peut conclure que :

$$\varphi_{n+1}(P_1) = \varphi_{n+1}(P_2) \qquad \text{(e 3.33)}$$

On sait que d'après:

$$T^S_{n+1}(P_2) = \frac{3^{M_n(P_2)}}{2^{n+1}}P_2 + \varphi_{n+1}(P_2)$$

$$T^S_{n+1}(P_1) = \frac{3^{M_n(P_1)}}{2^{n+1}}P_1 + \varphi_{n+1}(P_1)$$

Comme on a $\varphi_{n+1}(P_2) = \varphi_{n+1}(P_1)$ et $M_n(P_1) = M_n(P_2)$ on peut déduire que:

$$\frac{3^{M_n(P_1)}}{2^{n+1}}(P_2 - P_1) = T^S_{n+1}(P_2) - T^S_{n+1}(P_1) \qquad \text{(e 4.44)}$$

D'une part $(P_2 - P_1)$ et $(T^S_{n+1}(P_2) - T^S_{n+1}(P_1))$ sont deux entiers naturels non nuls et d'autre part $3^{M_n(P_1)}$ et $2^{n+1}$ sont premiers entre eux donc $(P_2 - P_1)$ doit être un multiple de $2^{n+1}$ ce qui nous permet de conclure que $(P_2 - P_1)$ peut s'écrire sous la forme :

$$P_2 - P_1 = 2^{n+1}\, m \qquad \text{(e 4.44)}$$

Avec m un entier naturel non nul

**Exemple 3.7**

Prenons le cas des entiers qui s'écrits sous la forme $7 + 2^6 m$

$N_1 = 7$ ; $N_2 = 7 + 2^6 = 71$; $N_3 = 7 + 2x2^6 = 135$ et $N_4 = 7 + 3x2^6 = 199$

Les suites (non générées) $S^r(7,5)$ , $S^r(71,5)$ , $S^r(135,5)$ et $S^r(199,5)$ sont isochromatiques

| | | | | | |
|---|---|---|---|---|---|
| 7 | 11 | 17 | 26 | 13 | 20 |
| 71 | 107 | 161 | 242 | 121 | 182 |
| 135 | 203 | 305 | 458 | 229 | 344 |
| 199 | 299 | 449 | 674 | 337 | 506 |

*Figure 9:* Quelques suites isoformes

La différence entre deux générateurs de deux suites n-isoformes et consécutifs est constance, elle est égale à $2^n$ avec n la longueur de ces suites.

**Lemme 3.2**

Soient $\tilde{S}^r(P_1, n)$ et $\tilde{S}^r(P_2, n)$ deux suites générées finies de Syracuse de même taille n, si ces deux suites sont isochromatiques (ou isoformes) donc leurs coefficients transformationnels globaux principaux sont égaux et aussi leurs coefficients transformationnels globaux secondaires sont égaux. Autrement s'il existe un entier naturel non nul m tel que $P_2 = (P_1 + 2^{n+1} m)$ donc les deux relations suivantes sont vérifiées:

$$\begin{cases} \dfrac{3^{\tilde{M}_n(P_1)}}{2^n} = \dfrac{3^{\tilde{M}_n(P_2)}}{2^n} \\ \tilde{\varphi}_n(P_1) = \tilde{\varphi}_n(P_2) \end{cases} \qquad \text{(e 3.33)}$$

# 4 Les super-suites

Les distributions génératrices de Collatz sont des suites arithmétiques des rasions arithmétiques bien déterminées permettant de générer des suites finies (ou infinies) de Collatz possédant des propriétés communes entre eux. On peut distinguer deux types fondamentaux des suites génératrices de Collatz

-Les super-suites qui sont les suites génératrices des chromologues de Collatz.

-Les suites génératrices parfaites qui sont les suites génératrices des chromoformes parfaits de Collatz.

**Définition 4.1**

La relation établit entre deux générateurs de deux suites générées isochromatiques de longueur n, nous conduit à la notion des super-suites qui sont appelées aussi les super-génératrices. Généralement, elles sont des suites arithmétiques permettant de créer des suites générées de Collatz parfaitement isochromatiques. Plus précisément, une super-suite est une suite arithmétique de raison arithmétique de la forme $2^{n+1}$ avec n un entier naturel non nul et de premier terme un entier naturel non nul P tel que ($P \leq 2^{n+1}$). Cette super-suite est notée comme suit :

$$X^S(P, n) = (P, P + 2^{n+1}, P + 2x2^{n+1}, \dots, P + jx2^{n+1}, \dots) \quad \text{(e 4.1)}$$

-P et n sont appelés les paramètres constructifs de la super-suite $X^S(P, n)$.

-n est appelé ordre constructif de la super-suite considérée.

**Remarque 4.1**

Le terme général de la super-suite peut s'écrire sous la forme suivante :

$$P_j = P + (j - 1)x2^{n+1} \quad \text{(e 4.2)}$$

Avec j un entier naturel non nul quelconque.

Une super-suite est entièrement déterminée en fonction de son premier terme P et de son ordre constructif n. Le plus souvent on caractérise une super-suite par leur ordre constructif au lieu de leur raison arithmétique.

Noter que la super-suite est de longueur infini, elle renferme tous les termes de la suite arithmétique de premier terme P (avec $P \leq 2^{n+1}$ ) et de raison arithmétique $2^{n+1}$.

**Définition 4.2**

On sait que les premiers termes des super-suites de même ordre constructif n sont des entiers naturels non nuls inférieurs ou égaux à $2^{n+1}$ donc on peut conclure qu'on peut définir exactement $2^{n+1}$ super-suites de même ordre constructif n dont les premiers termes sont les suivants :

$$P_1 = 1, P_2 = 2, P_3 = 3, P_4 = 4, \dots, P_{2^{n+1}-1} = 2^{n+1} - 1, P_{2^{n+1}} = 2^{n+1} \qquad \text{(e 4.3)}$$

Les super-suites qui correspondent aux ces premiers termes constitue un super-système qu'on le note comme suit:

$$Sy^r(n) = [X^S(1,n), X^S(2,n), X^S(3,n), X^S(4,n), \dots, X^S(2^{n+1} - 1, n), X^S(2^{n+1}, n)] \quad \text{(e 4.4)}$$

**Définition 4.3**

D'une règle générale, on fait séparer les super-suites dont les premiers termes sont des entiers pairs et les super-suites dont les premiers termes sont des entiers impairs. Le super-système déjà défini est subdivisé en deux sous super-système :

-Le premier super-système contient toutes les super-suites dont les termes sont des entiers naturels impairs:

$$Sy^r(1,n) = [X^S(1,n), X^S(3,n), X^S(5,n), \dots, X^S(2^{n+1} - 1, n)] \qquad \text{(e 4.5)}$$

-Le deuxième est le super-système qui contient toutes les super-suites dont les termes sont des entiers pairs:

$$Sy^r(2,n) = [X^S(2,n), X^S(4,n), X^S(6,n), \dots, X^S(2^{n+1}, n)] \qquad \text{(e 4.6)}$$

Chaque super-système contient le même nombre des super-suites qui égal à $2^n$ super-suites.

**Remarque 4.2**

Les deux super-systèmes présentent des propriétés identiques, généralement chaque super- système est représenté séparément à l'autre dans une matrice qui ne contient que l'un de deux.

## 4.1 Représentation matricielle d'un super-système :

Une super-matrice est la représentation matricielle des toutes les super-suites d'un même super-système donc pour un super-système d'ordre constructifs n, la super-matrice correspond à un tableau de $2^n$ lignes et d'une infinité des colonnes. Chaque ligne du tableau contient tous les termes d'une même super-suite du super-système considéré. Comme on a deux super-système donc on peut distinguer deux super-matrices, la première est notée $\mathbb{X}^S(1,n)$ relative au super-système $Sy^r(1,n)$ et la deuxième est notée $\mathbb{X}^S(2,n)$ relative au super-système $Sy^r(2,n)$. Le tableau suivant représente la super-matrice $\mathbb{X}^S(1,n)$ du super-système $Sy^r(1,n)$.

*Tableau 1:* la super-matrice $\mathbb{X}^S(1,n)$

| | | | | | | |
|---|---|---|---|---|---|---|
| $X^S(1,n)$ | 1 | $1+2^{n+1}$ | | | $1+2^{n+1}(i-1)$ | |
| $X^S(3,n)$ | 3 | $3+2^{n+1}$ | | | $3+2^{n+1}(i-1)$ | |
| $X^S(5,n)$ | 5 | $5+2^{n+1}$ | | | $5+2^{n+1}(i-1)$ | |
| | | | | | | |
| | | | | | | |
| $X^S(2j-1,n)$ | $2j-1$ | $2j-1+2^{n+1}$ | | | $2j-1+2^{n+1}(i-1)$ | |
| | | | | | | |
| $X^S(2^{n+1}-1,n)$ | $2^{n+1}-1$ | $2^{n+2}-1$ | | | $2^{n+1}(i-1)-1$ | |

Noter que la première colonne séparée du reste du tableau par un trait épais ne fait pas partie de la super-matrice. Elle est appelée colonne indicatrice, elle permet de faciliter l'interprétation et l'étude des matrices étudiée.

De même on peut définir la super-matrice $\mathbb{X}^S(2,n)$ relative au super-système $Sy^r(2,n)$ comme suit :

*Tableau 2:* La super-matrice $\mathbb{X}^S(2,n)$

| | | | | | | | |
|---|---|---|---|---|---|---|---|
| $X^S(2,n)$ | 2 | $2+2^{n+1}$ | | | | $2+2^{n+1}(i-1)$ | |
| $X^S(4,n)$ | 4 | $4+2^{n+1}$ | | | | $4+2^{n+1}(i-1)$ | |
| $X^S(6,n)$ | 6 | $6+2^{n+1}$ | | | | $6+2^{n+1}(i-1)$ | |
| | | | | | | | |
| | | | | | | | |
| $X^S(2j,n)$ | $2j$ | $2j+2^{n+1}$ | | | | $2j+2^{n+1}(i-1)$ | |
| | | | | | | | |
| $X^S(2^{n+1},n)$ | $2^{n+1}$ | $2^{n+2}$ | | | | $(i+1)x2^{n+1}$ | |

## 4.2 Représentation matricielle à deux indices d'un super-système

Pour des raisons de simplification, on fait recours à la représentation matricielle de la super-matrice avec des termes à deux indices. Cette représentation consiste à remplacer chaque terme de la super-matrice par une notation à deux indices.

Chaque terme du tableau s'écrit sous la forme d'une lettre P à deux indices $P_{i,j}$ tel que :

Pour tout entier non nul i et pour tout entier j vérifiant $1 \leq j \leq 2^n$, l'entier noté $P_{i,j}$ désigne l'élément de la colonne i et de la ligne j de la super-matrice.

La super-suite du rang j de premier terme $P_{1,j}$ s'écrit comme suit :

$$X^S(P_{1,j}, n) = (P_{1,j}, P_{2,j}, P_{3,j}, \ldots., P_{i,j}, \ldots.) \quad (4.7)$$

Tel que pour tout entier naturel i, le terme du rang i de la super-suite s'écrit comme suit :

$$P_{i,j} = P_{1,j} + (i-1)x2^{n+1} \quad (4.8)$$

Le tableau suivant représente une super-matrice à deux indices. La double indexation nous facilite l'étude des super-matrices et aussi des toutes autres types des matrices qu'on va les voir.

*Tableau 3:* La super-matrice $\mathbb{X}^S(P_{1,1}, n)$

| | | | | | | | |
|---|---|---|---|---|---|---|---|
| $X^S(P_{1,1}, n)$ | $P_{1,1}$ | $P_{2,1}$ | $P_{3,1}$ | | | $P_{i,1}$ | |
| $X^S(P_{1,2}, n)$ | $P_{1,2}$ | $P_{2,2}$ | $P_{3,2}$ | | | $P_{i,2}$ | |
| $X^S(P_{1,3}, n)$ | $P_{1,3}$ | $P_{2,3}$ | $P_{3,3}$ | | | $P_{i,3}$ | |
| | | | | | | | |
| $X^S(P_{1,j}, n)$ | $P_{1,j}$ | $P_{2,j}$ | $P_{3,j}$ | | | $P_{i,j}$ | |
| | | | | | | | |
| | | | | | | | |
| $X^S(P_{1,2^n}, n)$ | $P_{1,2^n}$ | $P_{2,2^n}$ | $P_{3,2^n}$ | | | $P_{i,2^n}$ | |

Chaque super-matrice contient $2^n$ super-suites deux à deux distinctes..

**Remarques 4.3**

-Une distribution entière verticale désigne toute suite constituée par les termes contenus dans une même colonne d'un même tableau (ou d'une même matrice).

-Une distribution entière horizontale désigne une suite dont les termes sont les éléments contenus dans une même ligne d'un même tableau (ou bien d'une même matrice).

Par exemple les distributions horizontales dans la super-matrice correspondent aux différentes super-suites déjà définies alors les distributions entières verticales dans cette même matrice correspondent à des suites arithmétiques de rasions 2 qu'on l'appelle les suites génératrices parfaites qu'on va les voir dans la deuxième partie de cette section.

# 5 Les suites génératrices parfaites du Collatz

**Définition 5.1**

Une suite génératrice parfaite de Collatz désigne toutes suites arithmétiques de raison arithmétique 2, de longueur $2^n$ et de premier terme un entier naturel non nul P. Egalement, on peut distinguer deux suites génératrices fondamentales de premier terme respectivement 1 et 2 et de même longueur $2^n$. Elles sont notées comme suit:

$$\begin{cases} Y^S(1,n) = (1,3,5,\dots,2^{n+1}-1) \\ Y^S(2,n) = (2,4,6,\dots,2^{n+1}) \end{cases} \quad \text{(e 5.1)}$$

-n est un entier naturel non nul, il est appelé ordre dimensionnel de la suite génératrice considérée.

-Le couple $(1,n)$constitue le couple des paramètres constructifs de la suite parfaite $Y^S(1,n)$

**Définition 5.2**

Les deux suites génératrices peuvent être définies par leurs termes généraux comme suit :

Pour tout entier non nul j tel que $1 \leq j \leq 2^n$, le terme du rang j peut s'écrire comme suit :

$$P_j = P_1 + 2(j-1) \quad \text{(e 5.2)}$$

Avec $P_1 = 1$ pour la première suite génératrice $Y^S(1,n)$ et $P_1 = 2$ si la suite considérée correspond à la suite $Y^S(2,n)$.

Chaque suite génératrice contient alors exactement $2^n$ termes et la réunion de deux suites correspond à l'ensemble constitué par tous les entiers naturels non nuls allant de 1 à $2^{n+1}$.

**Définition 5.3**

En se basant sur la définition de la suite génératrice parfaite d'ordre dimensionnel n, on peut déduire qu'on peut définir une infinité des suites génératrice parfaites de même ordre dimensionnel et qui sont deux à deux distincts. L'ensemble des suites génératrices parfaites constitue par le système générateur parfait est noté comme suit:

$$Sy^P(P_{1,1},n) = \left[Y^S(P_{1,1},n), Y^S(P_{2,1},n), Y^S(P_{3,1},n),, \dots, Y^S(P_{i,1},n), \dots.\right] \quad (5.3)$$

Tel que pour tout entier naturel non nul $k, P_{k,1}$ désigne le premier terme de la suite génératrice parfaite du rang k dans le système générateur parfait.

**Proposition 5.1**

l'expression du premier terme de la suite génératrice parfaite du rang k dans le système générateur parfait est comme suit :

$$P_{k,1} = P_{1,1} + 2^{n+1}(k-1) \text{ avec } k \geq 1 \quad \text{(e 5.4)}$$

$P_{1,1}$ est le premier terme de la première suite génératrice parfaite. Ce terme ne peut prendre que deux valeurs 1 ou bien 2

**Démonstration**

Si on considère une suite arithmétique de raison arithmétique 2 et de longueur infinie qui renferme par exemple tous les entiers naturels impairs et qu'on la note comme suit :

$$\mathbb{I} = (1,3,5,7, \ldots, 2j + 1, ..)$$

On peut subdiviser cette suite en une infinité des suites arithmétiques finies de même longueur $2^n$ et de même raison arithmétique 2 de plus elles sont deux à deux distinctes. L'ensemble des suites génératrices parfaites de même ordre dimensionnel constitue le système générateur parfait. Il est noté comme suit :

$$(1,3,5,7, \ldots) \rightarrow ((\underbrace{1,3, .., 2^{n+1} - 1}_{Y^S(1,n)}), (\underbrace{2^{n+1} + 1, \ldots, 2^{n+2} - 1}_{Y^S(2^{n+1}+1,n)}), (\underbrace{2^{n+2} + 1, \ldots, 3x2^{n+1} - 1}_{Y^S(2^{n+2}+1,n)}), \ldots, \ldots)$$

Toutes les suites sont des suites arithmétiques de rasions 2 et de même dimension $2^n$.

Pour que ces suites soient deux à deux distinctes, il faut que le preterme de la suite du rang (k-1) représente le premier terme de la suite du rang k.

Par exemple si on considère la suite génératrice parfaite $Y^S(1, n)$. Le preterme de cette suite c'est à dire le terme qui suit son dernier terme est définie comme suit :

$$(2^{n+1} - 1) + 2 = 2^{n+1} + 1$$

Ce dernier entier représente le premier terme de la deuxième suite génératrice parfaite notée $Y^S(2^{n+1} + 1, n)$ tel que :

$$Y^S(2^{n+1} + 1, n) = (2^{n+1} + 1, 2^{n+1} + 3, 2^{n+1} + 5, \ldots, 2^{n+2} - 1)$$

On se basant sur cette relation de récurrence on peut déduire l'expression du premier terme de la suite génératrice parfaite du rang k dans le système générateur parfait est comme suit :

$$P_{k,1} = P_{1,1} + 2^{n+1}(k - 1) \text{ avec } k \geq 1$$

$P_{1,1}$ est le premier terme de la première suite génératrice parfaite. Ce terme ne peut prendre que deux valeurs 1 ou bien 2 donc, on peut distinguer deux systèmes générateurs parfaits :

**Définition 5.4**

Le premier système parfait est constitué par toutes les suites génératrices parfaites dont les premiers termes sont des entiers impairs, il est défini comme suit :

$$Sy^P(1,n) = [Y^S(1,n), Y^S(1+2^{n+1},n), Y^S(1+2x2^{n+1},n),, \dots, Y^S(1+2^{n+1}k,n), \dots.] \quad (5.5)$$

**Définition 5.5**

Le deuxième système générateur parfait renferme toutes les suites génératrices parfaites dont les premiers termes sont des entiers pairs, on le note comme suit:

$$Sy^P(2,n) = [Y^S(2,n), Y^S(2+2^{n+1},n), Y^S(2+2x2^{n+1},n),, \dots, Y^S(2+2^{n+1}k,n), \dots.] \quad (e\ 5.6)$$

Les suites d'un même système parfait présentent les propriétés suivantes :

## 5.1 La matrice génératrice parfaite $\mathbb{Y}^S(P_{1,1}, \mathbf{n})$

La représentation du système générateur parfait sous forme matricielle correspond à un tableau constitué de $2^n$ colonnes et d'une infinité des lignes. Chaque ligne correspond à une suite génératrice parfaite bien déterminée. Cette matrice est appelée la matrice génératrice parfaite d'ordre dimensionnel n et de générateur fondamental $P_{1,1}$ elle est notée $\mathbb{Y}^S(P_{1,1}, n)$.

*Tableau 4 :*La matrice génératrice parfaite $\mathbb{Y}^S(P_{1,1}, n)$:

| | $X^S(P_{1,1},n)$ | $X^S(P_{1,2},n)$ | | $X^S(P_{1,j},n)$ | | $X^S(P_{1,2^n},n)$ |
|---|---|---|---|---|---|---|
| $Y^S(P_{1,1},n)$ | $P_{1,1}$ | $P_{1,2}$ | | $P_{1,j}$ | | $P_{1,2^n}$ |
| $Y^S(P_{2,1},n)$ | $P_{2,1}$ | $P_{2,2}$ | | $P_{2,j}$ | | $P_{2,2^n}$ |
| $Y^S(P_{3,1},n)$ | $P_{3,1}$ | $P_{3,2}$ | | $P_{3,j}$ | | $P_{3,2^n}$ |
| | | | | | | |
| | | | | | | |
| $Y^S(P_{i,1},n)$ | $P_{i,1}$ | $P_{i,2}$ | | $P_{i,j}$ | | $P_{i,2^n}$ |
| | | | | | | |

Noter que la ligne vide signifie que la matrice possède une infinité des lignes. Les distributions verticales de la matrice génératrice parfaites correspondent aux super-suites génératrices des chromologues par exemple la première colonne de la matrice parfait contient tous les termes de la super-suite $X^S(P_{1,1}, n)$.

**Lemme 5.1**

Comme chaque colonne de la matrice génératrice parfaite correspond a une super-suite donnée de la super-matrice $\mathbb{X}^S(P_{1,1}, n)$ donc on peut conclure que la matrice génératrice parfaite $\mathbb{Y}^S(P_{1,1}, n)$ correspond à la transposée de la super-matrice $\mathbb{X}^S(P_{1,1}, n)$ ceci est traduit par l'expression suivante :

$$\mathbb{Y}^S(P_{1,1}, n) = \left(\mathbb{X}^S(P_{1,1}, n)\right)^{Transp} \quad (5.7)$$

-Les entiers $P_{1,1}, n$ sont appelés les paramètres constructifs de la matrice parfaite.

-n est appelé l'ordre dimensionnel de la matrice parfaite.

-$P_{1,1}$ est appelé le constructeur fondamental de la matrice parfaite.

**Remarques 5.1**

-La ligne du rang k du tableau correspond à la suite parfaite de premier terme $P_{k,1}$ qui s'écrit comme suit :

$$P_{k,1} = P_{1,1} + 2^{n+1}(k - 1) \quad (e\ 5.8)$$

-La colonne du rang j de la matrice génératrice parfaite correspond à la super-suite $X^S(P_{1,j}, n)$ tel que :

$$P_{1,j} = P_{1,1} + 2(j - 1) \quad (e\ 5.9)$$

# 6 Les chromologues du Collatz

**Définition 6.1**

Un chromologue (on dit aussi super-chromologues) de Collatz est l'ensemble constitué par une infinité des suites isochromatiques (ou isoformes) de Collatz générées par une super-suite bien déterminée.

## 6.1 Caractérisation d'un Chromologue de Collatz

Considérons la super-suite suivante :

$$X^S(P_{1,j}, n) = (P_{1,j}, P_{2,j}, P_{3,j}, \dots, P_{,ij}, \dots) \quad (6.1)$$

Tel que pour tout entier naturel non nul i, l'expression de $P_{ij}$ de est comme suit :

$$P_{i,j} = P_{1,j} + (i - 1)x2^{n+1} \quad (6.2)$$

On fait associer à chaque terme de la super-suite ci-dessus, une suite générée de Syracuse de longueur n et de premier terme $P_{i,j}$ comme suit :

$$P_{i,j} \rightarrow \tilde{S}^r(P_{i,j}, n) \quad (6.3)$$

On obtient un ensemble infini des suites générées de Syracuse de même longueur n. Toutes les suites obtenues sont isochromatiques, cet ensemble infini des suites isochromatiques issues d'une même super-suite est appelé le super-chromologue de Syracuse des paramètres constructifs $(P_{1j}, n)$ qu'on l'appelle aussi monoforme de Syracuse, il est noté comme suit :

$$H^{T}(P_{1,j}, n) = \left[\tilde{S}^{r}(P_{1,j}, n), \tilde{S}^{r}(P_{2,j}, n), \tilde{S}^{r}(P_{3,j}, n), \dots, \tilde{S}^{r}(P_{i,j}, n) \dots\right] \qquad (6.4)$$

- Les entiers $P_{1,j}$, n sont appelés les paramètres constructifs du chromologue considéré.

-$P_{1,j}$ est appelé le générateur fondamental du chromologue.

-n est appelé l'extension de $H^{T}(P_{1,j}, n)$.

-$\tilde{S}^{r}(P_{1,j}, n)$ est appelée la suite générée fondamentale du $H^{T}(P_{1,j}, n)$, elle correspond à la suite du premier rang du super-chromologue.

-La super-suite $X^{S}(P_{1,j}, n)$ est appelée la suite génératrice de ce super-chromologue.

-La raison arithmétique de la super-suite constitue la période structurelle ou chromatique du chromologue considéré. Elle est notée $\Lambda_{n}(P_{1,j})$ et elle est donnée par la relation suivante :

$$\Lambda_{n}(P_{1,j}) = 2^{(n+1)} \qquad (6.5)$$

Elle correspond à la raison arithmétique de la super-suite génératrice du chromologue considère.

**Définition 6.2**

Comme chaque suite générée de Collatz possède une distribution structurelle linéaire bien déterminée donc le chromologue $H^{T}(P_{1,j}, n)$ admet un chromologue structurelle constitué par l'ensemble des distributions structurelles linéaires relatives aux toutes les suites du chromologue fondamental. Pour tout entier naturel non nul j, on peut associer à chaque suite de chromologue sa propre distribution structurelle linéaire comme suit :

$$T^{conv}(\tilde{S}^{r}(P_{1,j}, n),) = \tilde{L}^{\mathcal{S}}(P_{1,j}, n)$$

Le chromologue structurel relatif au chromologue fondamental $H^{T}(P_{1j}, n)$ est noté comme suit :

$$H^{S}(P_{1,j}, n) = \left[\tilde{L}^{\mathcal{S}}(P_{1,j}, n), \tilde{L}^{\mathcal{S}}(P_{2,j}, n), \tilde{L}^{\mathcal{S}}(P_{3,j}, n), \dots, \tilde{L}^{\mathcal{S}}(P_{i,j}, n) \dots\right] \qquad (e\ 6.6)$$

## 6.2 La distribution structurelle linéaire caractéristique du super-chromologue :

On sait que les suites générées de Collatz de longueur n et qui sont générées par une même super-suite sont des suites isoformes, elles possèdent la même distribution structurelle linéaire cette propriété se traduit par la relation suivante :

$$\tilde{L}^{s}(P_{1,j},n)=\tilde{L}^{s}(P_{1,j}+2^{n+1},n)=L^{s}(P_{1,j}+2x2^{n+1},n)=\cdots=L^{s}(P_{1,j}+ix2^{n+1},n)=\cdots \quad (e\ 6.7)$$

On peut conclure que chaque chromologue peut être caractérisé par sa propre distribution linéaire structurelle unique qu'on la note comme suit :

$$F^{s}\left(H^{T}(P_{1,j},n)\right)=L^{s}(P_{1,j},n) \quad (e\ 6.8)$$

## 6.3 Représentations matricielles d'un chromologue Collatz :

On considère le super-chromologue $H^{T}(P_{1j},n)$ suivant :

$$H^{T}(P_{1,j},n)=[\tilde{S}^{r}(P_{1,j},n),\tilde{S}^{r}(P_{2,j},n),\tilde{S}^{r}(P_{3,j},n),\ldots,\tilde{S}^{r}(P_{k,j},n)\ldots]$$

La représentation matricielle de super-chromologue correspond à un tableau de n colonnes et d'une infinité des lignes. La ligne du rang k du tableau contient tous les termes de la suite générée de Syracuse de premier terme $P_{k,j}$ et de longueur n notée $\tilde{S}^{r}(P_{k,j},n)$. On ajoute une colonne supplémentaire contient les termes de la super-suite $X^{S}(P_{1,j},n)$ qui représente la suite génératrice de cette matrice. Cette première colonne est appelée colonne génératrice de la matrice. La matrice qui correspond à ce chromologue est notée $\mathbb{M}^{H}(P_{1,j},n)$.

**Exemple 6.1**

L'exemple ci-dessous correspond aux deux chromologues de Collatz $H^{T}(7,4)$ et $H^{T}(5,4)$

| | | | | |
|---|---|---|---|---|
| 7 | 11 | 17 | 26 | 13 |
| 39 | 59 | 89 | 134 | 67 |
| 71 | 107 | 161 | 242 | 121 |
| 103 | 155 | 233 | 350 | 175 |
| 135 | 203 | 305 | 458 | 229 |
| 167 | 251 | 377 | 566 | 283 |
| 199 | 299 | 449 | 674 | 337 |
| 231 | 347 | 521 | 782 | 391 |
| 263 | 395 | 593 | 890 | 445 |
| 295 | 443 | 665 | 998 | 499 |
| 327 | 491 | 737 | 1106 | 553 |
| 359 | 539 | 809 | 1214 | 607 |
| 391 | 587 | 881 | 1322 | 661 |
| 423 | 635 | 953 | 1430 | 715 |
| 455 | 683 | 1025 | 1538 | 769 |
| 487 | 731 | 1097 | 1646 | 823 |
| 519 | 779 | 1169 | 1754 | 877 |
| 5 | 8 | 4 | 2 | 1 |
| 37 | 56 | 28 | 14 | 7 |
| 69 | 104 | 52 | 26 | 13 |
| 101 | 152 | 76 | 38 | 19 |
| 133 | 200 | 100 | 50 | 25 |
| 165 | 248 | 124 | 62 | 31 |
| 197 | 296 | 148 | 74 | 37 |
| 229 | 344 | 172 | 86 | 43 |
| 261 | 392 | 196 | 98 | 49 |
| 293 | 440 | 220 | 110 | 55 |
| 325 | 488 | 244 | 122 | 61 |
| 357 | 536 | 268 | 134 | 67 |
| 389 | 584 | 292 | 146 | 73 |
| 421 | 632 | 316 | 158 | 79 |
| 453 | 680 | 340 | 170 | 85 |
| 485 | 728 | 364 | 182 | 91 |

*Figure 10:* Représentation matricielle de deux chromologues de Collatz $H^{T}(7,4)$ et $H^{T}(5,4)$

Noter que le nombre des lignes est infini.

**Théorème 6.1**

Comme toutes les suites d'un super-chromologue sont isochromatiques donc elles possèdent le même vecteur des effets transformationnels globaux ce qui signifie qu'un super-chromologue peut être caractérisé par un couple unique des effets transformationnels globaux : un coefficient transformationnel principal global et autre secondaire :

$$\begin{cases} \frac{3^{\widetilde{M}_n(P_{1,j})}}{2^n} = \frac{3^{\widetilde{M}_n(P_{2,j})}}{2^n} = \frac{3^{\widetilde{M}_n(P_{3,j})}}{2^n} = \cdots = \frac{3^{\widetilde{M}_n(P_{i,j})}}{2^n} = \cdots \\ \widetilde{\varphi}_n(P_{1,j}) = \widetilde{\varphi}_n(P_{2,j}) = \widetilde{\varphi}_n(P_{3,j}) = \cdots = \widetilde{\varphi}_n(P_{i,j}) = \cdots \end{cases} \quad \text{(e 6.9)}$$

Pour chaque arrangement binaire avec répétition $V_i(n)$ de longueur n, elles existent une infinité des suites générées de Collatz de longueur n qui ont toutes la même distribution structurelle linaire identique à cet arrangement. L'ensemble des ces suites isochromatiques constituent le chromologue de Collatz d'ordre dimensionnel n.

**Définition 6.3**

On désigne par $C\mathbb{H}$ l'ensemble des chromologues d'ordre dimensionnel n et On définie l'application $T^H$ de $\mathbb{V}^b(n)$ dans $C\mathbb{H}$ comme suit :

$$\begin{aligned} T^H &: \mathbb{V}^b(n) \to C\mathbb{H} \\ V_i(n) &\to \{\widetilde{S}^r(P + 2^{n+1}j, n) \ \ j \in \mathbb{N}\} \end{aligned} \quad \text{(e 6.10)}$$

On peut écrire alors :

$$V_i(n) = \widetilde{L}^b(P, n) = \widetilde{L}^b(P + 2^{n+1}, n) = \widetilde{L}^b(P + 2x2^{n+1}, n) = \cdots \quad \text{(e 6.11)}$$

## 6.4 Classification des chromologues

-Si la première suite générée du super-chromologue est de type $A_n$ alors toutes les autres suites du chromologue sont de type $A_n$.

-Si la première suite du super-chromologue est de type $B_n$ alors toutes les autres suites du chromologue sont de type $B_n$.

On conclu qu'un chromologue ne peut contenir qu'un seul type des suites de Collatz. Les chromologue sont classés selon le type des leurs suites générées qui le constituent, on peut distinguer donc deux types des chromologues:

**Définition 6.4**

Si les suites appartenant à un super-chromologue sont de type $A_n$ ,on dit qu'il s'agit d'un super-chromologue de type $A_n$ autrement si le coefficient transformationnel global caractéristique de super-chromologue considère est strictement supérieur à 1 le chromologue est de type $A_n$.

**Définition 6.5**

Si le coefficient transformationnel principal global caractéristique du super-chromologue est strictement inférieur à 1 on dit qu'il s'agit d'un super-chromatisme de type $B_n$.

**Exemple 6.2**

Reprenons l'exemple précédent, le coefficient transformationnel principal global caractéristique du super-chromologue $H^T(17,5)$ est comme suit :

$$\frac{3^2}{2^5} = \frac{9}{32} < 1$$

On peut conclure que le super-chromologue $H^T(17,5)$ est de type $B_5$.

L'exemple suivant correspond au super-chromologue $H^T(7,4)$ leur coefficient transformationnel principal global caractéristique est comme suit :

$$\frac{3^3}{2^4} = \frac{27}{16} > 1$$

Il s'agit d'un super-chromologue de type $A_4$.

**Corollaire 6.1**

Toutes les suites générées d'un super-chromologue de type $A_n$ sont croissantes c'est à dire qu'elles ont un effet global transformationnel grossissant. Le preterme de chaque suite est strictement supérieur à son premier terme :

$$T^S_{n+1}(P_{1,j}) > T^S_1(P_{1,j}) \quad \text{(e 6.12)}$$

Cette propriété est évidente puisque toutes les suites sont de type A donc elles possèdent un coefficient transformationnel global principal strictement supérieur à un ce qui signifie que les premiers termes sont toujours strictement inférieurs aux pretermes des suites du chromologue.

## 6.5 Propriétés des super-chromologues de type $B_n$

On considère un super-chromologue de type B comme suit :

$$H^T(P_{1,j},n)=[\tilde{S}^r(P_{1,j},n),\tilde{S}^r(P_{2,j},n),\tilde{S}^r(P_{3,j},n),\ldots,\tilde{S}^r(P_{k,j},n)\ldots]$$

Pour un super-chromologue de type B, on sait que toutes les suites sont de type $B_n$ elles possèdent un coefficient transformationnel principal global strictement inférieur à un.

**Lemme 6.1**

Si la première suite d'un super-chromologue de type $B_n$ est une suite abaissante c'est à dire de type $B_n^-$ alors toutes les autres suites sont abaissantes donc de type $B_n^-$ et le super-chromologue ne contient aucune suite de type $B_n^+$

**Démonstration :**

On considère le chromologue $H^T(P_{1,j},n)$ de premier générateur $P_{1,j}$. On choisi un terme quelconque $P_{i,j}$ du rang i de la super-suite génératrice du chromologue $H^T(P_{1,j},n)$.

Tenons compte des relations suivantes :

$$\begin{cases} 3^{\tilde{M}_n(P_{i,j})}=3^{\tilde{M}_n(P_{1,j})} \\ \tilde{\varphi}_n(P_{1,j})=\tilde{\varphi}_n(P_{1,j}) \end{cases} \quad \text{(e 6.13)}$$

Donc les expressions de deux pretermes relatives aux suites sont comme suit :

$$\begin{cases} T_{n+1}^S(P_{1,j})=\dfrac{3^{\tilde{M}_n(P_{1j})}}{2^n}T_1^S(P_{1,j})+\tilde{\varphi}_n(P_{1,j}) \\ T_{n+1}^S(P_{i,j})=\dfrac{3^{\tilde{M}_n(P_{1,j})}}{2^n}T_1^S(P_{i,j})+\tilde{\varphi}_n(P_{1,j}) \end{cases} \quad \text{(e 6.14)}$$

Ces équations sont équivalentes aux équations suivantes :

$$\begin{cases} \dfrac{T_{n+1}^S(P_{1,j})}{T_1^S(P_{1,j})}=\dfrac{3^{\tilde{M}_n(P_{1,j})}}{2^n}+\dfrac{\tilde{\varphi}_n(P_{1,j})}{T_1^S(P_{1,j})} \\ \dfrac{T_{n+1}^S(P_{i,j})}{T_1^S(P_{i,j})}=\dfrac{3^{\tilde{M}_n(P_{1,j})}}{2^n}+\dfrac{\tilde{\varphi}_n(P_{1,j})}{T_1^S(P_{i,j})} \end{cases}$$

Ce qui nous permet d'écrire :

$$\frac{T_{n+1}^S(P_{i,j})}{T_1^S(P_{i,j})}-\frac{T_{n+1}^S(P_{1,j})}{T_1^S(P_{1,j})}=\frac{\tilde{\varphi}_n(P_{i,j})}{T_1^S(P_{i,j})}-\frac{\tilde{\varphi}_n(P_{1,j})}{T_1^S(P_{1,j})}$$

On peut déduire que:

$$\frac{T_{n+1}^S(P_{i,j})}{T_1^S(P_{i,j})}-\frac{T_{n+1}^S(P_{1,j})}{T_1^S(P_{1,j})}=\frac{\tilde{\varphi}_n(P_{1,j})}{T_1^S(P_{i,j})}-\frac{\tilde{\varphi}_n(P_{1,j})}{T_1^S(P_{1,j})}$$

$$= \widetilde{\varphi}_n(P_{1j})\left(\frac{1}{T_1^S(P_{i,j})} - \frac{1}{T_1^S(P_{1,j})}\right) \qquad \text{(e 6.15)}$$

Comme $T_1^S(P_{i,j}) > T_1^S(P_{1,j})$ (puisque $P_{i,j} > P_{1,j}$) alors on peut conclure que:

$$\frac{\widetilde{\varphi}_n(P_{1,j})}{T_1^S(P_{i,j})} < \frac{\widetilde{\varphi}_n(P_{1,j})}{T_1^S(P_{j,1})} \Rightarrow \frac{T_{n+1}^{S1}(P_{i,j})}{T_1^S(P_{i,j})} - \frac{T_{n+1}^{S1}(P_{1,j})}{T_1^S(P_{1,j})} < 0$$

Puisque la suite $\tilde{S}^r(P_{1,j}, n)$ est de type $B_n^-$ donc elle vérifie l'inéquation suivante :

$$T_{n+1}^{S1}(P_{1,j}) \leq T_1^S(P_{1,j})$$

On déduit que:

$$\frac{T_{n+1}^{S1}(P_{i,j})}{T_1^S(P_{i,j})} < \frac{T_{n+1}^{S1}(P_{1,j})}{T_1^S(P_{1,j})} \leq 1$$

On conclut que:

$$T_{n+1}^S(P_{i,j}) < T_1^S(P_{i,j}) \qquad \text{(e 6.16)}$$

On peut déduire que :

La suite générée $\tilde{S}^r(P_{i,j}, n)$ est une suite de type $B_n^-$. Cette propriété est vraie pour n'importe quelle suite du chromologue considéré et par suite ce chromologue ne contient que des suites de type $B_n^-$.

**Lemme 6.2**

Si la première suite d'un super-chromologue $B_n$ est une suite croissante c'est-à-dire qu'elle s'agit d'une suite de type $B_n^+$, donc il existe un entier fini bien déterminée m de la super-suite génératrice du super-chromologue tel que pour tout terme $P_{i,j} > P_{i_{max},j}$, la suite $\tilde{S}^r(P_{i,j}, n)$ est une suite abaissante c'est à dire qu'elle s'agit d'une suite de type $B_n^-$.

$$\begin{cases} T_1^S(P_{i,j}) > T_{n+1}^S(P_{i,j}) & \text{si } P_{1j} \leq P_{i,j} \leq P_{i_{max},j} \\ T_1^S(P_{i,j}) \leq T_{n+1}^S(P_{i,j}) & \text{si } P_{i,j} > P_{i_{max},j} \end{cases} \qquad \text{(e 6.17)}$$

$P_{i_{max}j}$ est appelé le point de renversement du super-chromologue de type $B_n$ considéré.

**Corollaire 6.2**

Un super-chromologue de type $B_n$ contient une infinité des suites abaissantes de type $B_n^-$ et un nombre fini des suites croissantes de type $B_n^+$ qui peut être nul ou non nul.

**Démonstration du lemme 6.2**

On sait que l'expression de preterme d'une suite générée quelconque du rang i du chromologue est de la forme :

$$T_{n+1}^{S}(P_{i,j}) = \frac{3^{\widetilde{M}_n(P_{1,j})}}{2^n} T_1^S(P_{i,j}) + \widetilde{\varphi}_n(P_{1,j})$$

Ce qui équivaut à:

$$\frac{T_{n+1}^{S}(P_{i,j})}{T_1^S(P_{i,j})} = \frac{3^{\widetilde{M}_n(P_{1,j})}}{2^n} + \frac{\widetilde{\varphi}_n(P_{1,j})}{T_1^S(P_{i,j})}$$

On sait que:

$$T_1^S(P_{i,j}) = \frac{3^{\mathbf{i}_0(P_{i,j})}}{2} P_{i,j} + \frac{1}{2}\mathbf{i}_0(P_{i,j})$$

Ce qui nous permet déduire la limite suivante :

$$\lim_{P_{i,j}\to+\infty} T_1^S(P_{i,j}) = +\infty$$

On encore on fait tender l'indice vers l'infini :

$$\lim_{i\to+\infty} T_1^S(P_{i,j}) = +\infty$$

Puisque $P_{1,j}$ et n sont constantes pour un super-chromologue bien déterminé donc les deux coefficients transformationnels globaux (principal et secondaire) sont constants pour un ce super-chromologue donc on peut déduire la limite suivante :

$$\lim_{i\to+\infty} \frac{\widetilde{\varphi}_n(P_{1,j})}{T_1^S(P_{i,j})} = 0$$

On peut conclure alors que:

$$\lim_{P_{,ij}\to+\infty} \frac{T_{n+1}^{S}(P_{i,j})}{T_1^S(P_{i,j})} = \lim_{P_{i,j}\to+\infty}\left(\frac{3^{\widetilde{M}_n(P_{i,j})}}{2^n} + \frac{\widetilde{\varphi}_n(P_{1,j})}{T_1^S(P_{i,j})}\right) = \frac{3^{\widetilde{M}_n(P_{1,j})}}{2^n} \qquad (e\ 6.18)$$

Comme on a:

$$\frac{3^{\widetilde{M}_n(P_{i,j})}}{2^n} < 1$$

On peut déduire qu'il existe un entier non nul $P_{i_{max},j}$ tel que pour tout entier $P_{i,j} > P_{i_{max},j}$

$$\frac{T_{n+1}^{S}(P_{i,j})}{T_1^S(P_{i,j})} < 1 \qquad (e\ 6.19)$$

Les suites du chromologue dont les générateurs possédant des rangs strictement supérieurs à $i_{max}$ sont de type $B_n^-$ autrement à partir de cette indice toutes les suites sont de ce type.

$P_{i_{max},j}$ est le plus grand terme de la suite génératrice du super-chromologue $B_n$ qui vérifie :

$$\frac{T_{n+1}^S(P_{i_{max},j})}{T_1^S(P_{i_{max},j})} > 1 \qquad (e\ 6.20)$$

On peut conclure que :

$$\begin{cases} \dfrac{T_{n+1}^S(P_{i,j})}{T_1^S(P_{i,j})} > 1\ si\ P_{1,j} \leq P_{ij} \leq P_{i_{max},j} \\ \dfrac{T_{n+1}^S(P_{i,j})}{T_m^S(P_{i,j})} < 1\ si\ P_{i,j} > P_{i_{max},j} \end{cases}$$

Puisque le nombre des générateurs est infini dans un chromologue infini peut contenir un nombre fini (il peut être  nul) des suites croissantes dont les générateurs sont $P_{1,j}, P_{2,j}, .., P_{i_{max}j}$. En contre partie, il contient une infinité des suites abaissante du rang allant de $(i_{max} + 1)$ à l'infini.

$P_{i_{max},j}$ est appelé le point d'inversement d'accroissement du chromologue considéré.

## 6.6 Sous classement des chromologues $B_n$

Les super-chromologues de type $B_n$ sont classées en deux catégories :

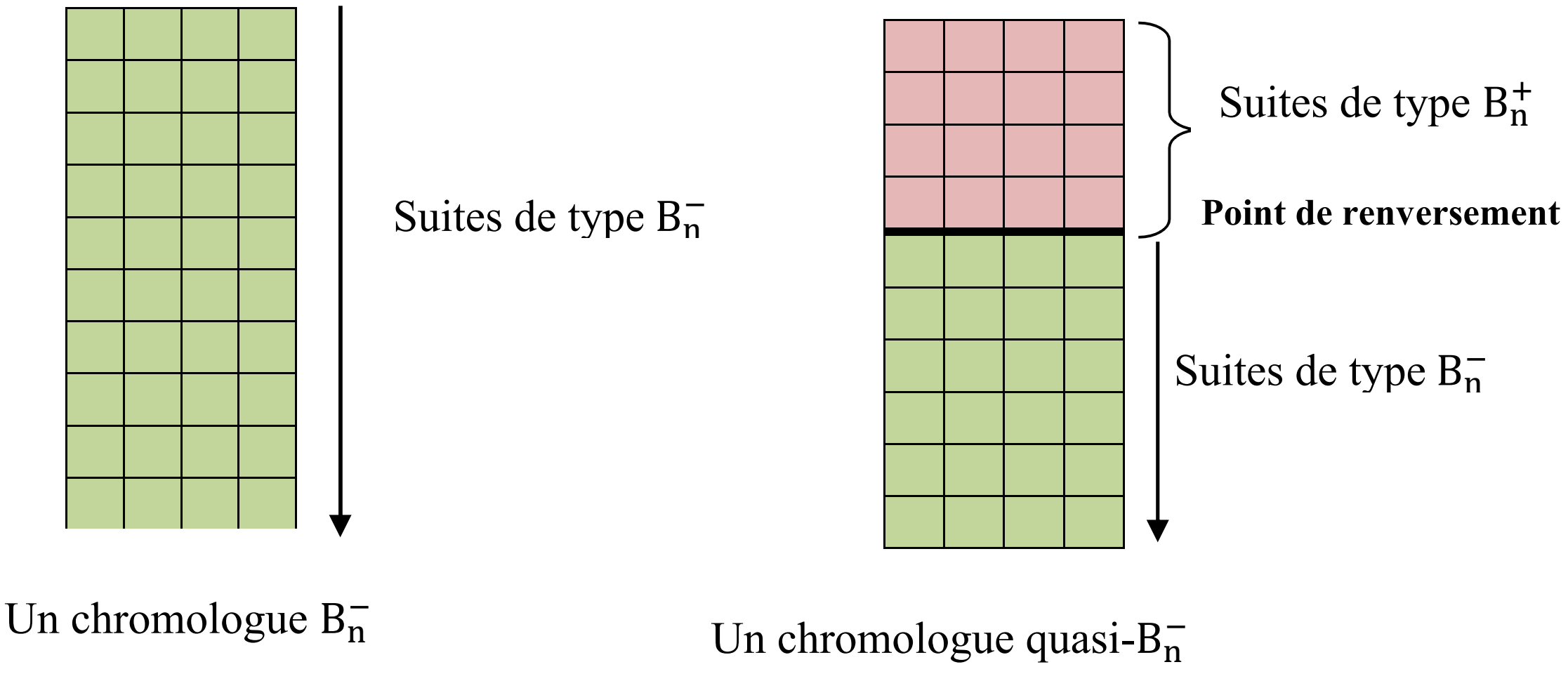

*Figure 11:* les deux catégories des chromologues $B_n$

-Les super-chromologues de type $B_n^-$ correspondent à des super-chromologues qui ne contiennent que des suites générées de type $B_n^-$

-Les super-chromologues quasi- $B_n^-$ correspondent aux super-chromologues qui contiennent un nombre fini des suites générées de type $B_n^+$ et un nombre infini des suites générées de type$B_n^-$

## 6.7 Prolongation parfaite d'un chromologue fractionnaire :

La prolongation parfaite d'un chromologue fractionnaire permet la construction d'un matrice qui contient tous les arrangements binaire ou chromatiques par prolongation d'un chromologue selon

- On considère un chromologue fractionnaire d'ordre dimensionnel n et qui contient $2^k$ lignes avec k un entier naturel non nul.
- On ajoute k colonnes a ce chromologue fractionnaire ce qui nous permet d'obtenir une matrice de k colonnes et $2^k$ lignes par prolongation.

Le tableau obtenu par prolongation renferme des suites de Collatz dont leurs structures constituent un tableau des arrangements chromatiques (ou binaires) complets.

**Exemple 6.3**

L'exemple suivant montre la prolongation parfaite du chromologue fractionnaire $H_f^T(7,4,2^4)$

| | n colonnes | | | | k colonnes | | | |
|---|---|---|---|---|---|---|---|---|
| 7 | 11 | 17 | 26 | 13 | 20 | 10 | 5 | 8 |
| 39 | 59 | 89 | 134 | 67 | 101 | 152 | 76 | 38 |
| 71 | 107 | 161 | 242 | 121 | 182 | 91 | 137 | 206 |
| 103 | 155 | 233 | 350 | 175 | 263 | 395 | 593 | 890 |
| 135 | 203 | 305 | 458 | 229 | 344 | 172 | 86 | 43 |
| 167 | 251 | 377 | 566 | 283 | 425 | 638 | 319 | 479 |
| 199 | 299 | 449 | 674 | 337 | 506 | 253 | 380 | 190 |
| 231 | 347 | 521 | 782 | 391 | 587 | 881 | 1322 | 661 |
| 263 | 395 | 593 | 890 | 445 | 668 | 334 | 167 | 251 |
| 295 | 443 | 665 | 998 | 499 | 749 | 1124 | 562 | 281 |
| 327 | 491 | 737 | 1106 | 553 | 830 | 415 | 623 | 935 |
| 359 | 539 | 809 | 1214 | 607 | 911 | 1367 | 2051 | 3077 |
| 391 | 587 | 881 | 1322 | 661 | 992 | 496 | 248 | 124 |
| 423 | 635 | 953 | 1430 | 715 | 1073 | 1610 | 805 | 1208 |
| 455 | 683 | 1025 | 1538 | 769 | 1154 | 577 | 866 | 433 |
| 487 | 731 | 1097 | 1646 | 823 | 1235 | 1853 | 2780 | 1390 |

$2^k$lignes

*Figure 12:* prolongation parfaite du chromologue fractionnaire $H_f^T(7,4,16)$

On peut utiliser les propriétés établies pour les tableaux des arrangements binaires complets pour déterminer les propriétés des chromologues lorsque l'ordre dimensionnel de la matrice de prolongation noté ici k tend vers l'infini.

**Corollaire 6.3** La proportion des suites abaissantes c'est à dire les suites générées de type $B_n^-$ dans un chromologue de type $B_n$ est pratiquement égal à 1.

**Remarque 6.1** Généralement la plupart des super-chromologues de type $B_n$ ne contiennent pas aucune suite croissante.

# 7 Les polychromologues de Collatz

**Définition 7.1**

On sait que chaque super-suite de la super-matrice représente une suite génératrice d'un chromologue bien déterminé. Ainsi la super-suite $X^S(P_{1,j}, n)$ représente la suite génératrice du chromologue $H^T(P_{1,j}, n)$ ci dessous :

$$H^T(P_{1,j}, n) = \left[\tilde{S}^r(P_{1,j}, n), \tilde{S}^r(P_{2,j}, n), \tilde{S}^r(P_{3,j} + 2x2^{n+1}, n), \dots, \tilde{S}^r(P_{k,j} + jx2^{n+1}, n) \dots\right]$$

Comme la super-matrice contient $2^n$ super-suites différentes donc on peut conclure qu'on peut définir $2^n$ chromologues différents à partir de toutes les super-suites d'une même super-matrice. L'ensemble constitué par les chromologues construits à partir d'un même super-système constitue un polychromologue homogène de Collatz qu'on le note comme suit :

$$\mathbb{H}^T(P_{1,1}, n) = \left[H^T(P_{1,1}, n), H^T(P_{1,2}, n), H^T(P_{1,3}, n), \dots., H^T(P_{1,2^n}, n)\right] \qquad \text{(e 7.1)}$$

Le couple $(P_{1,1}, n)$ represente les paramètres constructifs du polychromologue considéré.

La suite constituée par tous les générateurs fondamentaux de chaque chromologue représente la suite constructive du chromologue, elle est notée comme suit :

$$gr(P_{1,1}, n) = (P_{1,1}, P_{1,2}, P_{1,3}, \dots, P_{1,(2^n-1)}., P_{1,2^n}) \qquad \text{(e 7.2)}$$

-n représente l'ordre dimensionnel du polychromologue.

-$P_{1,1}$ : Générateur fondamental du polychromologue $\mathbb{H}^T(P_{1,1}, n)$.

La dimension d'un polychromologue est finie, elle correspond au nombre des chromologues constituant ce polychromologue. Ce nombre est égal à $2^n$ pour un polychromologue d'ordre dimensionnel n.

**Définition 7.2**

Pour une entier naturel n donné, on peut définir deux polychromologues :

Le polychromologue constitué par tous les chromologues dont leurs générateurs fondamentaux sont des entiers impairs.

$$\mathbb{H}^T(1,n) = \left[H^T(1,n), H^T(3,n), H^T(5,n), \dots, H^T\left(2^{n+1}-1, n\right)\right] \qquad \text{(e 7.3)}$$

Le deuxième polychromologue est constitué par tous les chromologues dont les générateurs fondamentaux sont des entiers pairs.

$$\mathbb{H}^T(2,n) = \left[H^T(2,n), H^T(4,n), H^T(6,n), \dots, H^T\left(2^{n+1}, n\right)\right] \qquad \text{(e 7.4)}$$

Chaque polychromologue est constitué exactement de $2^n$ chromologues.

**Définition 7.3**

Polychromologue fractionnaire des paramètres (n, q)

Un polychromologue complet est constituée de $2^n$ chromologues complets, chaque chromologue peut représenter sous forme d'un tableau à une infinité des lignes et de n colonnes. On sait qu'un chromologue fractionnaire est un chromologue qui contient un ensemble finie des suites générées et il est représenté sous d'un tableau de n colonnes et q lignes avec q un entier naturel non nul. Un polychromologue fractionnaire est un ensemble des $2^n$ chromologue fractionnaire tel que chaque chromologue constitué de q suites générées.

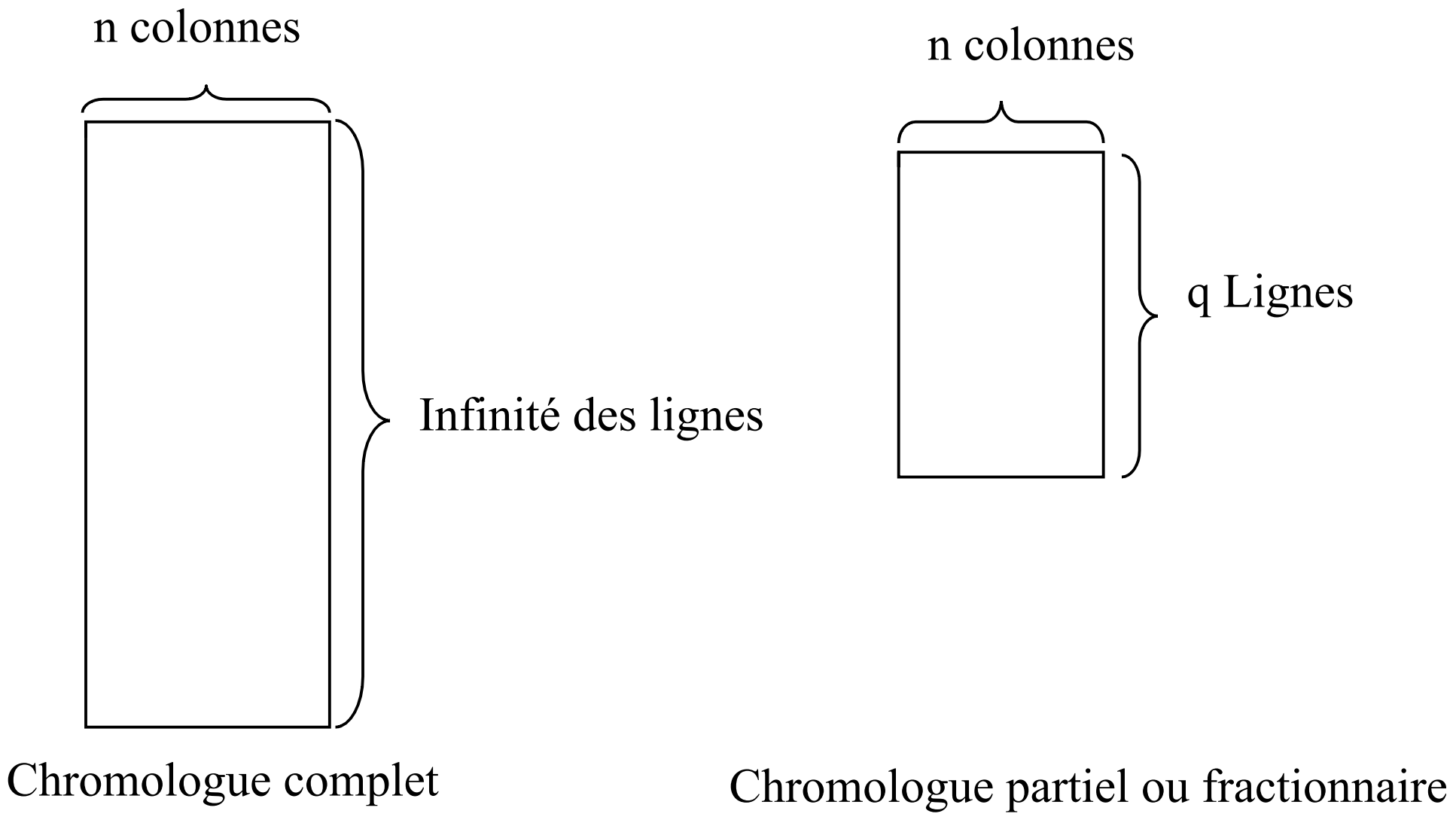

*Figure 13 :* chromollogue complet et chromologue fractionnaire

**Notations 7.1**

On considère un polychromologue complet de longueur n avec n un entier naturel non nul. On désigne par :

-$h_A(n)$ Le nombre des chromologues de type $A_n$ dans le polychromologue $\mathbb{H}^T(P,n)$.

-$h_B(n)$ Le nombre des chromologues de type $B_n$ dans $\mathbb{H}^T(P,n)$.

-$r_A(n)$ La proportion des chromologues de type $A_n$ dans $\mathbb{H}^T(P,n)$.

-$r_B(n)$ La proportion des chromologues de type $B_n$ dans $\mathbb{H}^T(P,n)$.

-Le polychromologue partiel des paramètres (P,n,q) est noté $\mathbb{H}_f^T(P,n,q)$

Le nombre total des chromologue du Collatz dans le polychromologue $\mathbb{H}^T(P,n)$ est comme suit :

$$h_A(n) + h_B(n) = 2^n \qquad \text{(e 7.5)}$$

Les deux proportions de deux types des chromologues dans le polychromologue complet considéré est comme suit :

$$\begin{cases} r_A(n) = \dfrac{h_A(n)}{2^n} \\ r_B(n) = \dfrac{h_A(n)}{2^n} \end{cases} \qquad \text{(e 7.6)}$$

La relation entre un polychromologue partiel des paramètres contractifs (P,n,q) et un polychromologue complet des paramètres constructifs (P, n) est comme suit :

$$\lim_{q\to+\infty} \mathbb{H}_f^T(P,n,q) = \mathbb{H}^T(P,n) \qquad \text{(e 7.7)}$$

**Théorème 7.1**

On considère un polychromologue complet de longueur n avec n un entier naturel non nul noté $\mathbb{H}^T(P,n)$. On désigne par :

-$r^+(n)$ la proportion des suites de type $S_n^+$ dans le polychromologue considéré.

-$r^-(n)$ la proportion des suites de type $S_n^-$ dans le polychromologue considéré.

Donc ces deux proportions vérifient les relations suivantes:

$$\begin{cases} r^+(n) = r_A(n) \\ r^-(n) = r_B(n) \end{cases} \qquad \text{(e 7.8)}$$

Autrement la proportion des suites générées finies strictement croissantes dans un polychromologue est égale à la proportion des suites de type $A_n$ dans ce même polychromologue alors que la proportion des suites générées décroissantes de longueur n est égale à la proportion des suites de type $B_n$.

**Démonstration**

On considère un polychromologue fractionnaire constitué de $2^n$ chromologues chaque chromologue renferme q suites générées de même longueur n qu'on le note $\mathbb{H}_f^T(P, n, q)$.

On désigne par :

-$k_{B^+}(n)$ Le nombre des suites de type $B_n^+$ dans $\mathbb{H}_f^T(P, n, q)$.

-$k_{B^-}(n)$ Le nombre des chromologue de type $B_n^-$ dans $\mathbb{H}_f^T(P, n, q)$.

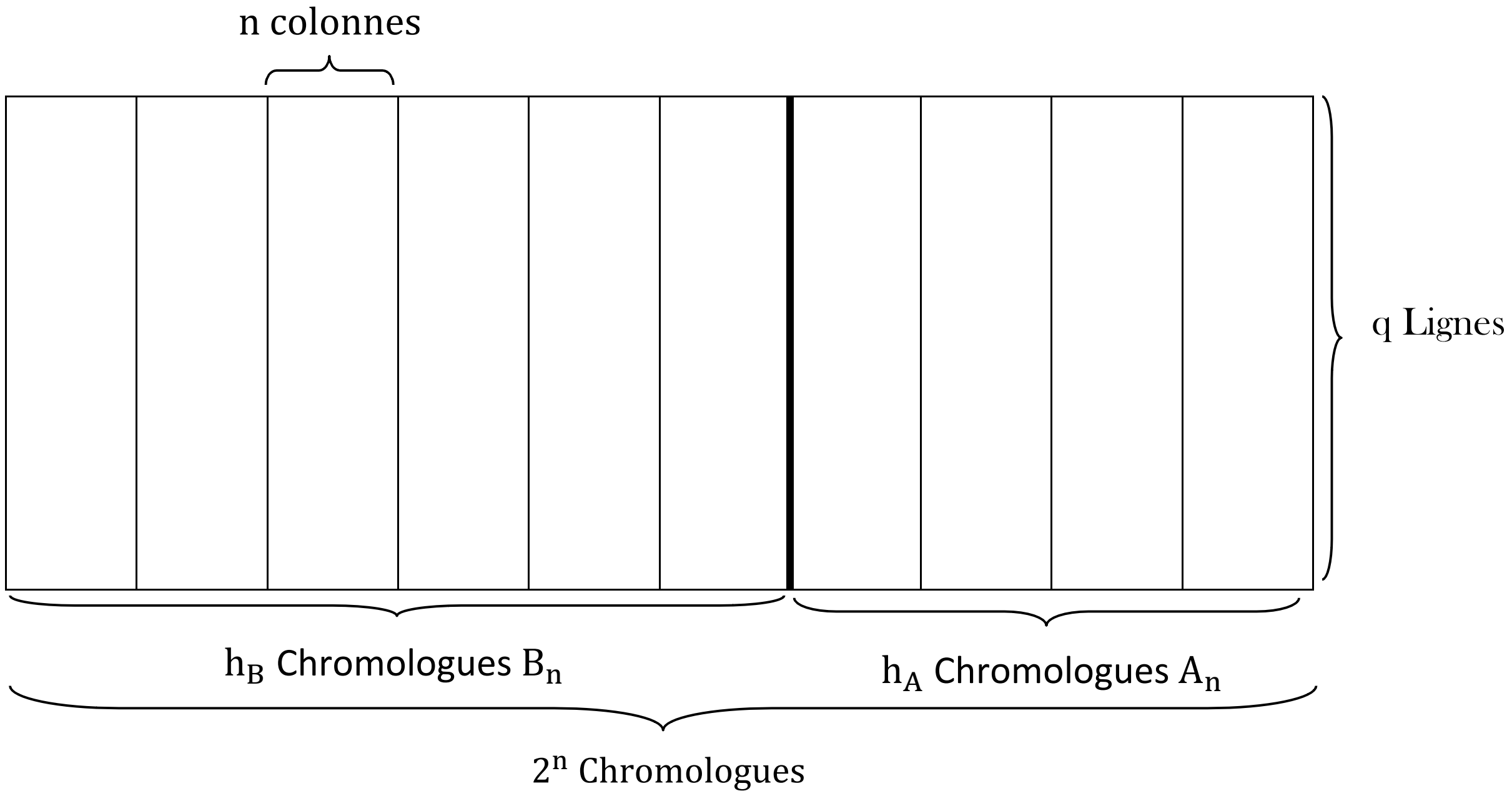

*Figure 14:* Structure d'un polychromologue fractionnaire de paramètres (n, q)

-Le nombre total des suites dans le polychromologue fractionnaire considéré est comme suit:

$$n_S = q(h_B(n) + h_A(n)) = 2^n q \qquad \text{(e 7.9)}$$

-Le nombre total des suites $S^+$:

Toutes les suites d'un chromologue de type $A_n$ sont croissantes donc elles sont de type $S_n^+$ alors qu'un chromologues de type $B_n$ peut contenir un nombre fini des suites croissantes donc le nombre total des suites croissantes dans un polychromologue fractionnaire d'extension longitudinale q est comme suit :

$$n_G \;= qh_A(n) + \; n_{B^+} \qquad \text{(e 7.10)}$$

Avec $n_{B^+}$ le nombre total des suites croissants dans tous les chromologue de type $B_n$, il est donné par:

$$n_{B^+} = \sum_{i=1}^{N_B} m_i \qquad \text{(e 7.11)}$$

La proportion des suites de type $A_n$ dans le polychromologue fractionnaire :

$$R_A(q,n) = \frac{qh_A(n)}{q(h_B(n) + h_A(n))} = \frac{h_A(n)}{h_B(n) + h_A(n)} = \frac{h_A(n)}{2^n} = r_A(n) \qquad \text{(e 7.12)}$$

La proportion des suites $S_n^+$ dans le polychromologue fractionnaire considéré :

$$r_G(q,n) = \frac{n_G}{n_S} = \frac{qh_A(n) + \; n_B^+}{q\big(h_B(n) + h_A(n)\big)}$$

$$= \frac{h_A(n)}{h_B(n) + h_A(n)} + \frac{n_B^+}{q(h_B(n) + h_A(n))}$$

$$= \frac{h_A(n)}{2^n} + \frac{n_{B^+}}{2^n q} = r_A(n) + \frac{n_{B^+}}{2^n q} \qquad \text{(e 7.13)}$$

On désigne par $R_G$ la proportion des suites croissantes dans le polychromologue complet donc cette proportion est correspond à la limite de la proportion de ces suites dans le polychromologue fractionnaire lorsque q tend vers l'infini :

$$R_G(n) = \lim_{q\to\infty} r_G(q,n)$$

$$= \lim_{q\to\infty} \left(\frac{h_A(n)}{2^n} + \frac{n_B^+}{2^n q}\right)$$

$$= \frac{h_A(n)}{2^n} = r_A(n) \qquad \text{(e 7.14)}$$

On peut conclure que la proportion des suites strictement croissantes dans un polychromologue complet d'extension n est égale à la proportion des suites de type $A_n$ c'est à dire à la proportion des suites de Collatz possédants des coefficients transformationnels globaux principaux strictement supérieurs à 1.

Si on désigne par $r^+(n)$ la proportion des suites strictement croissantes dans un polychromologue complet d'extension (ou ordre dimensionnel) n on peut écrire alors :

$$r^+(n) = r_A(n) = \frac{h_A(n)}{2^n} \qquad \text{(e 7.15)}$$

La proportion des suites abaissantes ou décroissantes dans un polychromologue complet est comme suit:

$$\begin{aligned} r^{-}(n) &= 1 - r^{+}(n) \\ &= 1 - \frac{h_A(n)}{2^n} \\ &= \frac{h_B(n)}{2^n} \\ &= r_B(n) \end{aligned} \qquad (e\ 7.16)$$

On peut conclure que la proportion des suites décroissantes (de type $S_n^-$) dans un polychromologue complet d'extension n est égale à la proportion des suites de type $B_n$ dans ce même polychromologue.

# 8 Les chromoformes parfaits de Collatz

**Définition 8.1**

On fait correspondre à chaque terme $P_j$ de la suite génératrice $Y^S(P_1, n)$ une suite générée de Collatz de générateur $P_j$ et de longueur n on obtient un ensemble constitué de $2^n$ suites générées de Collatz de même longueur n. L'ensemble obtenu est appelé le chromoforme parfait de Syracuse-Collatz d'ordre dimensionnel n. Il est noté $Z^T(P_1, n)$. Il est comme suit:

$$Z^T(P_1, n) = \left[\tilde{S}^r(P_1, n), \tilde{S}^r(P_2, n), \tilde{S}^r(P_3, n), \dots, \tilde{S}^r\left(P_j, n\right), \dots, \tilde{S}^r(P_{2^n}, n)\right] \qquad (e\ 8.1)$$

-La suite $Y^S(P_1, n)$ est appelée la suite génératrice du chromoforme $Z^T(P_1, n)$.

-n est l'ordre dimensionnel du chromoforme.

-$P_1$ est le générateur fondamental du chromoforme.

-Les entiers n et $P_1$ sont appelés les paramètres constructifs du chromoforme.

**Notations 8.1**

Le chromoforme parfait de Collatz construit à partir de la suite génératrice parfaite $Y^S(1, n)$ est comme suit:

$$Z^T(1, n) = \left[\tilde{S}^r(1, n), \tilde{S}^r(3, n), \tilde{S}^r(5, n), \dots, \tilde{S}^r(2^{n+1} - 1, n)\right] \qquad (e\ 8.2)$$

Le chromoforme parfait de Collatz construit à partir de la suite génératrice parfaite $Y^S(2, n)$ est comme suit:

$$Z^T(2, n) = \left[\tilde{S}^r(2, n), \tilde{S}^r(4, n), \tilde{S}^r(6, n), \dots, \tilde{S}^r(2^{n+1}, n)\right] \qquad (e\ 8.3)$$

**Définition 8.2**

On définit, le chromoforme structurel relatif au chromoforme parfait comme étant l'ensemble constitué par toutes les distributions structurelles linéaires des toutes les suites générées constituants le chromoforme numérique $Z^T(P_1, n)$. Il est noté comme suit:

$$Z^S(P_1, n) = [\tilde{L}^S(P_1, n), \tilde{L}^S(P_2, n), \tilde{L}^S(P_3, n), \dots, \tilde{L}^S(P_{2^n}, n)] \quad \text{(e 8.4)}$$

## 8.1 La matrice générée parfaite de Collatz

La représentation matricielle du chromoforme parfait $Z^T(P_1, n)$ de Collatz dont la suite génératrice est $Y^S(P_1, n)$ correspond à un tableau de $2^n$ lignes et de n colonnes. La ligne du rang i du tableau renferme les n premiers termes d'une suite générée de Collatz de générateur $P_i$ c'est à dire la suite générée $\tilde{S}^r(P_i, n)$ . Généralement, on peut ajouter une ligne en haut du tableau pour indiquer l'ordre transformationnel relatif à chaque terme de la suite comme montre le tableau ci-dessous. La matrice parfaite de Collatz générée par $Y^S(P_1, n)$ est notée $\mathbb{M}^T(P_1, n)$ , elle est comme suit:

*Tableau5:* Matrice parfaite $\mathbb{M}^T(P_1, n)$ relative au chromoforme $Z^T(P_1, n)$

| | $T_1^S(P_i)$ | $T_2^S(P_i)$ | | | | $T_{n-1}^S(P_i)$ | $T_n^S(P_i)$ |
|---|---|---|---|---|---|---|---|
| $P_1$ | $T_1^S(P_1)$ | $T_2^S(P_1)$ | | | | $T_{n-1}^S(P_1)$ | $T_n^S(P_1)$ |
| $P_2$ | $T_1^S(P_2)$ | $T_2^S(P_2)$ | | | | $T_{n-1}^S(P_2)$ | $T_n^S(P_2)$ |
| | | | | | | | |
| | | | | | | | |
| $P_{2^n}$ | $T_1^S(P_{2^n})$ | $T_2^S(P_{2^n})$ | | | | $T_{n-1}^S(P_{2^n})$ | $T_n^S(P_{2^n})$ |

La première colonne et la première ligne qui sont séparées du reste du tableau par deux lignes gras ne font pas partie de la matrice parfaite générée relative au chromoforme considéré.

**Exemple 8.1**

La matrice $\mathbb{M}^T(1,4)$ correspond à un tableau de 4 colonnes et de $2^4$ lignes. La première colonne contient tous les termes de la suite parfaite $Y^S(1,4)$ , elle constitue la suite génératrice de la matrice considérée.

4 colonnes

| | | | | |
|---|---|---|---|---|
| 1 | 2 | 1 | 2 | 1 |
| 3 | 5 | 8 | 4 | 2 |
| 5 | 8 | 4 | 2 | 1 |
| 7 | 11 | 17 | 26 | 13 |
| 9 | 14 | 7 | 11 | 17 |
| 11 | 17 | 26 | 13 | 20 |
| 13 | 20 | 10 | 5 | 8 |
| 15 | 23 | 35 | 53 | 80 |
| 17 | 26 | 13 | 20 | 10 |
| 19 | 29 | 44 | 22 | 11 |
| 21 | 32 | 16 | 8 | 4 |
| 23 | 35 | 53 | 80 | 40 |
| 25 | 38 | 19 | 29 | 44 |
| 27 | 41 | 62 | 31 | 47 |
| 29 | 44 | 22 | 11 | 17 |
| 31 | 47 | 71 | 107 | 161 |

$2^4$Lignes

*Figure 15:* La matrice parfaite $\mathbb{M}^{T}(1,4)$ relative au chromoforme parfait $Z^{T}(1,4)$

## 8.2 La représentation matricielle chromatisée d'un chromoforme parfait de Collatz : La chromatrice parfaite de Collatz:

Les chromatrices parfaites de Collatz sont des tableaux colorés obtenues par la coloration de différentes cases d'une matrice générée parfaite de Collatz suivant le principe suivant :

- Une case contient un entier impair sera colorée en bleu.
- Une case contient un entier pair sera colorée par une autre couleur ou bien

On effectue cette coloration sans élimination des entiers ou des termes des suites contenus dans les différentes cases. Les chromatrices nous permettent de visualiser simultanément les distributions numériques de différents termes des suites générées de Collatz et les distributions structurelles qui lui correspondent. La chromatrice relative au chromoforme parfait $Z^{T}(P_1, n)$ est notée $\mathbb{M}^{V}(P_1, n)$.

**Exemple 8.2**

La figure suivante correspond aux chromatrices parfaites $\mathbb{M}^{r}(1,4)$ et $\mathbb{M}^{r}(2,4)$ qui correspondent respectivement aux chromoformes parfaits $Z^{T}(1,4)$ et $Z^{T}(2,4)$

| | | | | |
|---|---|---|---|---|
| 1 | 2 | 1 | 2 | 1 |
| 3 | 5 | 8 | 4 | 2 |
| 5 | 8 | 4 | 2 | 1 |
| 7 | 11 | 17 | 26 | 13 |
| 9 | 14 | 7 | 11 | 17 |
| 11 | 17 | 26 | 13 | 20 |
| 13 | 20 | 10 | 5 | 8 |
| 15 | 23 | 35 | 53 | 80 |
| 17 | 26 | 13 | 20 | 10 |
| 19 | 29 | 44 | 22 | 11 |
| 21 | 32 | 16 | 8 | 4 |
| 23 | 35 | 53 | 80 | 40 |
| 25 | 38 | 19 | 29 | 44 |
| 27 | 41 | 62 | 31 | 47 |
| 29 | 44 | 22 | 11 | 17 |
| 31 | 47 | 71 | 107 | 161 |
| 2 | 1 | 2 | 1 | 2 |
| 4 | 2 | 1 | 2 | 1 |
| 6 | 3 | 5 | 8 | 4 |
| 8 | 4 | 2 | 1 | 2 |
| 10 | 5 | 8 | 4 | 2 |
| 12 | 6 | 3 | 5 | 8 |
| 14 | 7 | 11 | 17 | 26 |
| 16 | 8 | 4 | 2 | 1 |
| 18 | 9 | 14 | 7 | 11 |
| 20 | 10 | 5 | 8 | 4 |
| 22 | 11 | 17 | 26 | 13 |
| 24 | 12 | 6 | 3 | 5 |
| 26 | 13 | 20 | 10 | 5 |
| 28 | 14 | 7 | 11 | 17 |
| 30 | 15 | 23 | 35 | 53 |
| 32 | 16 | 8 | 4 | 2 |

*Figure 16:* Chromatrices $\mathbb{M}^{r}(1,4)$ et $\mathbb{M}^{r}(2,4)$

Les deux premières colonnes colorées en rouge brique correspondent aux suites génératrices noter que ces deux colonnes ne font pas parties du chromatrices.

## 8.3 La matrice structurelle binaire parfaite de Collatz relative à un chromoforme parfait

On fait correspondre à chaque matrice générée parfaite de Collatz un tableau contient le même nombre des colonnes et le même nombre des lignes tout en effectuant les modifications suivantes :

-Chaque entier pair est remplacé par 0.

-Chaque entier impair est remplacé par 1.

On obtient un tableau qui ne contient que des arrangements avec répétition de deux entiers 0 et 1. Le tableau obtenu est appelée la matrice structurelle binaire de Collatz relative au chromoforme parfait considéré.

La matrice structurelle binaire relative au chromoforme parfait $Z^{T}(P_1, n)$ est notée $\mathbb{M}^{b}(P_1, n)$.

*Tableau 6:* Matrice structurelle binaire relative au chromoforme parfait $Z^T(P_1, n)$

| | $\mathbf{i}_1(P_i)$ | $\mathbf{i}_2(P_i)$ | | | $\mathbf{i}_{n-1}(P_i)$ | $\mathbf{i}_n(P_i)$ |
|---|---|---|---|---|---|---|
| $\tilde{L}^{s}(P_1, n)$ | $\mathbf{i}_1(P_1)$ | $\mathbf{i}_2(P_1)$ | | | $\mathbf{i}_{n-1}(P_1)$ | $\mathbf{i}_n(P_1)$ |
| $\tilde{L}^{s}(P_2, n)$ | $\mathbf{i}_1(P_2)$ | $\mathbf{i}_2(P_2)$ | | | $\mathbf{i}_{n-1}(P_2)$ | $\mathbf{i}_n(P_2)$ |
| $\tilde{L}^{s}(P_3, n)$ | $\mathbf{i}_1(P_3)$ | $\mathbf{i}_3(P_3)$ | | | $\mathbf{i}_{n-1}(P_3)$ | $\mathbf{i}_n(P_3)$ |
| | | | | | | |
| | | | | | | |
| $\tilde{L}^{s}(P_{2^n}, n)$ | $\mathbf{i}_1(P_{2^n})$ | $\mathbf{i}_2(P_{2^n})$ | | | $\mathbf{i}_{n-1}(P_{2^n})$ | $\mathbf{i}_n(P_{2^n})$ |

La première colonne et la première ligne qui sont séparées du reste du tableau par deux traits gras ne font pas partie de la matrice structurelle.

## 8.4 La matrice structurelle chromatique parfaite relative au chromoforme parfait

La matrice structurelle chromatique représente la conversion chromatique de la matrice générée parfaite elle est équivalente à la matrice structurelle binaire. Au lieu de remplacer les termes par 0 ou 1 selon leurs parités, on fait colorer chaque case tout en éliminant les termes contenus dans les différentes cases  selon les règles suivantes :

-Une case contient un entier impair sera colorée en bleu.

-Une case contient un entier pair on le laisse sans aucune coloration ou on dit qu'elle est colorée en blanc.

La matrice structurelle chromatique relative au chromoforme $Z^T(P_1, n)$ est notée $\mathbb{M}^c(P_1, n)$

**Exemple 8.3**

La figure suivante représente un exemple d'une matrice structurelle chromatique et d'une matrice structurelle binaire.

| | | | |
|---|---|---|---|
| 1 | 2 | 1 | 2 |
| 3 | 5 | 8 | 4 |
| 5 | 8 | 4 | 2 |
| 7 | 11 | 17 | 26 |
| 9 | 14 | 7 | 11 |
| 11 | 17 | 26 | 13 |
| 13 | 20 | 10 | 5 |
| 15 | 23 | 35 | 53 |

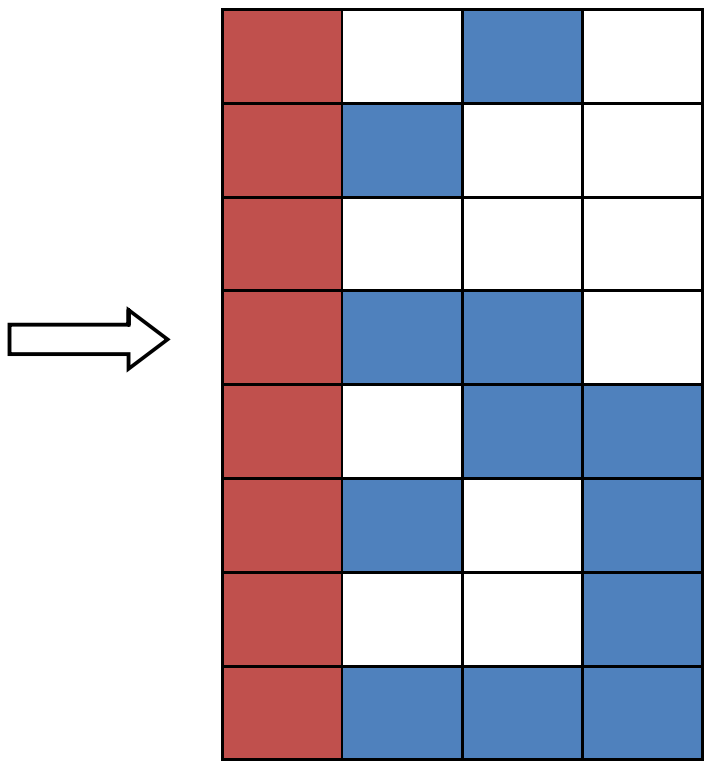

| | | | |
|---|---|---|---|
| 1 | 0 | 1 | 0 |
| 3 | 1 | 0 | 0 |
| 5 | 0 | 0 | 0 |
| 7 | 1 | 1 | 0 |
| 9 | 0 | 1 | 1 |
| 11 | 1 | 0 | 1 |
| 13 | 0 | 0 | 1 |
| 15 | 1 | 1 | 1 |

*Figure 17:* Les deux matrices structurelles relatives au chromoforme $Z^T(1,3)$.

On remarque bien que les deux tableaux constituent deux tableaux d'arrangements complets de base d'arrangement (0,1) et de base chromatique (bleu, blanc).

Un autre exemple illustre les matrices structurelles (binaire et chromatique) relatives à la matrice du chromoforme parfait $Z^T(2,3)$.

| | | | |
|---|---|---|---|
| 2 | 1 | 2 | 1 |
| 4 | 2 | 1 | 2 |
| 6 | 3 | 5 | 8 |
| 8 | 4 | 2 | 1 |
| 10 | 5 | 8 | 4 |
| 12 | 6 | 3 | 5 |
| 14 | 7 | 11 | 17 |
| 16 | 8 | 4 | 2 |

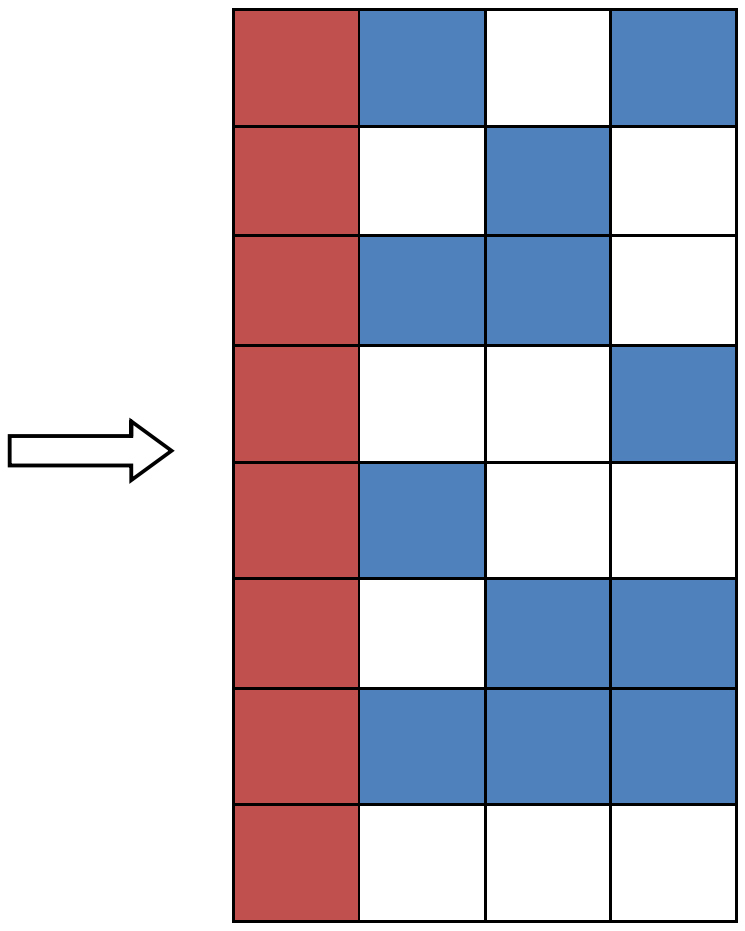

| | | | |
|---|---|---|---|
| 2 | 1 | 0 | 1 |
| 4 | 0 | 1 | 0 |
| 6 | 1 | 1 | 0 |
| 8 | 0 | 0 | 1 |
| 10 | 1 | 0 | 0 |
| 12 | 0 | 1 | 1 |
| 14 | 1 | 1 | 1 |
| 16 | 0 | 0 | 0 |

*Figure 18:* Les deux matrices structurelles relatives au chromoforme $Z^T(2,3)$.

De même dans ce cas, on remarque bien que les deux tableaux de la conversion structurelle constituent deux tableaux d'arrangements complets de base d'arrangement (0,1) et de base chromatique (bleu, blanc).

**Remarque importante 8.1**

Les deux types de conversion structurelle binaire ou bien chromatique reposent sur le même principe, ils nous permettent d'obtenir des tableaux qui ne contiennent que des arrangements avec répétition et ils traduits les répartitions des entiers pairs et impairs dans une distribution générée donnée. L'utilisation de l'une de deux configurations pour étudier le comportement structurelle des suites de Collatz conduit toujours aux mêmes résultats et aux mêmes propriétés autrement les deux représentations structurelles relatives à un même chromoforme sont strictement équivalentes.

**Lemme 8.1**

Toutes les suites générées de même longueur n d'un même chromoforme parfait de Collatz d'ordre dimensionnel n sont deux à deux non isoformes (ou non chromatiques).

**Démonstration**

D'après le théorème(e 3.3) pour que deux suites générées $\tilde{S}^r(P_{j_1}, n)$ et $\tilde{S}^r(P_{j_2}, n)$ de même longueur n soient isochromatiques, il faut que la différence en valeur absolue entre leurs générateurs est égale à un multiple de $2^{n+1}$ c'est à dire il faut que :

$$|P_{j_1} - P_{j_2}| = 2^{n+1}k \quad \text{(e 8.5)}$$

Avec k un entier naturel non nul.

En contre partie, on sait que la différence entre deux générateurs quelconques $P_{j_1}$ et $P_{j_2}$ d'une même suite génératrice parfaite $Y^S(P_1, n)$ d'ordre dimensionnel n est majorée par $(2^{n+1} - 2)$ puisque cette suite correspond à une suite arithmétique finie de longueur $2^n$ et de raison 2 et par suite la différence entre deux termes quelconques de la suite génératrice est inférieure ou égale à la différence entre le premier terme et le dernier terme de cette suite génératrice.

Pour la suite génératrice $Y^S(1, n)$ le premier terme est égale à 1 et le dernier terme correspond à $(2^{n+1} - 1)$ alors que pour la suite $Y^S(2, n)$ le premier terme est 2 et le dernier est $2^{n+1}$ ce qui nous permet d'écrire:

$$|P_{j_1} - P_{j_2}| \leq 2^{n+1} - 2 < 2^{n+1} \text{ avec } j_1 \neq j_2 \quad \text{(e 8.6)}$$

On peut conclure que les deux suites $\tilde{S}^r(P_{j_1}, n)$ et $\tilde{S}^r(P_{j_2}, n)$ ne sont pas isochromatiques. Cette propriété est vraie pour tous couple de la suite génératrice ce qui nous permet de conclure que les suites générées d'un même chromoforme sont deux à deux non isoformes.

**Théorème 8.1**

Pour tout entier naturel non nul n, les matrices structurelles d'ordre dimensionnel n relatives aux chromoforme s parfaits de Syracuse-Collatz d'ordre dimensionnel n sont des tableaux des arrangements complets autrement tout chromoforme parfait de Collatz d'ordre dimensionnel n admet une matrice structurelle binaire qui correspond à un tableau d'arrangements binaires complets et admet une matrice structurelle chromatique qui correspond à un tableau des arrangements chromatiques complets.

**Exemple 8.4**

la conversion des termes du la matrice $\mathbb{M}^T(1,4)$ en 0 et 1 selon leurs parités conduit à un tableau des arrangements binaires complets qui contient tous les arrangements binaires avec répétition de longueur 4 comme montre la figure suivante :

| | | | | |
|---|---|---|---|---|
| 2 | 1 | 2 | 1 | 2 |
| 4 | 2 | 1 | 2 | 1 |
| 6 | 3 | 5 | 8 | 4 |
| 8 | 4 | 2 | 1 | 2 |
| 10 | 5 | 8 | 4 | 2 |
| 12 | 6 | 3 | 5 | 8 |
| 14 | 7 | 11 | 17 | 26 |
| 16 | 8 | 4 | 2 | 1 |
| 18 | 9 | 14 | 7 | 11 |
| 20 | 10 | 5 | 8 | 4 |
| 22 | 11 | 17 | 26 | 13 |
| 24 | 12 | 6 | 3 | 5 |
| 26 | 13 | 20 | 10 | 5 |
| 28 | 14 | 7 | 11 | 17 |
| 30 | 15 | 23 | 35 | 53 |
| 32 | 16 | 8 | 4 | 2 |
| 2 | 1 | 0 | 1 | 0 |
| 4 | 0 | 1 | 0 | 1 |
| 6 | 1 | 1 | 0 | 0 |
| 8 | 0 | 0 | 1 | 0 |
| 10 | 1 | 0 | 0 | 0 |
| 12 | 0 | 1 | 1 | 0 |
| 14 | 1 | 1 | 1 | 0 |
| 16 | 0 | 0 | 0 | 1 |
| 18 | | 0 | 1 | 1 |
| 20 | 0 | 1 | 0 | 0 |
| 22 | 1 | 1 | 0 | 1 |
| 24 | 0 | 0 | 1 | 1 |
| 26 | 1 | 0 | 0 | 1 |
| 28 | 0 | 1 | 1 | 1 |
| 30 | 1 | 1 | 1 | 1 |
| 32 | 0 | 0 | 0 | 0 |

*Figure 19:* Chromatrice parfaite de Collatz d'ordre 4 et sa conversion structurelle

**Démonstration**

On considère le chromoforme parfait de Collatz et le chromoforme structurel qui lui correspond :

$$Z^T(P_1, n) = [\tilde{S}^r(P_1, n), \tilde{S}^r(P_2, n), \tilde{S}^r(P_3, n), \dots, \tilde{S}^r(P_j, n), \dots, \tilde{S}^r(P_{2^n}, n)]$$

Les suites générées d'un même chromoforme sont deux à deux non isochromatique (ou non isoformes) c'est à dire que chaque suites générée de longueur n possède sa propre distribution structurelle linéaire qui est différente à toutes les autres distributions structurelles linaires relatives aux autres suites de chromoforme.

Pour tous entiers naturels non nuls $j_1$ et $j_2$ tel que $1 \leq j_2 \leq 2^n$ et $1 \leq j_1 \leq 2^n$

$$\text{si } j_1 \neq j_2 \text{ alors } \tilde{L}^{s}(P_{j_1}, n) \neq \tilde{L}^{s}(P_{j_2}, n) \qquad \text{(e 8.7)}$$

Ce qui nous permet de conclure que le chromoforme structurel parfait relatif au chromoforme fondamental de Collatz contient exactement $2^n$ arrangements avec répétition de longueur n qui sont deux à deux distincts.

$$Z^S(P_1, n) = [\tilde{L}^{s}(P_1, n), \tilde{L}^{s}(P_2, n), \tilde{L}^{s}(P_3, n), \dots, \tilde{L}^{s}(P_j, n), \dots, \tilde{L}^{s}(P_{2^n}, n)]$$

Comme ce nombre correspond au nombre total des tous les arrangements avec répétition dans un tableau d'arrangements binaires complets d'ordre dimensionnel n, nous pouvons conclure que le chromoforme structurel parfait de Syracuse contient tous les arrangements possibles avec répétition de longueur n de deux entiers 1 et 0. Comme la matrice structurelle binaire est une représentation matricielle du chromoforme structurel on peut déduire que cette matrice correspond à un tableau d'arrangements binaires complets.

**Corollaire 8.1**

Soient n et k deux entiers naturels non nuls tel que $(k \leq n)$. Le nombre des suites de Collatz qui contiennent exactement k entiers impairs dans une matrice générée parfaite de Collatz d'ordre dimensionnel n est donné par la formule suivante :

$$N_I(n,k) = C_n^k = \frac{n!}{k!\,(n-1)!} \qquad \text{(e 8.8)}$$

**Démonstration**

La matrice structurelle est un tableau des arrangements binaires complets avec répétition, il contient tous les arrangements possibles avec répétition de deux indicateurs de parité (0,1). En utilisant les formules et les lois relatives aux arrangements avec répétition et aux combinaisons on peut déterminer le nombre de tous les arrangements conditionnels possibles définie comme suit: arrangement avec répétition de deux objets rangés dans n cases sachant que le premier objet occupe k cases parmi les n cases vides donc ce nombre correspond à la combinaison de k éléments parmi n éléments comme suit :

$$C_n^k = \frac{n!}{k!\,(n-1)!}$$

**Lemme 8.2**

Soit $\mathbb{M}^T(P,2n)$ la matrice parfaite de Collatz d'ordre dimensionnel 2n. On désigne par $a_{2n}$ et $b_{2n}$ respectivement le nombre des suites de type $A_{2n}$ et le nombre des suites de type $B_{2n}$ dans la matrice $\mathbb{M}^T(P,2n)$ alors $a_{2n}$ et $b_{2n}$ s'écrits comme suit :

$$\begin{cases} b_{2n} = a_{2n} + C_{2n}^{n} + 2\sum_{k=1}^{e} C_{2n}^{n+k} \\ a_{2n} = \sum_{k=0}^{n-e-1} C_{2n}^{k} \end{cases} \qquad \text{(e 8.9)}$$

Avec :

$$e = \alpha_{2n} - n$$

$$\alpha_{2n} = E\left(\frac{\ln(2)}{\ln(3)}(2n)\right) \text{est le point de renversement de } \mathbb{Z}^T(P, 2n)$$

**Démonstration**

Considérons le tableau des arrangements binaires complets d'ordre dimensionnel 2n. On désigne par k le nombre des entiers impairs dans une suite quelconque de la matrice $\mathbb{M}^T(P, 2n)$. Les valeurs de k varient de 0 à 2n comme suit :

$$k \in \{0,1,2,3, \dots, \alpha_{2n}, \alpha_{2n} + 1, \dots, 2n - 1, 2n\}$$

Les suites de type $B_{2n}$ possèdent un nombre des entiers impairs inférieur ou égal au coefficient de renversement $\alpha_{2n}$ donc le nombre total des suites $B_{2n}$ est donné par :

$$b_{2n} = C_{2n}^0 + C_{2n}^1 + \cdots + C_{2n}^{\alpha_{2n}-1} + C_{2n}^{\alpha_{2n}} = \sum_{k=0}^{\alpha_{2n}} C_{2n}^k \tag{e 8.10}$$

Les suites de type $A_{2n}$ possèdent un nombre des entiers impairs strictement supérieur à $\alpha_{2n}$ donc le nombre des suites de type $A_{2n}$ dans la matrice considérée est comme suit :

$$a_{2n} = C_{2n}^{\alpha_{2n}+1} + C_{2n}^{\alpha_{2n}+2} + \cdots C_{2n}^{2n-1} + C_{2n}^{2n} = \sum_{k=\alpha_{2n}+1}^{2n} C_{2n}^k \tag{e 8.11}$$

En utilisant la relation de combinaison suivante : $C_{2n}^k = C_{2n}^{n-k}$ avec $k \leq n$ on peut déduire :

$$a_{2n} = C_{2n}^0 + C_{2n}^1 + \cdots + C_{2n}^{2n-\alpha_{2n}-2} + C_{2n}^{2n-\alpha_{2n}-1} \tag{e 8.12}$$

On sait que $\alpha_{2n} > n$ donc on pose $e = (\alpha_{2n} - n)$, on peut tirer la relation suivante :

$$2n - \alpha_{2n} = n - e \tag{e 8.13}$$

Remplaçons $(2n - \alpha_{2n})$ par $(n - e)$ dans l'expression de $a_{2n}$ ce qui nous permet d'écrire:

$$\begin{aligned} a_{2n} &= C_{2n}^0 + C_{2n}^1 + \cdots + C_{2n}^{2n-\alpha_{2n}-2} + C_{2n}^{2n-\alpha_{2n}-1} \\ &= C_{2n}^0 + C_{2n}^1 + \cdots + C_{2n}^{n-e-2} + C_{2n}^{n-e-1} \end{aligned} \tag{e 8.14}$$

On remplace $\alpha_{2n}$ par m+n dans l'expression de $b_{2n}$, on obtient:

$$\begin{aligned} b_{2n} &= C_{2n}^0 + C_{2n}^1 + \cdots + C_{2n}^{e+n-1} + C_{2n}^{e+n} \\ &= (C_{2n}^0 + C_{2n}^1 + \cdots + C_{2n}^{n-e-2} + C_{2n}^{n-e-1}) + (C_{2n}^{n-e} + C_{2n}^{n-e+1} + \cdots + C_{2n}^{n} + \cdots + C_{2n}^{e+n}) \\ &= a_{2n} + (C_{2n}^{n-e} + \dots + C_{2n}^{n-1} + C_{2n}^{n} + C_{2n}^{n+1} \dots + C_{2n}^{n+e}) \end{aligned} \tag{e 8.15}$$

En utilisant la relation $C_{2n}^k = C_{2n}^{2n-k}$ on peut écrire:

$$C_{2n}^{n-e} + C_{2n}^{n-e+1} + \cdots + C_{2n}^{n-1} = C_{2n}^{n+e} + C_{2n}^{n+e-1} + \cdots + C_{2n}^{n+1}$$

$$= \sum_{k=1}^{e} C_{2n}^{n+k} \qquad \text{(e 8.16)}$$

Ce qui nous permet d'écrire :

$$C_{2n}^{n-e} + .. + C_{2n}^{n-1} + C_{2n}^{n} + C_{2n}^{n+1} + \cdots + C_{2n}^{n+e} = C_{2n}^{n} + 2(C_{2n}^{n+1} + \cdots + C_{2n}^{n+e})$$

$$= C_{2n}^{n} + 2\sum_{k=1}^{e} C_{2n}^{n+k}$$

L'expression de $b_{2n}$ s'écrit comme suit :

$$b_{2n} = a_{2n} + C_{2n}^{n} + 2\sum_{k=1}^{e} C_{2n}^{n+k} \qquad \text{(e 8.17)}$$

**Théorème 8.2**

Soit $\mathbb{M}^{T}(P, 2n)$ la matrice parfaite de Collatz d'ordre dimensionnel 2n et on désigne par $a_{2n}$ et $b_{2n}$ respectivement le nombre des suites de type $A_{2n}$ et le nombre des suites de type $B_{2n}$ dans la matrice considérée. Lorsque l'ordre dimensionnel de la matrice tend vers l'infini, Les proportions de ces deux types des suites vérifient les limites suivantes :

$$\begin{cases} \lim\limits_{n \to +\infty} r_B(2n) = 1 \\ \lim\limits_{n \to +\infty} r_A(2n) = 0 \end{cases} \qquad \text{(e 8.18)}$$

Avec:

$$\begin{cases} r_B(2n) = \dfrac{b_{2n}}{2^{2n}} \text{ Proportion des suites de type } B_{2n} \text{ dans } \mathbb{M}^{T}(P, 2n) \\ r_A(2n) = \dfrac{a_{2n}}{2^{2n}} \text{ Proportion des suites de type } A_{2n} \text{ dans } \mathbb{M}^{T}(P, 2n) \end{cases} \qquad \text{(e 8.19)}$$

**Stratégie de la démonstration**

On cherche à démontrer que:

$$\lim_{n \to +\infty} \left(\frac{a_{2n}}{2^{2n}}\right) = 0 \qquad \text{(e 8.20)}$$

Équivaut à:

$$\lim_{n \to +\infty} \left(\frac{2^{2n}}{a_{2n}}\right) = +\infty \qquad \text{(e 8.21)}$$

Avec :

$2^{2n}$ est le nombre total des suites générées dans la matrice $\mathbb{M}^{T}(P, 2n)$.

$a_{2n}$ est le nombre total des suites générées de type $A_{2n}$ dans la matrice $\mathbb{M}^{T}(P, 2n)$.

On sait que :

$$\frac{2^{2n}}{a_{2n}} = \frac{a_{2n} + b_{2n}}{a_{2n}} = \frac{2a_{2n} + C_{2n}^{n} + 2(C_{2n}^{n-e} + C_{2n}^{n-e+1} + \cdots + C_{2n}^{n-1})}{a_{2n}}$$

$$= 2 + \frac{C_{2n}^{n}}{a_{2n}} + 2\,\frac{C_{2n}^{n-e} + C_{2n}^{n-e+1} \ldots + C_{2n}^{n-1}}{a_{2n}}$$

Donc il suffit de démontrer que :

$$\lim_{n\to+\infty} \left(\frac{C_{2n}^{n}}{a_{2n}} + \frac{C_{2n}^{n-e} + C_{2n}^{n-e+1} \ldots + C_{2n}^{n-1}}{a_{2n}}\right) = +\infty$$

On pose :

$$f(2n) = C_{2n}^{n-e} + C_{2n}^{n-e+1} + \cdots + C_{2n}^{n-1} + C_{2n}^{n} \sum_{i=0}^{e} C_{2n}^{n-i}$$

$f(2n)$ est composée de $(e+1)$ termes sous la forme de $C_{2n}^{n-i}$

D’autre part, on sait que :

$$a_{2n} = C_{2n}^{0} + C_{2n}^{1} + \cdots + C_{2n}^{n-e-2} + C_{2n}^{n-e-1}$$

On détermine une partie de $a_{2n}$ qu’on la note $g(2n)$ qui vérifie les conditions suivantes :
-$g(2n)$ contient exactement $(e+1)$ termes sous la forme $C_{2n}^{k}$
-$g(2n)$ n’est pas négligeable par rapport $a_{2n}$ donc quelque soit la valeur de n, elle existe une constante non nul C tel que :

$$C\,g(2n) \geq a_{2n} \qquad \text{(e 8.23)}$$

Dans ce cas si on montre que :

$$\lim_{n\to+\infty} \frac{f(2n)}{g(2n)} = +\infty$$

On peut déduire que :

$$\lim_{n\to+\infty} \frac{f(2n)}{a_{2n}} = +\infty$$

Ce qui implique :

$$\lim_{n\to+\infty} \frac{2^{2n}}{a_{2n}} = +\infty$$

On peut vérifier que la partie recherchée est correspond à la fonction suivante :

$$g(2n) = \sum_{k=n-2e-1}^{n-e-1} C_{2n}^{k} \qquad \text{(e 8.24)}$$

En effet, on effectue une partition de $a_{2n}$ comme suit :

$$a_{2n} = \sum_{k=0}^{n-3e-3} C_{2n}^{k} + \sum_{k=n-3e-2}^{n-2e-2} C_{2n}^{k} + \sum_{k=n-2e-1}^{n-e-1} C_{2n}^{k} \qquad (e\ 8.25)$$

On pose :

$$\begin{cases} g_1(2n) = \sum_{k=0}^{n-3e-3} C_{2n}^{k} \\ g_2(2n) = \sum_{k=n-3e-2}^{n-2e-2} C_{2n}^{k} \\ g(2n) = \sum_{k=n-2e-1}^{n-e-1} C_{2n}^{k} \end{cases} \qquad (e\ 8.26)$$

On sait que pour tous entiers naturels non nuls $k_2$ et $k_1$ tel que $0 \leq k_1 \leq k_2 \leq n$, on a :

$$C_{2n}^{k_2} > C_{2n}^{k_1}$$

Comme g(2n)et $g_2(2n)$ contiennent le même nombre des éléments qui sont sous la forme $C_{2n}^{k}$ si on fait comparer ces termes deux à deux en utilisant l'inégalité précédente on obtient l'inégalité suivante :

$$\sum_{k=n-2e-1}^{n-e-1} C_{2n}^{k} > \sum_{k=n-3e-2}^{n-2e-2} C_{2n}^{k} \qquad (e\ 8.27)$$

Équivaut à :

$$g(2n) > g_2(2n)$$

$g_1(2n)$ contient $(n-3e-2)$ termes or on sait que :

$$e > 0.25n$$

Ce qui nous permet de déduire :

$$n - 3e - 2 < 0.25n - 2 < e + 1$$

$g_1(2n)$ contient un nombre des éléments inférieur à $(e+1)$ donc en utilisant l'inégalité on peut déduire que :

$$\sum_{k=n-2e-1}^{n-e-1} C_{2n}^{k} > \sum_{k=0}^{n-3e-3} C_{2n}^{k}$$

C'est à dire que:

$$g(2n) > g_1(2n)$$

On peut déduire les inégalités suivantes de tous ce qui précède :

$$g(2n) > g_2(2n) > g_1(2n) \Longrightarrow 3g(2n) > g_1(2n) + g_2(2n) + g(2n)$$

Comme on a :

$$a_{2n} = g_1(2n) + g_2(2n) + g(2n)$$

Donc on peut déduire que:

$$3g(2n) > a_{2n}$$

Si on montre que :

$$\lim_{n\to+\infty} \frac{f(2n)}{g(2n)} = +\infty$$

On peut déduire la limite recherchée.

$$\frac{f(2n)}{g(2n)} = \frac{C_{2n}^{n} + C_{2n}^{n-1} + \cdots + C_{2n}^{n-e}}{C_{2n}^{n-e-1} + C_{2n}^{n-e-2} + \cdots + C_{2n}^{n-2e-1}}$$

$f(2n)$ et $g(2n)$ peuvent s'écrire comme suit :

$$f(2n) = \sum_{i=0}^{i=e} C_{2n}^{n-i} \ \text{ et } \ g(2n) = \sum_{i=0}^{i=e} C_{2n}^{n-e-1-i}$$

Pour montrer la limite ci- dessus, On montre tout d'abord que pour tout entier naturel i tel que $0 \leq i \leq e$:

$$\lim_{n\to+\infty} \frac{C_{2n}^{n-i}}{C_{2n}^{n-e-1-i}} = +\infty$$

Une fois la limite précédente est démontrée on peut procéder comme suit pour démontrer la limite recherchée. Si l'équation est vraie donc ceci signifie que pour tout entier naturel $M > 0$ il existe un entier naturel non nul $n_i$ tel que pour tout entier naturel $n \geq n_i$

$$\frac{C_{2n}^{n-i}}{C_{2n}^{n-n-1-i}} > M \Longrightarrow C_{2n}^{n-i} > M\, C_{2n}^{(n-e-1)-i}$$

Dans ce cas, on peut déduire que pour tout entier naturel $M > 0$ il existe un entier naturel non nul $N_0$ qui vérifie :

$$N_0 = \max\{n_0, n_1, \ldots., n_e\}$$

Tel que Pour tout entier naturel $n \geq N_0$ et pour tout entier $0 \leq i \leq e$, on peut écrire :

$$C_{2n}^{n-i} > M\, C_{2n}^{n-e-1-i}$$

Ce qui nous permet d'écrire:

$$\sum_{i=0}^{e} C_{2n}^{n-i} > M \sum_{i=0}^{e} C_{2n}^{n-e-1-i}$$

Équivaut à :

$$f(2n) > M\, g(2n) \Longrightarrow \frac{f(2n)}{g(2n)} > M$$

On peut conclure par suite que :

$$\lim_{n\to+\infty} \frac{f(2n)}{g(2n)} = +\infty$$

**Démonstration**

f(2n) est composée de $(e + 1)$ termes, chaque terme est sous la forme $C_{2n}^{n-i}$

$$f(2n) = \sum_{i=0}^{i=e} C_{2n}^{n-i}$$

g(2n) est composée de $(e + 1)$ termes, chaque terme est sous la forme $C_{2n}^{(n-e-1)-i}$

$$g(2n) = \sum_{i=0}^{i=e} C_{2n}^{n-e-1-i}$$

On considère les différents quotients suivants :

$$j_0(2n) = \frac{C_{2n}^{n}}{C_{2n}^{n-e-1}}, j_1(2n) = \frac{C_{2n}^{n-1}}{C_{2n}^{n-e-2}}, \dots, j_k(2n) = \frac{C_{2n}^{n-k}}{C_{2n}^{n-e-1-k}}, \dots, j_e(2n) = \frac{C_{2n}^{n-e}}{C_{2n}^{n-2e-1}} \qquad (8.28)$$

On commence par la détermination de la limite de $J_0(2n)$ à l' infini:

$$J_0(2n) = \frac{C_{2n}^{n}}{C_{2n}^{n-e-1}}$$

On sait que :

$$\begin{cases} C_{2n}^{n} = \dfrac{(2n)!}{n!\, n!} \\ C_{2n}^{n-e-1} = \dfrac{(2n)!}{(n-e-1)!\,(n+e+1)!} \end{cases}$$

Ce qui nous permet d'écrire:

$$J_0(2n) = \frac{C_{2n}^{n}}{C_{2n}^{n-e-1}} = \frac{(n+e+1)!\,(n-e-1)!}{n!\, x\, n!}$$

On remplace $(n+e+1)!$ et n! dans $J_0(2n)$ par les expressions suivantes:

$$\begin{cases}(n+e+1)! = (n+e+1)(n+e)x..x(n+1)x\, n! \\ n! = nx(n-1)(n-2)x..x(n-e)x\,(n-e-1)!\end{cases}$$

On peut écrire:

$$J_0(2n) = \frac{(n+e+1)(n+e)x..x(n+1)x\, n!}{n!} x \frac{(n-e-1)!}{nx(n-1)..x(n-e)x\,(n-e-1)!}$$

$$= \frac{(n+e+1)(n+e)\dots(n+2)(n+1)}{n(n-1)(n-2)\dots(n-e+1)(n-e)}$$

$$= \left(\frac{n+e+1}{n}\right) x \left(\frac{n+e}{n-1}\right) x \left(\frac{n+e-1}{n-2}\right) x \dots x \left(\frac{n+2}{n-e+1}\right) x \left(\frac{n+1}{n-e}\right)$$

$$= \left(1+\frac{e+1}{n}\right)\left(1+\frac{e+1}{n-1}\right) x \left(1+\frac{e+1}{n-2}\right) \dots x \left(1+\frac{e+1}{n-e+1}\right)\left(1+\frac{e+1}{n-e}\right) \quad (e\ 8.29)$$

La dernière expression obtenue contient (e+1) termes, de plus chaque terme vérifie l'inéquation suivante :

$$\left(1+\frac{e+1}{n-k}\right) \geq \left(1+\frac{e+1}{n}\right) \text{ pour tout entier k tel que } 0 \leq k \leq e \quad (e\ 8.30)$$

On peut déduire l'inéquation suivante :

$$\left(1+\frac{e+1}{n}\right) x \left(1+\frac{e+1}{n-1}\right) x \dots x \left(1+\frac{e+1}{n-e+1}\right) x \left(1+\frac{e+1}{n-e}\right) \geq \left(1+\frac{e+1}{n}\right)^{e+1} \quad (e\ 8.31)$$

En utilisant, l'inéquation suivante qui correspond à des valeurs de n suffisamment grandes:

$$0.27n > e > 0.25n$$

On déduire que :

$$e > (0.25n - 1)$$

Ce qui implique:

$$(e+1) > 0.25n$$

En déduire que:

$$\left(1+\frac{e+1}{n}\right) > 1.25$$

Ce qui implique:

$$\left(1+\frac{e+1}{n}\right)^{e+1} > (1{,}25)^{e+1} > (1{,}25)^{0{,}25n} \quad (e\ 8.32)$$

Revenons aux expressions de départ, en exploitant les relations ci-dessus, on peut écrire :

$$J_0(2n) = \frac{(n+e+1)!\,(n-e-1)!}{n!\,x\,n!} > \left(1+\frac{e+1}{n}\right)^{e+1} > (1,25)^{e+1} > (1,25)^{0,25n} \quad (e\ 8.33)$$

On déduit que:

$$\lim_{n\to+\infty} \frac{(n+e+1)!\,(n-e-1)!}{n!\,x\,n!} \geq \lim_{n\to+\infty} (1,25)^{0,25n} = +\infty \qquad (e\ 8.34)$$

On peut conclure que:

$$\lim_{n\to+\infty} \frac{(n+e+1)!\,(n-e-1)!}{n!\,x\,n!} = +\infty$$

Equivaut à:

$$\lim_{n\to+\infty} (J_0(2n)) = \lim_{n\to+\infty} \left(\frac{C_{2n}^{n}}{C_{2n}^{n-e-1}}\right) = +\infty \qquad (e\ 8.35)$$

Détermination de la limite de $J_k(2n)$ lorsque 2n tend vers l'infini pour tout entier k tel que

$1 \leq k \leq e$

Pour tout entier naturel k tel que $k < 2n$ on peut vérifier que :

$$\frac{C_{2n}^{k+1}}{C_{2n}^{k}} = \frac{2n-k}{k+1}$$

On applique cette relation aux deux entiers naturel non nuls q et P quelconque tel que $q < P < 2n$, on peut déterminer la différence suivante:

$$\frac{C_{2n}^{q+1}}{C_{2n}^{q}} - \frac{C_{2n}^{p+1}}{C_{2n}^{p}} = \frac{2n-q}{q+1} - \frac{2n-p}{p+1}$$

$$= \frac{2np+2n-pq-q}{(q+1)(P+1)} - \frac{2nq-pq+2n-p}{(q+1)(P+1)}$$

$$= \frac{2n(p-q)+(p-q)}{(q+1)(P+1)}$$

$$= \frac{(p-q)(2n+1)}{(q+1)(P+1)} \qquad (e\ 8.36)$$

Comme $P > q$ ( $p - q > 0$) on peut déduire que :

$$\frac{C_{2n}^{q+1}}{C_{2n}^{q}} - \frac{C_{2n}^{p+1}}{C_{2n}^{p}} > 0 \Longrightarrow \frac{C_{2n}^{q+1}}{C_{2n}^{q}} > \frac{C_{2n}^{p+1}}{C_{2n}^{p}}$$

Ce qui implique :

$$\frac{C_{2n}^{p}}{C_{2n}^{q}} > \frac{C_{2n}^{p+1}}{C_{2n}^{q+1}} \Longrightarrow \frac{C_{2n}^{p+1}}{C_{2n}^{q+1}} < \frac{C_{2n}^{p}}{C_{2n}^{q}} \qquad (e\ 8.37)$$

On peut la généraliser comme suit :

$$\frac{C_{2n}^{p+1}}{C_{2n}^{q+1}} < \frac{C_{2n}^{p}}{C_{2n}^{q}} < \frac{C_{2n}^{p-1}}{C_{2n}^{q-1}} < \frac{C_{2n}^{p-2}}{C_{2n}^{q-2}} < \cdots \qquad (e\ 8.38)$$

On applique cette inégalité aux différents quotients $J_k(2n)$ on peut déduire :

$$\frac{C_{2n}^{n-e}}{C_{2n}^{n-2e-1}} > \frac{C_{2n}^{n-e+1}}{C_{2n}^{n-2e}} > \cdots > \frac{C_{2n}^{n-2}}{C_{2n}^{n-e-3}} > \frac{C_{2n}^{n-1}}{C_{2n}^{n-e-2}} > \frac{C_{2n}^{n}}{C_{2n}^{n-e-1}} \qquad (e\ 8.39)$$

Les expressions de la limite relative à chaque terme lorsque n tend vers l'infini vérifie la relation ci-dessous :

$$\lim_{n\to+\infty} \frac{C_{2n}^{n-e}}{C_{2n}^{n-2e-1}} \geq \lim_{n\to+\infty} \frac{C_{2n}^{n-e+1}}{C_{2n}^{n-2e}} \geq \cdots \geq \lim_{n\to+\infty} \frac{C_{2n}^{n-1}}{C_{2n}^{n-e-2}} \geq \lim_{n\to+\infty} \frac{C_{2n}^{n}}{C_{2n}^{n-e-1}} \qquad (e\ 8.40)$$

Comme on a:

$$\lim_{n\to+\infty} \frac{C_{2n}^{n}}{C_{2n}^{n-e-1}} = +\infty$$

Ce qui nous permet de conclure que:

$$\lim_{n\to+\infty} \frac{C_{2n}^{n}}{C_{2n}^{n-m-1}} = \lim_{n\to+\infty} \frac{C_{2n}^{n-1}}{C_{2n}^{n-e-2}} = \lim_{n\to+\infty} \frac{C_{2n}^{n-2}}{C_{2n}^{n-e-3}} = \cdots = \lim_{n\to+\infty} \frac{C_{2n}^{n-e}}{C_{2n}^{n-2e-1}} = +\infty \qquad (e\ 8.41)$$

Détermination de la limite de $\frac{f(2n)}{g(2n)}$ lorsque n tend vers l'infini

On sait que :

$$\frac{f(2n)}{g(2n)} = \frac{C_{2n}^{n} + C_{2n}^{n-1} + \cdots + C_{2n}^{n-e}}{C_{2n}^{n-e-1} + C_{2n}^{n-e-2} + \cdots + C_{2n}^{n-2e-1}}$$

Première méthode:

On peut utiliser la méthode décrite dans la stratégie de la démonstration. Cette méthode déjà décrite dans le paragraphe consacrée pour la stratégie de cette démonstration.

Deuxième méthode:

On suppose que:

$$\lim_{n\to+\infty} \frac{f(2n)}{g(2n)} = L \Longrightarrow \lim_{n\to+\infty} \frac{C_{2n}^{n} + C_{2n}^{n-1} + \cdots + C_{2n}^{n-e}}{C_{2n}^{n-e-1} + C_{2n}^{n-e-2} + \cdots + C_{2n}^{n-2e-1}} = L \qquad (e\ 8.42)$$

L'équation ci-dessus s'écrit sous la forme suivante :

$$\lim_{n\to+\infty} [(C_{2n}^{n} + C_{2n}^{n-1} + \cdots + C_{2n}^{n-e}) - L(C_{2n}^{n-e-1} + C_{2n}^{n-e-2} + \cdots + C_{2n}^{n-2e-1})] = 0$$

On sait que:

$$C_{2n}^{n-k} = J_k(2n)\, C_{2n}^{n-e-k-1}$$

f(2n) peut s'exprimer comme suit :

$$f(2n) = \sum_{k=0}^{k=e} J_k(2n)\, C_{2n}^{n-e-k-1}$$

Remplaçons f(2n) par son expression ci-dessus  nous permet d'écrire :

$$\lim_{n\to+\infty}[\sum_{k=n-2e-1}^{n-e-1}(J_k(2n)C_{2n}^{k}) - L(\sum_{k=n-2e-1}^{n-e-1} C_{2n}^{k})] = 0$$

Ou encore:

$$\lim_{n\to+\infty}[\sum_{k=n-2e-1}^{n-e-1}((J_k(2n) - L)C_{2n}^{k})] = 0 \qquad \text{(e 8.44)}$$

On peut déduire que:

$$\lim_{n\to+\infty}((J_k(2n) - L)C_{2n}^{k}) = 0 \qquad \text{(e 8.45)}$$

Comme $C_{2n}^{k}$ tend vers l'infinie lorsque  n tend vers l'infini pour un entier non nul k donné, on déduire que :

$$\lim_{n\to+\infty}(J_k(2n) - L) = 0 \Longrightarrow \lim_{n\to+\infty} J_k(2n) = L \qquad \text{(e 8.46)}$$

Comme on a  montré que :

$$\lim_{n\to+\infty} J_k(2n) = +\infty$$

On peut conclure que:

$$\lim_{n\to+\infty}\frac{f(2n)}{g(2n)} = \lim_{n\to+\infty}\frac{C_{2n}^{n} + C_{2n}^{n-1} + \cdots + C_{2n}^{n-e}}{C_{2n}^{n-e-1} + C_{2n}^{n-e-2} + \cdots + C_{2n}^{n-2e-1}} = +\infty \qquad \text{(e 8.47)}$$

En déduit alors:

$$\lim_{n\to\infty}\frac{a_{2n}}{f(2n)} = 0 \qquad (8.48)$$

Comme on a:

$$\frac{a_{2n}}{f(2n)} > \frac{a_{2n}}{b_{2n}} > \frac{a_{2n}}{2^{2n}}$$

On peut déduire la limite à l' infini de $r_A(2n)$ comme suit

$$\lim_{n\to+\infty}\frac{a_{2n}}{f(2n)} = \lim_{n\to+\infty}\frac{a_{2n}}{2^{2n}} = 0 \qquad \text{(e 8.49)}$$

On peut déduire aussi la limite à l' infini de $r_B(2n)$ à partir de la relation suivante :

$$b_{2n} = 2^{2n} - a_{2n}$$

Ceci implique que :

$$\lim_{n\to+\infty}\left(\frac{b_{2n}}{2^{2n}}\right) = \lim_{n\to+\infty}\left(1 - \frac{a_{2n}}{2^{2n}}\right) = 1 \qquad \text{(e 8.50)}$$

On peut conclure que lorsque n tend vers l'infini, la proportion des suites générées de type $A_{2n}$ de longueurs 2n tend vers 0, alors que la proportion des suites générées de type $B_{2n}$ tend vers 1.

## 8.5 Décomposition parfaite d'une chromatrice de Collatz

La suite génératrice parfait est une suite arithmétique de raison 2, elle permet de créer une matrice constitue par un ensemble des suites de Collatz dont leurs structures constituent un tableau des arrangements binaires (ou chromatique) complets. En réalité, toutes suites arithmétiques de raison un entier naturel pair non nul quelconque peut générer des chromatrices parfaites du Collatz qui admettent des tableaux structurels qui renferment tous les arrangements binaires avec répétition.

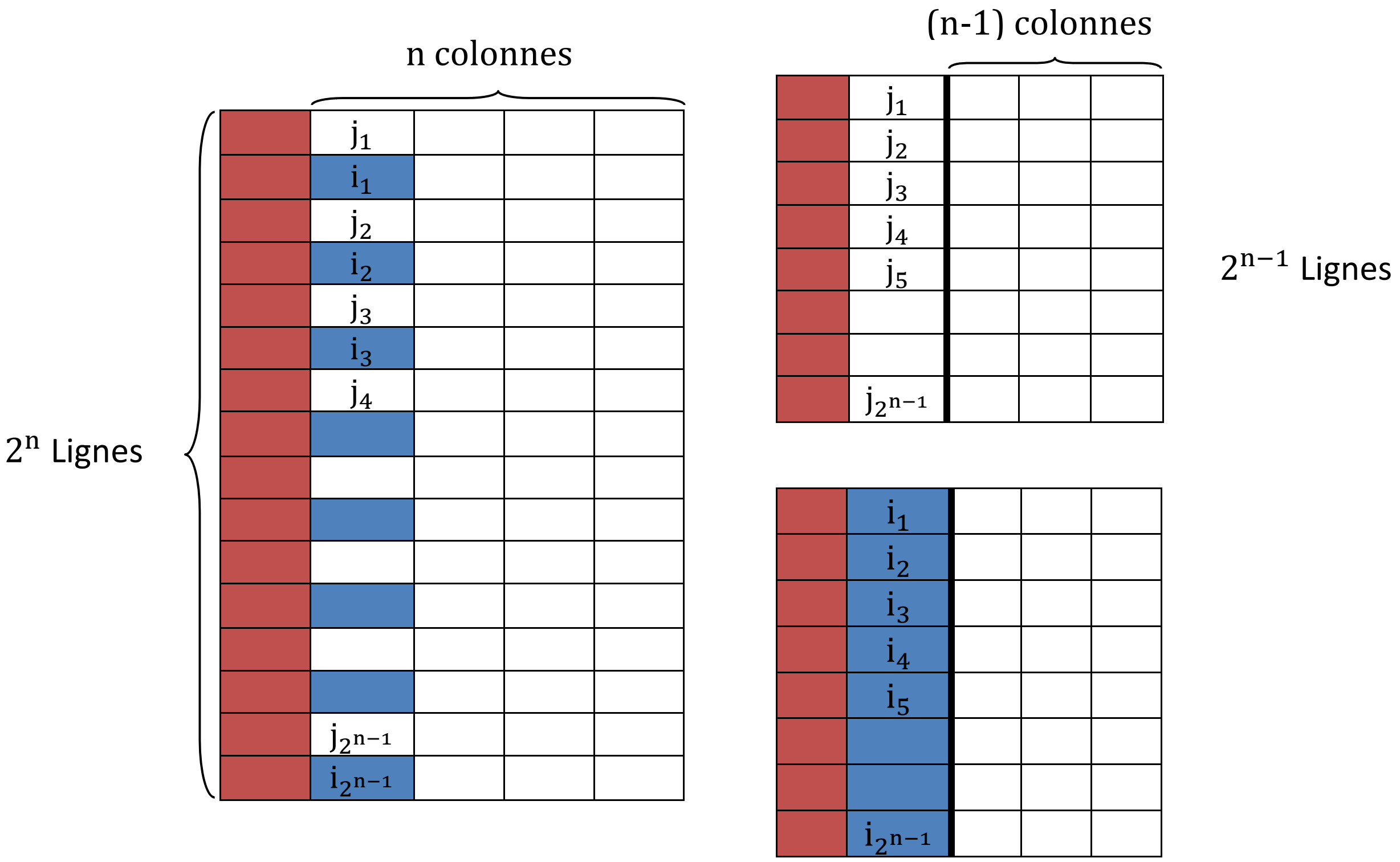

*Figure 20:* Principe de la décomposition parfaite d'une chromatrice de Collatz

Si on considère la chromatrice de Collatz d'ordre dimensionnel n. On sait que la première colonne (colonne supplémentaire ne fait pas partie de la chromatrice) renferme tous les termes de la suite génératrice de ce tableau. La deuxième colonne contient des entiers pairs qui constituent une suite arithmétique de raison un entier pair non nul aussi les termes impairs dans cette même colonne constituent une suite arithmétique de raison un entier pair non nul. La décomposition parfaite consiste a construire deux matrices (n-1 colonnes, $2^{n-1}$ lignes) à partir de la matrice de départ tel que la deuxième colonne dans l'une de ces deux matrice est constituée uniquement des entiers impairs alors que la deuxième colonne de l'autre matrice constituée uniquement des entier pairs.

**Exemple 8.5**

La figure suivante montre un exemple d'une décomposition parfaite de la chromatrice $\mathbb{M}^{T}(1,4)$ :

| | | | | |
|---|---|---|---|---|
| 1 | 2 | 1 | 2 | 1 |
| 3 | 5 | 8 | 4 | 2 |
| 5 | 8 | 4 | 2 | 1 |
| 7 | 11 | 17 | 26 | 13 |
| 9 | 14 | 7 | 11 | 17 |
| 11 | 17 | 26 | 13 | 20 |
| 13 | 20 | 10 | 5 | 8 |
| 15 | 23 | 35 | 53 | 80 |
| 17 | 26 | 13 | 20 | 10 |
| 19 | 29 | 44 | 22 | 11 |
| 21 | 32 | 16 | 8 | 4 |
| 23 | 35 | 53 | 80 | 40 |
| 25 | 38 | 19 | 29 | 44 |
| 27 | 41 | 62 | 31 | 47 |
| 29 | 44 | 22 | 11 | 17 |
| 31 | 47 | 71 | 107 | 161 |
| 1 | 2 | 1 | 2 | 1 |
| 5 | 8 | 4 | 2 | 1 |
| 9 | 14 | 7 | 11 | 17 |
| 13 | 20 | 10 | 5 | 8 |
| 17 | 26 | 13 | 20 | 10 |
| 21 | 32 | 16 | 8 | 4 |
| 25 | 38 | 19 | 29 | 44 |
| 29 | 44 | 22 | 11 | 17 |
| 3 | 5 | 8 | 4 | 2 |
| 7 | 11 | 17 | 26 | 13 |
| 11 | 17 | 26 | 13 | 20 |
| 15 | 23 | 35 | 53 | 80 |
| 19 | 29 | 44 | 22 | 11 |
| 23 | 35 | 53 | 80 | 40 |
| 27 | 41 | 62 | 31 | 47 |
| 31 | 47 | 71 | 107 | 161 |

*Figure 21:* Décomposition parfaite de la chromatrice $\mathbb{M}^{T}(1,4)$ en deux chromatrices

On remarque les deux chromatrices contiennent tous les arrangements chromatiques de deux couleurs (blanc et bleu) cette propriété est générale quelque soit l'ordre dimensionnel de la chromatrice de départ.

## 8.6 Les chromoformes parfaits infinis de Collatz

**Définition 8.3**

Le chromoforme parfait infini de Collatz correspond au chromoforme parfait d'ordre dimensionnel infini, autrement il correspond à la limite du chromoforme fini d'ordre dimensionnel n lorsque n tend vers l'infini:

$$\mathbb{Z}^{T}(1,\dots) = \lim_{n\to+\infty} \mathbb{Z}^{T}(1,n) \qquad \text{(e 8.51)}$$

Dans ce cas on obtient un nombre infini des suites générées de longueurs aussi infinies.

Le chromoforme parfait infini dont la suite génératrice est la suite constituée par tous les entiers impairs est comme suit :

$$\mathbb{Z}^{T}(1,\dots) = \left(\tilde{S}^{R}(1), \tilde{S}^{R}(2), \tilde{S}^{r}(5), \tilde{S}^{R}(7), \dots., \tilde{S}^{R}(2k+1), \dots\right) \qquad \text{(e 8.52)}$$

**Définition 8.4**

Aussi la limite de $\mathbb{Z}^{T}(2,n)$ lorsque n tend vers l'infini correspond au chromoforme parfait infini dont la suite génératrice est l'ensemble constitué par tous les entiers naturel pairs non nul :

$$\mathbb{Z}^{T}(2,\dots) = \lim_{n\to+\infty} \mathbb{Z}^{T}(2,n) \qquad \text{(e 8.53)}$$

Autrement :

$$\mathbb{Z}^{T}(2,\dots) = \left(\tilde{S}^{R}(2), \tilde{S}^{R}(4), \tilde{S}^{R}(6), \dots., \tilde{S}^{R}(2k), \dots\right) \qquad \text{(e 8.54)}$$

**Remarque 8.3**

Les représentations matricielles de ces deux chromoformes infinis correspondent à deux matrices des dimensions infinies. Chaque matrice contient une infinité des colonnes et une infinité des lignes. La première colonne contient tous les entiers naturels impairs alors que pour la première colonne de la deuxième chromatrice, contient tous les entiers naturels pairs non nuls.

**Notations 8.2**

Rappeler que :

$\tilde{L}^{b}(P,n)$ désigne une distribution structurelle binaire finie relative à une suite générée finie.

$\tilde{L}^{c}(P,n)$ désigne une distribution structurelle chromatique finie relative à une suite générée finie.

On désigne par :

$\tilde{L}^{B}(P)$ La distribution structurelle binaire infinie relative à une suite générée infinie.

$\tilde{L}^{Cr}(P)$ La distribution structurelle chromatique infinie relative à une suite générée infinie.

$\tilde{L}^{S}(P)$ La distribution structurelle de longueur infinie (chromatique ou binaire) relative à une suite générée de longueur infinie notée $\tilde{S}^{R}(P)$, alors que la notation $\tilde{L}^{s}(P,n)$ avec petite s est relative à une suite générée de longueur finie.

**Définition 8.5**

On fait correspondre à chaque chromoforme parfait de Collatz un chromoforme structurel parfait d'ordre dimensionnel infini ce chromoforme structurel infini est notée $\mathbb{Z}^{S}(1,\dots)$ et il contient toutes les distributions structurelles linéaires relatives aux toutes suites générées de longueurs infinies du chromoforme parfait infini $\mathbb{Z}^{T}(1,\dots)$:

$$\mathbb{Z}^{S}(1,\dots) = \left(\tilde{L}^{S}(1), \tilde{L}^{S}(5), \tilde{L}^{S}(7), \dots., \tilde{L}^{S}(2k+1), \dots\right) \qquad \text{(e 8.55)}$$

**Définition 8.6**

De même on peut définir le chromoforme structurel infini relatif au chromoforme parfait infini de Collatz $\mathbb{Z}^{T}(2,\dots)$ comme suit :

$$\mathbb{Z}^{S}(2,\dots) = \left(\tilde{L}^{S}(2), \tilde{L}^{S}(4), \tilde{L}^{S}(8), \dots., \tilde{L}^{S}(2k), \dots\right) \qquad \text{(e 8.56)}$$

Les deux chromoformes parfaits d'ordre dimensionnel infini peuvent être regroupés dans un seul ensemble qu'on le note $CLZ^{\infty}$ c'est l'ensemble constitué par toutes les suites générées de longueurs infinies il est comme suit :

$$CLZ^{R} = \left(\tilde{S}^{R}(1), \tilde{S}^{R}(2), \tilde{S}^{R}(3), \dots., \tilde{S}^{R}(P), \dots\right) \qquad \text{(e 8.57)}$$

La suite génératrice de l'ensemble $CLZ^{\infty}$ correspond donc à l'ensemble des entiers naturels non nuls $\mathbb{N}^{*}$.

**Définition 8.7**

Le chromorfisme est l'ensemble constitué par toutes les distributions structurelles linéaires de longueurs infinies qui correspondent aux toutes les suites générées infinies de l'ensemble $CLZ^{\infty}$. Il est noté comme suit:

$$CLZ^{S} = \left(\tilde{L}^{S}(1), \tilde{L}^{S}(2), \tilde{L}^{S}(3), \dots, \tilde{L}^{S}(P), \dots.\right) \qquad \text{(e 8.58)}$$

**Théorème 8.3**

Toutes les suites de Collatz appartenant à un même chromoforme parfait infini ($\mathbb{Z}^{T}(1, \dots)$ ou bien $\mathbb{Z}^{T}(2, \dots)$ ) sont deux à deux non isoformes c'est à dire qu'ils n'existent pas deux entiers naturels non nuls qui ont la même parité et qui peuvent générer deux suites de Collatz de longueurs infinies possédants la même distribution structurelle ceci se traduit par :

$\forall$ entiers naturels non nuls $P_1$ et $P_2$ tel que $i_0(P_1) = i_0(P_2)$ :

$$\text{si } P_1 \neq P_2 \text{ alors } \tilde{L}^{S}(P_1) \neq \tilde{L}^{S}(P_2) \qquad \text{(e 8.59)}$$

**Remarque 8.4**

Ce théorème appelé aussi la propriété de la singularité structurelle relative aux suites de Collatz.

**Démonstration**

On suppose qu'elles existent deux suites générées et de longueurs infinies isochromatiques générées respectivement par deux entiers naturels non nul $P_1$ et $P_2$ notées respectivement $\tilde{S}^{R}(P_1)$ et $\tilde{S}^{R}(P_2)$ .

Comme $\tilde{S}^{R}(P_1)$ et $\tilde{S}^{R}(P_2)$ sont isoformes donc pour tout entier naturel non nul n les deux suites générées finies $\tilde{S}^{r}(P_1, n)$ et $\tilde{S}^{r}(P_2, n)$ sont n-isoformes.

De plus comme on donc il existe un entier naturel non nul $n_0$

$$\sup\{P_1, P_2\} < 2^{n+1} \qquad \text{(e 8.60)}$$

Dans ce cas les deux générateurs $P_1$ et $P_2$ appartenant à la même suite génératrice parfaite $Y^{S}(P, n)$ avec P peut être égal à 1 ou bien 2 selon la parité de ces deux entiers. De même dans ce cas les deux suites $\tilde{S}^{r}(P_1, n)$ et $\tilde{S}^{r}(P_2, n)$ appartient à la même chromatrice $\mathbb{M}^{T}(1, n)$ ceci est impossible puisque toutes chromatrices parfaites de Collatz ne peuvent pas contenir deux suites isochromatiques (ou isoformes) ce qui nous permet de conclure qu'il est impossible que deux entiers naturels différents quelconques peuvent générer deux suites isoformes de longueurs infinies.

# 9 Les polychromoformes parfaits de Collatz

**Définition 9.1**

Un polychromoforme parfait de Collatz est l'ensemble constitué par tous les chromoformes parfaits construits à partir d'un même système générateur parfait.

## 9.1 Construction d'un polychromoforme de Collatz

On considère le chromoforme parfait de Collatz construit à partir de la suite génératrice parfaite $Y^S(P_{i,1}, n)$ comme suit:

$$Z^T(P_{i,1}, n) = [\tilde{S}^r(P_{i,1}, n), \tilde{S}^r(P_{i,2}, n), \tilde{S}^r(P_{i,3}, n), \tilde{S}^r(P_{i,4}, n), \dots, \tilde{S}^r(P_{i,2^n}, n)] \quad (e\ 9.1)$$

On sait qu'on a une infinité des suites génératrices parfaites dans un système générateur parfait donc on peut construire une infinité des chromoformes parfaits de Collatz.

$$Sy^P(P_{1,1}, n) = [Y^S(P_{1,1}, n), Y^S(P_{2,1}, n), Y^S(P_{3,1}, n), , \dots, Y^S(P_{i,1}, n), \dots.] \quad (e\ 92)$$

Avec $P_{k,1} = (P_{1,1} + 2^{n+1}k)$ le premier terme de la suite génératrice parfaite du rang k dans le système générateur considéré.

On fait correspondre à chaque suite génératrice parfaite du système générateur parfait son chromoforme parfait de Collatz, on obtient le polychromoforme parfait du Collatz. Il est noté comme suit ;

$$\mathbb{Z}^{\mathbf{T}}(P_{1,1}, n) = [Z^T(P_{1,1}, n), Z^T(P_{2,1}, n), Z^T(P_{3,1}, n), \dots, . Z^T(P_{i,1}, n), \dots] \quad (e\ 9.3)$$

-Les entiers $P_{1,1}, n$ représentent les paramètres constructifs du polychromoforme.

-n est appelé l'ordre dimensionnel du polychromoforme considéré.

**Définition 9.2**

Dans l'ensemble on peut définir deux polychromoformes dans $\mathbb{N}^*$ qui sont les suivants :

Le polychromoforme construit à partir du système générateur dont les termes sont des entiers impairs qui correspond au système générateur $Sy^P(1, n)$. Le polychromoforme généré par ce système est le suivant:

$$\mathbb{Z}^{\mathbf{T}}(1, n) = [Z^T(1, n), Z^T(1 + 2^{n+1}, n), Z^T(1 + 2x2^{n+1}, n), \dots, . Z^T(1 + 2^{n+1}k, n), \dots.] \quad (e\ 9.4)$$

Le polychromoforme générés par le système générateur $Sy^P(2, n)$, il est comme suit :

$$\mathbb{Z}^{\mathbf{T}}(2, n) = [Z^T(2, n), Z^T(2 + 2^{n+1}, n), Z^T(2 + 2x2^{n+1}, n), \dots, . Z^T(2 + 2^{n+1}k, n), \dots.] \quad (e\ 9.5)$$

## 9.2 Les représentations matricielles de polychromoforme du Collatz

La représentation matricielle condensée du polychromoforme d'ordre dimensionnel n correspond à un tableau constitué de $2^n$ colonnes et d'une infinité des lignes. Chaque ligne du tableau correspond à un chromoforme parfait de Collatz bien déterminé. Noter que les cases contiennent uniquement les symboles relatifs aux suites générées de Collatz constituants le polychromoforme. On ajoute une ligne qui indique la nature des distributions verticales (la suite constituée par les termes d'une même colonne) et une colonne qui nous indique la nature des distributions horizontales (la suite constituée par les termes d'une même ligne). Noter que la ligne et la colonne supplémentaires ne font pas partie du tableau, elles sont séparées du reste du tableau par deux traits en gras. La matrice principale ou fondamentale relative au polychromoforme considéré est notée $\mathbb{G}^{M}(P_{1,1},n)$.

*Tableau7:* représentation matricielle du polychromoforme $\mathbb{Z}^{T}(P_{1,1},n)$:

| | $H^{T}(P_{1,1},n)$ | $H^{T}(P_{1,2},n)$ | $H^{T}(P_{1,3},n)$ | | | $H^{T}(P_{1,2^n},n)$ |
|---|---|---|---|---|---|---|
| $Z^{T}(P_{1,1},n)$ | $\tilde{S}^{r}(P_{1,1},n)$ | $\tilde{S}^{r}(P_{1,2},n)$ | $\tilde{S}^{r}(P_{1,3},n)$ | | | $\tilde{S}^{r}(P_{1,2^n},n)$ |
| $Z^{T}(P_{2,1},n)$ | $\tilde{S}^{r}(P_{2,1},n)$ | $\tilde{S}^{r}(P_{2,2},n)$ | $\tilde{S}^{r}(P_{2,3},n)$ | | | $\tilde{S}^{r}(P_{2,2^n},n)$ |
| $Z^{T}(P_{3,1},n)$ | $\tilde{S}^{r}(P_{3,1},n)$ | $\tilde{S}^{r}(P_{3,2},n)$ | $\tilde{S}^{r}(N_{3,3},n)$ | | | $\tilde{S}^{r}(P_{3,2^n},n)$ |
| | | | | | | |
| | | | | | | |
| $Z^{T}(P_{i,1},n)$ | $\tilde{S}^{r}(P_{i,1},n)$ | $\tilde{S}^{r}(P_{i,2},n)$ | $\tilde{S}^{r}(P_{i,3},n)$ | | | $\tilde{S}^{r}(P_{i,2^n},n)$ |
| | | | | | | |

Les distributions verticales constituées par une infinité des suites qui sont sous la forme suivante:

$$(\tilde{S}^{r}(P_{1,j},n),\tilde{S}^{r}(P_{2,j},n),\tilde{S}^{r}(N_{3,j},n),\ldots.,\tilde{S}^{r}(P_{i,j},n),\ldots)$$

Avec j un entier qui vérifie : $1 \leq j \leq 2^{n}$

**Lemme 9.1**

On sait que cette distribution correspond au chromologue $H^{T}(P_{1,j},n)$, on peut conclure que chaque colonne de la polymatrice $\mathbb{G}^{M}(P_{1,1},n)$ du polychromoforme $\mathbb{Z}^{\mathbf{T}}(P_{1,1},n)$ correspond à un chromologue bien déterminé du polychromologue $\mathbb{H}^{T}(P_{1,1},n)$ ce qui nous permet de conclure que $\mathbb{G}^{M}(P_{1,1},n)$ est la transposée de la matrice $\mathbb{H}^{M}(P_{1,1},n)$ de polychromologue $\mathbb{H}^{T}(P_{1,1},n)$ ceci se traduit par l'expression suivante :

$$\mathbb{G}^{M}(P_{1,1},n) = \left(\mathbb{H}^{M}(P_{1,1},n)\right)^{\text{Transp}} \qquad \text{(e 9.6)}$$

La polymatrice contient tous les chromologues dans leurs colonnes, elle renferme donc $2^{n}$ chromologues appartenant à un même polychromologue qui est construit à partir d'un même super-système.

**Lemme 9.2**

Chaque suite du premier chromoforme $Z^T(P_{1,1}, n)$ du polychromoforme $\mathbb{Z}^{\mathbf{T}}(P_{1,1}, n)$ représente la première suite d'un chromologue donné de polychromologue $\mathbb{H}^T(P_{11}, n)$. Plus précisément la suite $\tilde{S}^r(P_{k,1}, n)$ du rang k du chromoforme $\mathbb{Z}^{\mathbf{T}}(P_{1,1}, n)$ représente la première suite du chromolgue $H^T(P_{1,k}, n)$ du rang k dans le polychromologue $\mathbb{H}^T(P_{1,1}, n)$.

**Démonstration**

Comme la matrice du polychromoforme est la transposée de la matrice du polychromologue et chaque colonne représente un chromologue donné donc la suite

## 9.3 La matrice structurelle d'un polychromoforme

Elle correspond à la représentation matricielle du polychromoforme structurel $\mathbb{Z}^S(P_{1,1}, n)$ relatif au polychromoforme $\mathbb{Z}^T(P_{1,1}, n)$. Elle est notée $\mathbb{G}^S(P_{1,1}, n)$

*Tableau 8:* Matrice structurelle du polychromoforme $\mathbb{Z}^S(P_{1,1}, n)$

| | $H^S(P_{1,1}, n)$ | $H^S(P_{1,2}, n)$ | $H^S(P_{1,3}, n)$ | | | $H^S(P_{1,2^n}, n)$ |
|---|---|---|---|---|---|---|
| $Z^S(P_{1,1}, n)$ | $L^{\mathcal{S}}(P_{1,1}, n)$ | $L^{\mathcal{S}}(P_{1,2}, n)$ | $L^{\mathcal{S}}(P_{1,3}, n)$ | | | $L^{\mathcal{S}}(P_{1,2^n}, n)$ |
| $Z^S(P_{2,1}, n)$ | $L^{\mathcal{S}}(P_{2,1}, n)$ | $L^{\mathcal{S}}(P_{2,2}, n)$ | $L^{\mathcal{S}}(P_{2,3}, n)$ | | | $L^{\mathcal{S}}(P_{2,2^n}, n)$ |
| $Z^S(P_{3,1}, n)$ | $L^{\mathcal{S}}(P_{3,1}, n)$ | $L^{\mathcal{S}}(P_{3,2}, n)$ | $L^{\mathcal{S}}(N_{3,3}, n)$ | | | $L^{\mathcal{S}}(P_{3,2^n}, n)$ |
| | | | | | | |
| | | | | | | |
| $Z^S(P_{i,1}, n)$ | $L^{\mathcal{S}}(P_{i,1}, n)$ | $L^{\mathcal{S}}(P_{i,2}, n)$ | $L^{\mathcal{S}}(P_{i,3}, n)$ | | | $L^{\mathcal{S}}(P_{i,2^n}, n)$ |
| | | | | | | |

## 9.4 Représentation matricielle à trois niveaux constructifs de la polychromatrice du polychromoforme

Cette représentation matricielle fait apparaitre les différentes valeurs des différents termes de chaque suite. La polymatrice est constituée d'une infinité des chromatrices disposées l'une au dessous de l'autre tel que la matrice du rang k dans la polymatrice, en allant du haut vers le bas correspond à la matrice du chromoforme $Z^T(P_{k,1})$ du rang k de polychromoforme considéré.

**Définition 9.3**

La polychromatrice est un assemblage vertical des toutes les chromatrices relatives aux chromoformes constituants le polychromoformes, elle est notée comme suit :

$$\mathbb{M}^{TT}(P_{1,1},n) = Asv[\mathbb{M}^{T}(P_{1,1},n), \mathbb{M}^{T}(P_{2,1},n), \mathbb{M}^{T}(P_{3,1},n), \dots, \mathbb{M}^{T}(P_{i,1},n), \dots] \quad \text{(e 9.7)}$$

Dans la matrice $\mathbb{G}^{M}(P_{1,1},n)$ chaque chromoforme est représenté par une ligne alors que dans la matrice $\mathbb{M}^{TT}(P_{1,1},n)$ chaque chromoforme est représenté par une matrice.

La figure ci-dessous représente les trois première chromatrices de la polychromatrice du polychromoforme.

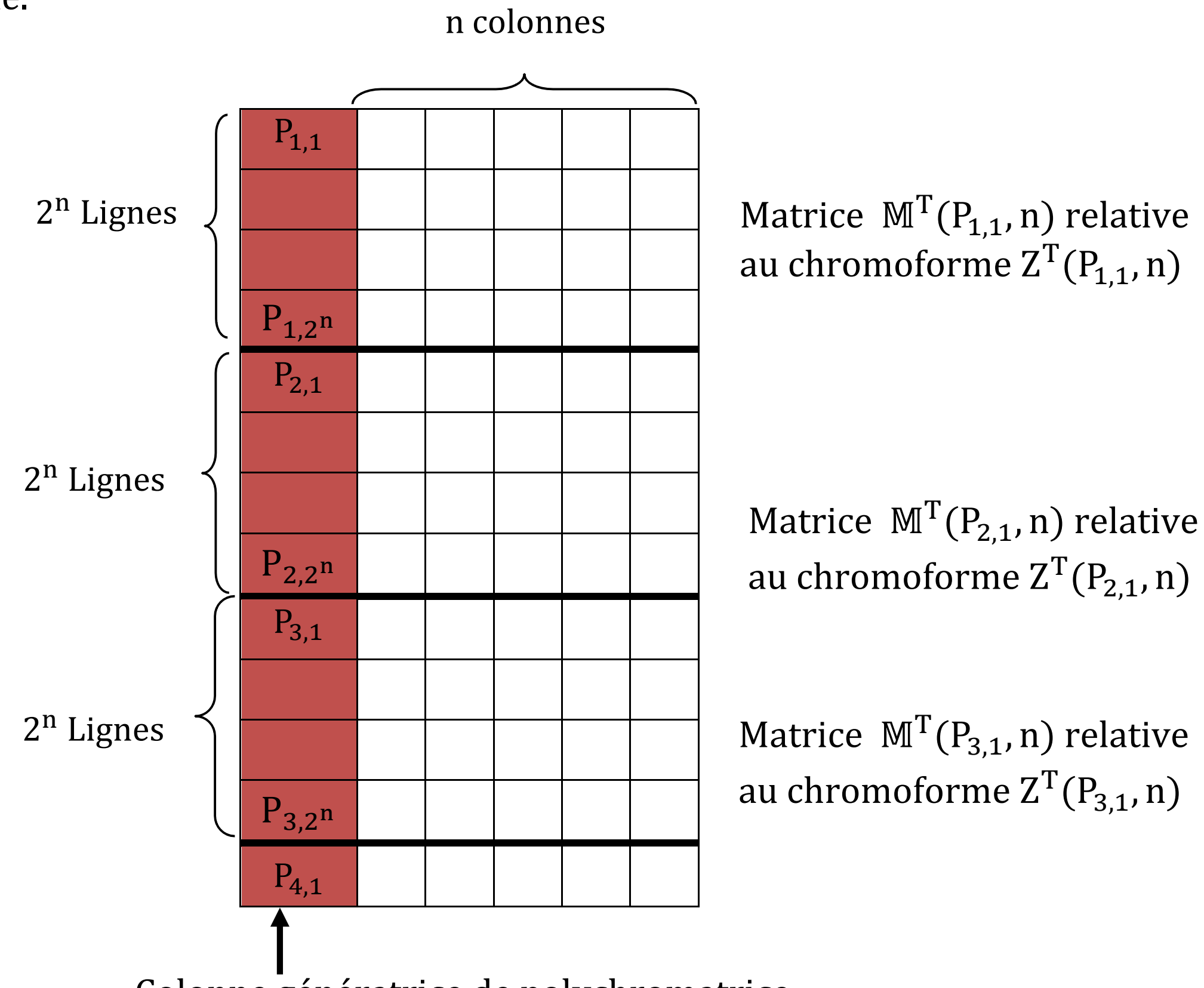

*Figure 22:* Assemblage vertical des différentes chromatrices dans la polychromatrice d'un polychromoforme parfait d'ordre dimensionnel n

- Toute les chromatrices possèdent le même nombre des colonnes et le même nombre des lignes. Chaque chromoforme correspond à la matrice fondamentale d'un chromoforme donné appartenant au polychromoforme globale.
- La première colonne en rouge brique est la colonne génératrice contient toutes les suites génératrice parfaites de l'un de deux système générateurs parfaits ainsi elle correspond à une suite arithmétique de raison 2. Elle contient tous les entiers impairs si le terme $P_{1,1}$correspond à 1 ou bien tous les entiers naturels non nuls pairs si le terme $P_{1,1}$ égal à 2.

**Exemple 9.1**

La figure suivante correspond à la polychromatrice $\mathbb{M}^{TT}(1,3)$ :

| | | | |
|---|---|---|---|
| 1 | 2 | 1 | 2 |
| 3 | 5 | 8 | 4 |
| 5 | 8 | 4 | 2 |
| 7 | 11 | 17 | 26 |
| 9 | 14 | 7 | 11 |
| 11 | 17 | 26 | 13 |
| 13 | 20 | 10 | 5 |
| 15 | 23 | 35 | 53 |
| 17 | 26 | 13 | 20 |
| 19 | 29 | 44 | 22 |
| 21 | 32 | 16 | 8 |
| 23 | 35 | 53 | 80 |
| 25 | 38 | 19 | 29 |
| 27 | 41 | 62 | 31 |
| 29 | 44 | 22 | 11 |
| 31 | 47 | 71 | 107 |
| 33 | 50 | 25 | 38 |
| 35 | 53 | 80 | 40 |
| 37 | 56 | 28 | 14 |
| 39 | 59 | 89 | 134 |
| 41 | 62 | 31 | 47 |
| 43 | 65 | 98 | 49 |
| 45 | 68 | 34 | 17 |
| 47 | 71 | 107 | 161 |
| 49 | 74 | 37 | 56 |
| 51 | 77 | 116 | 58 |
| 53 | 80 | 40 | 20 |
| 55 | 83 | 125 | 188 |
| 57 | 86 | 43 | 65 |
| 59 | 89 | 134 | 67 |
| 61 | 92 | 46 | 23 |
| 63 | 95 | 143 | 215 |

Suites isoformes d'un même chromologue

$\mathbb{M}^{T}(1,3)$

$\mathbb{M}^{T}(17,3)$

$\mathbb{M}^{T}(33,3)$

$\mathbb{M}^{T}(47,3)$

*Figure 23:* les quatre premières chromatrices de la polymatrice $\mathbb{M}^{TT}(1,3)$

On remarque bien que ces quatre chromatrices sont parfaitement isoformes elles possèdent la même conversion chromatique obtenue par coloration de chaque case qui contient un entier impair par une coloration bleu. En réalité, on peut montrer que toutes les chromatrices constituants cette polychromatrice sont isoformes. De plus deux générateurs de deux suites isochromatiques consécutives appartenant à cette polymatrice admettent une différence $2^4$. Puisque les suites d'une même chromatrice sont deux à deux non informes donc chaque suite appartient à un chromologue bien déterminé différent de tous les autres chromologues relatives aux autres suites de cette même chromatrice.

**Définition 9.4**

La matrice structurelle correspond au polychromoforme structurel est définie comme étant l'assemblage vertical de toutes les matrices structurelles relatives aux différentes chromoformes structurels. Elle est notée comme suit :

$$\mathbb{M}^{SS}(P_{1,1},n) = \mathrm{Asv}\left[\mathbb{M}^{S}(P_{1,1},n),\mathbb{M}^{S}(P_{2,1},n),\mathbb{M}^{S}(P_{3,1},n),\dots,\mathbb{M}^{S}(P_{i,1},n),\dots\right] \qquad (e\ 9.8)$$

Asv: assemblage vertical c'est à dire disposition des matrices l'une au dessous de l'autre.

**Théorème 9.1**

Tous les chromoformes d'un même polychromoforme sont parfaitement isochromatiques autrement toutes les matrices structurelles relatives aux tous les chromoformes d'un même polychromoforme sont identiques ceci se traduit par :

$$\mathbb{M}^{S}(P_{1,1},n) = \mathbb{M}^{S}(P_{2,1},n) = \mathbb{M}^{S}(P_{3,1},n) = \cdots = \mathbb{M}^{S}(P_{k,1},n) = \cdots \qquad (e\ 9.9)$$

Autrement, toutes les suites générées de deux chromoformes d'un même polychromoforme sont deux à deux isoformes.

**Démonstration**

Considérons deux chromoformes quelconques $Z^T(P_{i,1},n)$ et $Z^T(P_{j,1},n)$du polychromoforme $\mathbb{Z}^T(P_{1,1},n)$ définies comme suit :

$$Z^T(P_{i,1},n) = [\tilde{S}^r(P_{i,1},n),\tilde{S}^r(P_{i,2},n),\dots,\tilde{S}^r(P_{i,2^n},n)]$$
$$Z^T(P_{j,1},n) = [\tilde{S}^r(P_{j,1},n),\tilde{S}^r(P_{j,2},n),\dots,\tilde{S}^r(P_{j,2^n},n)]$$

Soient $P_{i,k}$ et $P_{j,k}$ les termes du rang k respectivement dans la suite génératrice parfaite $Y^S(P_{i,1},n)$ et la suite génératrice parfaite $Y^S(P_{j,1},n)$ on suppose de plus que $P_{j,k} > P_{i,k}$ donc on peut écrire pour tout entier k tel que $1 \leq k \leq 2^n$ :

$$\begin{cases} P_{i,k} = P_{i,1} + 2(k-1) \\ P_{j,k} = P_{j,1} + 2(k-1) \end{cases}$$

On déduit que :

$$P_{j,k} - P_{i,k} = P_{j,1} - P_{i,1}$$

Comme $P_{j,1}$ et $P_{i,1}$ sont deux termes de la même super-suite $X^S(P_{1,1}, n)$ donc on peut écrire :

$$P_{j,1} - P_{i,1} = 2^{n+1}(j - i) \qquad (e\ 9.10)$$

On déduit que:

$$P_{j,k} = P_{i,k} + 2^{n+1}(j - i)$$

On peut conclure que les deux suites générées $\tilde{S}^r(P_{j,k}, n)$ et $\tilde{S}^r(P_{i,k}, n)$ sont isoformes. Cette propriétés est vraie pour tout entier naturel k compris entre 1 et $2^n$ donc les suites de deux chromoformes sont deux à deux isochromatiques ce qui signifie que les deux chromoformes sont isoformes.

**Exemple 9.2**

Toutes les chromatrices de la polychromatrice $\mathbb{M}^{TT}(1,3)$ admettent la même matrice structurelle représentée ci-dessous :

$Y^S(P_{1,k}, 3)$

| | | | |
|---|---|---|---|
| | 0 | 1 | 0 |
| | 1 | 0 | 0 |
| | 0 | 0 | 0 |
| | 1 | 1 | 0 |
| | 0 | 1 | 1 |
| | 1 | 0 | 1 |
| | 0 | 0 | 1 |
| | 1 | 1 | 1 |

*Figure 24 :* chromatrice structurelle relative aux toutes les chromatrices de $\mathbb{M}^{TT}(1,3)$

Pour tout entier naturel non nul k, les suite génératrices parfaites $Y^S(P_{1,k}, 3)$ génèrent la même chromatrice structurelle. On dit que les différentes chromatrices $\mathbb{M}^T(P_{1,k}, 3)$ sont isoformes (ou isochromatiques).

## 9.5 La super-décomposition d'une polychromatrice parfaite de Collatz

La super-décomposition consiste à décomposer un polychromoforme donné d'ordre dimensionnel égal à n en $2^n$ chromologues.

On sait que chaque colonne de la matrice condensée de $\mathbb{G}^T(P_{1,1}, n)$ relative au polychromoforme $\mathbb{Z}^T(P_{1,1}, n)$ constitue un chromologue bien déterminé appartenant au polychromologue $\mathbb{H}^T(P_{1,1}, n)$. Comme la matrice $\mathbb{G}^T(P_{1,1}, n)$ renferme $2^n$ colonne donc on peut conclure qu on peut décomposer le polychromoforme $\mathbb{Z}^{T}(P_{1,1}, n)$ en $2^{n}$ chromologues. L'ensemble de ces chromologues obtenus par la décomposition $\mathbb{Z}^{T}(P_{1,1}, n)$constituent le polychromologue $\mathbb{H}^{T}(P_{1,1}, n)$. Les différents chromologues obtenus par la décomposition du polychromoforme $\mathbb{Z}^{T}(P_{1,1}, n)$ sont les suivants :

$$\mathbb{Z}^{T}(P_{1,1}, n) \rightarrow \underbrace{\left[H^{T}(P_{1,1}, n), H^{T}(P_{1,2}, n), H^{T}(P_{1,3}, n), \ldots, H^{T}(P_{1,2^{n}}, n)\right]}_{\mathbb{H}^{T}(P_{1,1}, n)} \qquad \text{(e 9.11)}$$

**Exemple 9.3**

Reprenons le cas de la polychromatrice $\mathbb{M}^{TT}(1,3)$ de la polychromoforme $\mathbb{Z}^{T}(1,3)$, une super-décomposition de cette matrice en des chromologues nous permet d'obtenir $2^{3}$ chromologues differentes. Chaque suite de la première chromatrice constitue la première suite d'un chromologue donné. La figure ci-dessous représente les différents chromologues obtenus par super-décomposition de la polymatrice $\mathbb{M}^{TT}(1,3)$.

| $H^{T}(1,3)$ | | | | $H^{T}(3,3)$ | | | | $H^{T}(5,3)$ | | | | $H^{T}(7,3)$ | | | |
|---|---|---|---|---|---|---|---|---|---|---|---|---|---|---|---|
| 1 | 2 | 1 | 2 | 3 | 5 | 8 | 4 | 5 | 8 | 4 | 2 | 7 | 11 | 17 | 26 |
| 17 | 26 | 13 | 20 | 19 | 29 | 44 | 22 | 21 | 32 | 16 | 8 | 23 | 35 | 53 | 80 |
| 33 | 50 | 25 | 38 | 35 | 53 | 80 | 40 | 37 | 56 | 28 | 14 | 39 | 59 | 89 | 134 |
| 49 | 74 | 37 | 56 | 51 | 77 | 116 | 58 | 53 | 80 | 40 | 20 | 55 | 83 | 125 | 188 |
| 65 | 98 | 49 | 74 | 67 | 101 | 152 | 76 | 69 | 104 | 52 | 26 | 71 | 107 | 161 | 242 |
| | | | | | | | | | | | | | | | |

| $H^{T}(9,3)$ | | | | $H^{T}(11,3)$ | | | | $H^{T}(13,3)$ | | | | $H^{T}(15,3)$ | | | |
|---|---|---|---|---|---|---|---|---|---|---|---|---|---|---|---|
| 9 | 14 | 7 | 11 | 11 | 17 | 26 | 13 | 13 | 20 | 10 | 5 | 15 | 23 | 35 | 53 |
| 25 | 38 | 19 | 29 | 27 | 41 | 62 | 31 | 29 | 44 | 22 | 11 | 31 | 47 | 71 | 107 |
| 41 | 62 | 31 | 47 | 43 | 65 | 98 | 49 | 45 | 68 | 34 | 17 | 47 | 71 | 107 | 161 |
| 57 | 86 | 43 | 65 | 59 | 89 | 134 | 67 | 61 | 92 | 46 | 23 | 63 | 95 | 143 | 215 |
| 73 | 110 | 55 | 83 | 75 | 113 | 170 | 85 | 77 | 116 | 58 | 29 | 79 | 119 | 179 | 269 |
| | | | | | | | | | | | | | | | |

*Figure 25:* La super-décomposition de la polychromatrice $\mathbb{M}^{TT}(1,3)$ en $2^{3}$ chromologues différents

**Notations 9.1**

On désigne par

-$a_n(P_{1,k})$ le nombre des suites de type $A_n$ dans la chromatrice $\mathbb{M}^T(P_{1,k}, n)$(ou bien le chromoforme $Z^T(P_{1,k}, n)$).

-$b_n(P_{1,k})$ le nombre des suites de type $B_n$ dans la chromatrice $\mathbb{M}^T(P_{1,k}, n)$(ou encore dans le chromoforme $Z^T(P_{1,k}, n)$).

-$f_A(P_{1,k}, n)$ la proportion des suites de type $A_n$ dans la chromatrice $\mathbb{M}^T(P_{1,k}, n)$(ou bien le chromoforme $Z^T(P_{1,k}, n)$).

-$f_B(P_{1,k}, n)$ la proportion des suites de type $B_n$ dans la chromatrice $\mathbb{M}^T(P_{1,k}, n)$(ou encore dans le chromoforme $Z^T(P_{1,k}, n)$).

-$q_A(P_{1,1}, n)$ la proportion des suites de type $A_n$ dans la polychromatrice $\mathbb{M}^{TT}(P_{1,1}, n)$(ou bien dans le polychromoforme $\mathbb{Z}^T(P_{1,1}, n)$).

-$q_B(P_{1,1})$ la proportion des suites de type $B_n$ dans la polychromatrice $\mathbb{M}^{TT}(P_{1,1}, n)$(ou encore dans le polychromoforme $\mathbb{Z}^T(P_{1,1}, n)$).

-$z_A(n)$ le nombre des chromologues de type A dans l'ensemble des chromologues obtenus par la super- décomposition du polychromoforme d'ordre dimensionnel n.

-$z_B(n)$ le nombre des chromologues de type B dans l'ensemble des chromologues obtenus par la super- décomposition du polychromoforme d'ordre dimensionnel n.

On peut admettre les égalités suivantes qui traduisent l'équivalence entre les polychromologues et les polychromoformes.

$$\begin{cases} z_A(n) = h_A(n) \\ z_B(n) = h_B(n) \end{cases} \qquad \text{(e 9.12)}$$

**Corollaire 9.1**

Tous les chromoformes parfaits d'un même polychromoforme parfait de Collatz d'ordre dimensionnel n possèdent le même nombre des suites de type $A_n$ et le même nombre des suites de type $B_n$ c'est à dire que :

$$\begin{cases} a_n(P_{1,1}) = a_n(P_{2,1}) = \cdots = a_n(P_{j,1}) = \cdots = \cdots \\ b_n(P_{1,1}) = b_n(P_{2,1}) = \cdots = b_n(P_{j,1}) = \cdots = \cdots \end{cases} \qquad \text{(e 9.13)}$$

On peut remplacer les notations $a_n(P_{k,1})$ par $a_n$ et $b_n(P_{k,1})$ par $b_n$ qui désignent le nombre des suites de type $A_n$ (ou bien $B_n$) dans une chromatrice quelconque d'ordre dimensionnel un entier naturel non nul n.

**Démonstration**

Considérons les deux matrices $\mathbb{M}^T(P_{1,1}, n)$ et $\mathbb{M}^T(P_{i,1}, n)$. On sait que ces deux matrices sont isochromatiques c'est à dire que leurs suites sont deux a deux isochromatiques. On sait aussi que deux suites isochromatiques admettant le même coefficient transformationnel principal global et par suite elles sont de même type.

**Corollaire 9.2**

La relation entre le nombre des suites de type $B_n$ et le nombre des suites de type $A_n$ dans une même chromatrice d'ordre dimensionnel n est comme suit :

$$b_n = 2^n - a_n \qquad \text{(e 9.14)}$$

**Théorème 9.2**

Le nombre des chromologues de type $A_n$ dans le polychromologue $\mathbb{H}^T(P_{1,1}, n)$ est égal au nombre des suites de type $A_n$ dans le premier chromoforme $Z^T(P_{1,1}, n)$ et Le nombre des chromologues de type $B_n$ dans ce même polychromologue est égal au nombre des suites de type $B_n$ dans ce premier chromoforme parfait $Z^T(P_{1,1}, n)$ ceci se traduit par les équations suivantes :

$$\begin{cases} h_A(n) = a_n \\ h_B(n) = b_n \end{cases} \qquad \text{(e 9.15)}$$

**Démonstration**

On a montré que chaque suite du premier chromoforme $Z^T(P_{1,1}, n)$ représente la première suite d'un chromologue donnée du polychromologue $\mathbb{H}^T(P_{1,1}, n)$ donc nécessairement le nombre des suites de type $A_n$ dans ce premier chromoforme correspond au nombre des chromologues de type $A_n$ dans le polychromologue $\mathbb{H}^T(P_{1,1}, n)$.

**Corollaire 9.3**

La proportion des chromologue de type $A_n$ dans le polychromologue $\mathbb{H}^T(P_{1,1}, n)$ est égale à la proportion de suites de type $A_n$ dans le chromoforme $Z^T(P_{1,1}, n)$. De même la proportion des chromologues de type $B_n$ dans le polychromologue $\mathbb{H}^T(P_{1,1}, n)$ est égale à la proportion des suites de type $B_n$ dans le chromoforme $Z^T(P_{1,1}, n)$ ce qui se traduit par les deux équations suivantes :

$$\begin{cases} r_A(n) = f_A(n) \\ r_B(n) = f_B(n) \end{cases} \qquad \text{(e 9.15)}$$

**Corollaire 9.4**

La proportion des chromologue de type $S_n^+$ dans le polychromologue $\mathbb{H}^T(P_{1,1}, n)$ est égale à la proportion de suites de type $S_n^+$ dans le polychromoforme $\mathbb{Z}^T(P_{1,1}, n)$. De même la proportion des chromologues de type $S_n^-$ dans le polychromologue $\mathbb{H}^T(P_{1,1}, n)$ est égale à la proportion des suites de type $S_n^-$ dans le polychromoforme $\mathbb{Z}^T(P_{1,1}, n)$ ce qui se traduit par les deux équations suivantes :

$$\begin{cases} r^+(n) = f_n^+ \\ r^-(n) = f_n^- \end{cases} \qquad \text{(e 9.16)}$$

Puisque un polychromoforme peut être décomposée entièrement en un polychromologue donc les polychromologues et les polychromologues sont deux représentations ou des configurations d'un même ensemble infini des suites des Collatz.

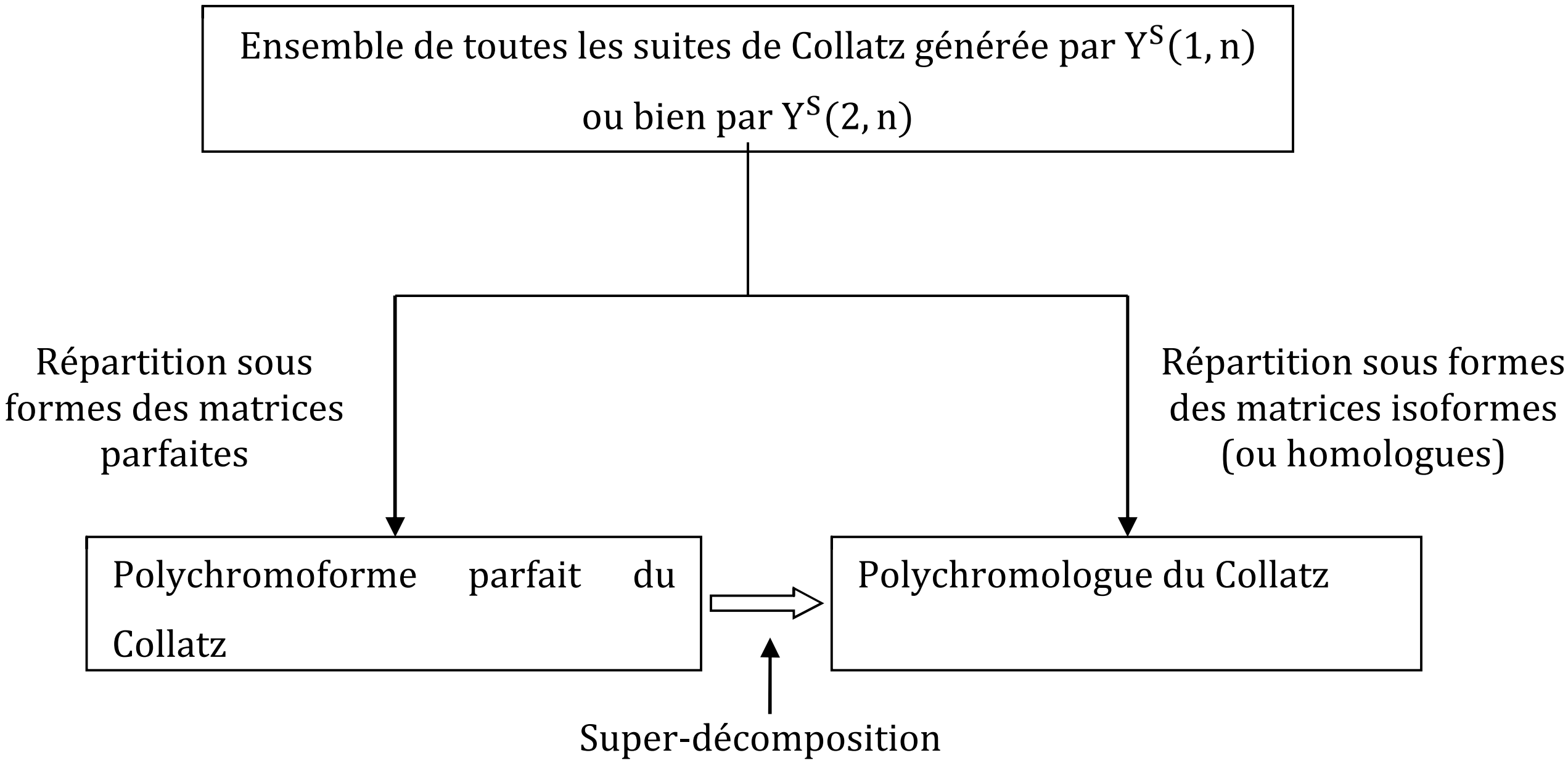

*Figure 26:* Equivalence entre les polychromologues et les polychromoformes de Collatz

Ces suites sont regroupes sous formes des matrices parfaites dans le cas d'un polychromoforme alors pour la deuxième ces même suites sont regroupées sous formes des chromologues ce qui nous permet de conclure que les proportions des suites de type $S_n^-$ sont égales et aussi évidement les deux proportions de type $S_n^+$ dans le polychromoforme et le polychromologue sont aussi égales.

**Corollaire 9.5**

On considère le polychromologue $\mathbb{H}^T(P_{1,1}, n)$ Les deux proportions des chromologues de types $A_n$ et $B_n$ vérifient les limites suivantes lorsque l'ordre dimensionnel des polychromologue tend vers l'infini :

$$\begin{cases} \lim_{n \to +\infty} r_A(n) = 0 \\ \lim_{n \to +\infty} r_B(n) = 1 \end{cases} \qquad \text{(e 9.17)}$$

**Démonstration :**

C'est conséquence directe du corollaire 9.3 et du théorème 8.2

**Théorème 9.3**

Soit $\mathbb{M}^T(P, 2n)$ la chromatrice parfaite de Collatz d'ordre dimensionnel 2n, on désigne par $f_{2n}^-$ la proportion des suites décroissantes de Collatz et par $f_{2n}^+$ la proportion des suites croissantes dans la chromatrice considérée alors $f_{2n}^+$ et $f_{2n}^-$ vérifient les limites suivantes :

$$\begin{cases} \lim_{n \to +\infty} f_{2n}^+ = 0 \\ \lim_{n \to +\infty} f_{2n}^- = 1 \end{cases} \qquad \text{(e 9.18)}$$

Autrement, lorsque l'ordre dimensionnel de la chromatrice tend vers l'infini la proportion des suites strictement croissantes tend vers 0 alors que la proportion des suites décroissantes tend vers 1 donc si on considère une chromatrice qui contient des suites de Collatz des longueurs infinies alors la quasi-totalité des ces suites sont décroissantes.

**Démonstration**

On sait que les deux proportions des deux types de chromologue $B_n$ et $A_n$dans le polychromologue est égale à la proportion des suites de type b dans la distribution parfaite de Collatz $Z^T(P_{1,1}, n)$

$$r^+(n) = r_A(n)$$

La proportion des suites abaissantes dans un chromoforme complet d'ordre dimensionnel égal à n :

$$r^{-}(n) = 1 - r^{+}(n) = 1 - \frac{h_A(n)}{2^n}$$

On sait que:

$$\begin{cases} r_A(n) = f_A(n) \\ r_B(n) = f_B(n) \end{cases}$$

On peut conclure que :

$$r^{+}(2n) = \lim_{n\to+\infty} r_A(2n) \tag{e 9.19}$$

$$= \lim_{n\to+\infty} \frac{h_A(2n)}{2^{2n}}$$

$$= \lim_{n\to+\infty} \frac{a_{2n}}{2^{2n}}$$

Ce qui nous permet de déduire :

$$\lim_{n\to+\infty} f_{2n}^{-} = \lim_{n\to+\infty} r^{-}(2n) = 1 \tag{e 9.20}$$

# 10 Le Phénomène du glissement séquentiel structurel des suites de Collatz

**Notations 10.1**

L'ensemble constitué par tous les permutations binaires avec répétition et de longueurs infinies est noté $\mathbb{D}^b$. Un élément de cet ensemble est noté $D_i^b$ avec i un entier naturel non nul quelconque.

$$\mathbb{D}^{\mathbf{b}} = \{D_1^b, D_2^b, D_3^b, D_4^b, \dots\} \tag{e 10.1}$$

-L'ensemble constitué par tous les arrangements avec répétition de longueurs infinies de deux couleurs (blanc et bleu) est noté $\mathbb{D}^{\mathbf{c}}$. Un élément de l'ensemble est noté $D_i^c$ avec i un entier naturel non nul quelconque.

$$\mathbb{D}^{\mathbf{c}} = \{D_1^c, D_2^c, D_3^c, \dots\} \tag{e 10.2}$$

Les relations entre $\mathbb{V}^b(n)$ et $\mathbb{S}^b$ et entre $\mathbb{W}^c(n)$ et $\Omega^{\mathbf{c}}$ se traduisent par les limites suivantes :

$$\lim_{n\to+\infty} \mathbb{V}^b(n) = \mathbb{D}^b \text{ et } \lim_{n\to+\infty} \mathbb{W}^c(n) = \mathbb{D}^c \tag{e 10.3}$$

## 10.1 Processus inverse de Collatz: La conversion inverse des distributions structurelles

Chaque suite générée de longueur infinie possède une unique distribution structurelle de longueur infinie que se comporte comme un arrangement binaire avec répétition de longueur infinie**.** On définie l'application comme suit :

$$T^{\gamma}:\ CLZ^{R} \to CLZ^{S}$$

$$\tilde{S}^R(P) \rightarrow \tilde{L}^B(P) \qquad \text{(e 10.4)}$$

Donc on fait convertir une suite numérique de Collatz en une suite structurelle qui décrit la distribution des entiers pairs et impairs dans la suite considérée.

Suite numérique → suite structurelle (arrangement binaire avec répétition)

Le processus inverse consiste a considéré un arrangement binaire avec répétition de longueur infinie qui peut se représenter sous forme d'une ligne composée d'une infinité des cases, la question que se pose existe-t-il un entier naturel non nul P qui peut générer une suite de Collatz de longueur infinie admettant une distribution structurelle identique à l'arrangement considéré.

Suite structurelle (de longueur infinie) → suite numérique de Collatz

En général, on peut distinguer deux cas différents :

-Si elle existe une suite générée $\tilde{S}^R(P)$ de longueur infinie tel que :

$$D_i^b = \tilde{L}^B(P) \qquad (10.5)$$

On dit que la structure binaire $D_i^b$ est réalisable dans $\mathbb{N}^*$ et le générateur P de la suite considérée constitue aussi le générateur du $D_i^b$.

-Si elle n'existe pas aucune suite générée de longueur infinie admettant une distribution structurelle identique à l'arrangement considéré on dit que cet arrangement (ou cette structure binaire) est non réalisable dans $\mathbb{N}^*$ ceci se traduit par:

$$D_i^B \neq \tilde{L}^B(P)\ \forall\ P \in \mathbb{N}^* \qquad (10.6)$$

**Définition 10.1**

Notion des chromolignes et des structures binaires

Un chromoligne est un arrangement chromatique de longueur infinie, elle est représentée sous forme d'une ligne à une infinité de cases, les cases sont colorées par deux couleurs dans un ordre quelconque. Les chromolignes sont donc les éléments de l'ensemble $\mathbb{D}^c$ déjà définie alors que les suites binaire infinie (appelées aussi structures binaires complètes ou polyformes binaires) sont les éléments de l'ensemble $\mathbb{D}^b$.

Puisque on ne peut pas représenter toutes les cases, on ajoute une case allongée à la fin de la ligne peut être colorée par une couleur différente.

Le chromoligne peut être remplacée par la suite binaire qui lui correspond, dans une case colorée en bleu on place 1 alors que pour une case colorée en blanc on met 0.

**Exemple 10.1**

Exemple d’une chromoligne et le vecteur de parité qui lui correspond:

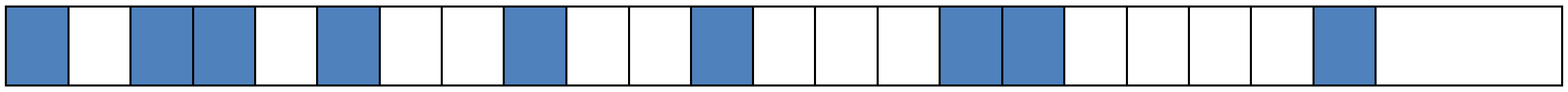

| | | | | | | | | | | | | | | | | | | | | | | |
|---|---|---|---|---|---|---|---|---|---|---|---|---|---|---|---|---|---|---|---|---|---|---|
| 1 | 0 | 1 | 1 | 0 | 1 | 0 | 0 | 1 | 0 | 0 | 1 | 0 | 0 | 0 | 1 | 1 | 0 | 0 | 0 | 0 | 1 | |

La suite binaire de longueur infinie est appelée une structure binaire infinie (ou complète). Les suites binaires infinies et les chromolignes sont donc des arrangements avec répétition de deux objets et de longueurs infinies. Ils sont appelés les polyformes autrement un polyforme peut designer un chromoligne ou bien une structure binaire de longueur infinie.

**Définition 10.2**

Les séquences d’ordre n d’une chromoligne donnée:

Une partie de longueur n constituée par les n premières cases d’un chromoligne donné ou bien d’une suite binaire infinie donnée est appelée séquence d’ordre n de cette suite ou de ce chromoligne.

La séquence chromatique de longueur n de d’une chromoligne D considérée est notée $w(D, n)$.

La séquence binaire de longueur n de la suite binaire infinie considérée est notée $v(L, n)$

**Exemple 10.2**

On considère le chromoligne suivant (supposé de longueur infinie)

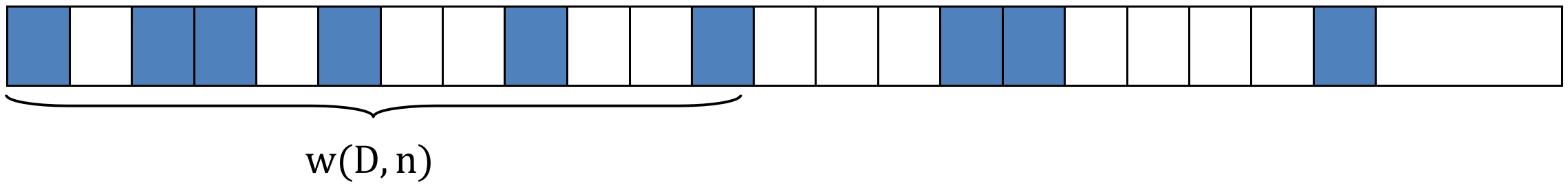

| | | | | | | | | | | | | | | | | | | | | | | |
|---|---|---|---|---|---|---|---|---|---|---|---|---|---|---|---|---|---|---|---|---|---|---|
| 1 | 0 | 1 | 1 | 0 | 1 | 0 | 0 | 1 | 0 | 0 | 1 | 0 | 0 | 0 | 1 | 1 | 0 | 0 | 0 | 0 | 1 | |

$v(L, n)$

Le tableau ci-dessous représente un ensemble des séquences extraites de la suite binaire infinie considérée :

*Tableau 9 :* Exemple d'une Série séquentielle relative au chromoligne considérée.

| | | | | | | | | | | | | | | | |
|---|---|---|---|---|---|---|---|---|---|---|---|---|---|---|---|
| v(L, 3) | 1 | 0 | 1 | | | | | | | | | | | | |
| v(L, 6) | 1 | 0 | 1 | 1 | 0 | 1 | | | | | | | | | |
| v(L, 9) | 1 | 0 | 1 | 1 | 0 | 1 | 0 | 1 | 1 | | | | | | |
| v(L, 12) | 1 | 0 | 1 | 1 | 0 | 1 | 0 | 1 | 1 | 1 | 1 | 0 | | | |
| v(L, 15) | 1 | 0 | 1 | 1 | 0 | 1 | 0 | 1 | 1 | 1 | 1 | 0 | 1 | 1 | 1 |

## 10.2 Classement des séquences structurelles

On sait qu'une séquence structurelle de longueur finie possède une infinité des conversions entières(ou représentations entières) qui constituent un chromologue de Collatz de même longueur et qui admettent une distribution structurelle linaire identique à cette séquence. On peut classer les séquences structurelles selon la nature de ces suites croissantes ou décroissante.

**Définition 10.3**

Les séquences structurelles de type $D^{+}$

Une séquence structurelle est dite de type $D^{+}$ si elle correspond à une distribution structurelle d'une suite strictement croissante de Collatz.

**Exemple 10.3**

On considère la distribution chromatique suivante:

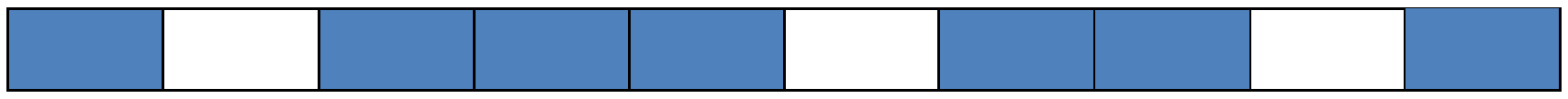

Cette distribution est type $D^{+}$ en effet elle admet une représentation numérique qui correspond à la suite croissante $S^{r}(37,9)$ représentée ci-dessous :

| | | | | | | | | | |
|---|---|---|---|---|---|---|---|---|---|
| 137 | 206 | 103 | 155 | 233 | 350 | 175 | 263 | 398 | 199 |

**Définition 10.4**

Les séquences structurelles de type $D^{-}$

Une séquence structurelle est dite de type $D^{-}$ si sa conversion entière correspond à une suite décroissante de Collatz.

**Exemple 10.4**

On considère la séquence finie suivante:

| | | | | | | | | | |
|---|---|---|---|---|---|---|---|---|---|
| | | | | | | | | | |

Une conversion inverse en suite de Collatz de cet arrangement donne toujours une suite décroissante de Collatz. Par exemple une des ces conversion numérique de la séquence considérée est comme suit :

| | | | | | | | | | |
|---|---|---|---|---|---|---|---|---|---|
| 147 | 221 | 332 | 166 | 83 | 125 | 188 | 94 | 47 | 71 |

## 10.3 Les chromolignes convertibles et les chromolignes non convertibles dans $\mathbb{N}^*$

**Définition 10.5**

Application indicatrice relative aux chromolignes convertibles et aux chromolignes non convertibles:

On définit l'application $\mathcal{R}$e pour tous chromolignes de l'ensemble $\mathbb{D}^C$ comme suit :

$$\mathcal{R}e: \mathbb{D}^b \longrightarrow \{0,1\}$$

$$D_i^b \longrightarrow \mathcal{R}e\left(D_i^b\right) \qquad \text{(e 10.7)}$$

On peut distinguer deux cas différents représentés par le système suivant :

$$\mathcal{R}e\left(D_i^b\right) = \begin{cases} 1 & \text{si L est convertible dans } \mathbb{N}^* \\ 0 & \text{si non} \end{cases} \qquad \text{(e 10.8)}$$

**Définition 10.6**

Série séquentielle d'une chromoligne donnée

Soit W une chromoligne donnée on sait qu'il correspond à une ligne à une infinité des cases colorées par deux couleurs différentes dans un ordre quelconque. On considère une série constituée d'une infinité des séquences extraites de la chromoligne **D**. les séquences chromatiques ont longueurs différentes qu'on peut le trier dans un ordre croissant comme suit :

$$n_1 < n_2 < n_3 < \cdots < n_{k-1} < n_k < \cdots \qquad (e\ 10.9)$$

Chaque séquence est incluse dans une autre comme suit:

$$w(D, n_1) \subset w(D, n_2\ ) \subset w(D, n_3\ ) \subset \cdots \subset w(D, n_4) \subset \cdots \subset w(D, n_1) \subset \cdots \qquad \text{(e 10.10)}$$

On dit que la limite de cette série tend vers la chromoligne **D** lorsque la longueur des séquences tend vers l'infini ceci se traduit par l'équation suivante:

$$\lim_{n_i \to +\infty} w(L, n_i) = \mathbf{D} \quad \text{(e 10.11)}$$

On sait que chaque séquence correspond à un arrangement avec répétition de longueur fini bien déterminée donc il existe un entier naturel P qui peut générer une suite de même longueur que cette séquence et qui admet une distribution linéaire structurelle qui lui identique.

$$P_i = gen\big(w(\mathbf{D}, n_i\,)\big) \quad \text{(e 10.12)}$$

On dit que la chromoligne **D** considérée est non réalisable dans $\mathbb{N}^*$si la condition suivante est vérifiée :

$$\lim_{n \to +\infty} gen(w(\mathbf{D}, n_i\,)) = +\infty \quad \text{(e 10.13)}$$

C'est à dire qu'il n'admet pas un générateur fini dans $\mathbb{N}^*$ donc elle n'existe pas aucune suite générée de longueur infinie de générateur un entier naturel non nul P dont la distribution structurelle est identique à D ceci se traduit par:

$$D \neq \tilde{L}^B(P)\ \forall \text{ entier naturel non nul } P \quad \text{(e 10.14)}$$

**Exemple 10.5**

On considère la chromoligne infinie suivante qu'on le note L :

| | | | | | | | | | | | | | | | | | | | |
|---|---|---|---|---|---|---|---|---|---|---|---|---|---|---|---|---|---|---|---|
| 1 | 1 | 1 | 0 | 1 | 0 | 0 | 1 | 0 | 1 | 0 | 1 | 0 | 1 | 0 | 1 | 0 | 1 | 0 | 1 |

*Figure 27:* Exemple d'une chromoligne convertible en suite de Collatz

La chromoligne considérée correspond à la distribution structurelle linéaire de la suite générée de longueur infinie $\tilde{S}^R(7)$ suivante :

| | | | | | | | | | | | | | | | | | |
|---|---|---|---|---|---|---|---|---|---|---|---|---|---|---|---|---|---|
| 7 | 11 | 17 | 26 | 13 | 20 | 10 | 5 | 8 | 4 | 2 | 1 | 2 | 1 | 2 | 1 | 2 | 1 |

On peut écrire alors:

$$L = \tilde{L}^S(P)$$

L'entier 7 représente le générateur de la chromoligne infinie L et la suite $\tilde{S}^R(7)$ est sa représentation numérique (ou conversion entière) ce qui nous permet de conclure que:

$$\mathcal{Re}(L) = 1$$

**Exemple 10.6**

L'exemple le plus simple, des chromolignes non réalisables correspond a une suite binaire ne contient que le nombre 1 répété une infinité des fois. On peut vérifier qu'elle n'existe aucun entier naturel non nul qui peut générer une suite de Collatz de longueur infinie dont la distribution structurelle linaire correspond à cette chromoligne.

| | | | | | | | | | | | | | | | | | |
|---|---|---|---|---|---|---|---|---|---|---|---|---|---|---|---|---|---|
| 1 | 1 | 1 | | | | | | | | | | | 1 | 1 | 1 | | |

On sait que si la longueur de cet arrangement est égale à n donc $P = (2^{n+1} - 1)$ est un générateur de cet arrangement et la suite générée de longueur n notée $\tilde{S}^{r}(P, n)$ ne contient que des termes impairs ce qui signifie que sa distribution structurelle linéaire correspond à un arrangement avec répétition qui ne contient que l'entier 1. Si la longueur de cet arrangement tend vers l'infini donc le générateur est aussi devient infini donc on peut conclure qu'il n'existe pas aucun un entier naturel non nul bien déterminé qui peut générer une suite de longueur infinie qui admet cette chromoligne comme distribution structurelle linaire et par conséquent cet arrangement est non réalisable dans l'ensemble des entiers naturels non nuls.

| | | | | | | | | | | |
|---|---|---|---|---|---|---|---|---|---|---|
| $w(B, 2)$ | 7 | 11 | 17 | | | | | | | |
| $w(B, 4)$ | 31 | 47 | 71 | 107 | 161 | | | | | |
| $w(B, 6)$ | 127 | | | | | | | | | |
| $w(B, 8)$ | 511 | | | | | | | | | |
| $w(B, 9)$ | 1023 | | | | | | | | | |

*Figure 28:* Le glissement de la chromoligne bleue et variation de leurs différents générateurs minimums en fonction de sa longueur

**Exemple 10.7**

Le deuxième exemple des arrangements non réalisables (ou non convertibles) correspondent à un arrangement structurel cyclique. Une chromoligne cyclique est caractérisée par une périodicité structurelle, elle est composée d'un nombre infini d'une même séquence comme montre l'exemple suivant :

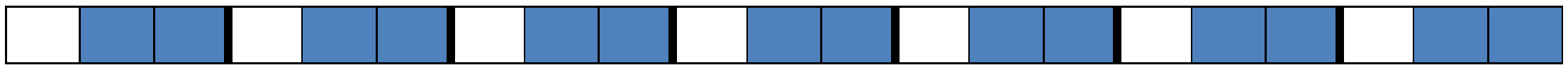

*Figure 29:* Une chromoligne cyclique constituée par un enchainement d'une même séquence

Cet arrangement est composé par la même séquence constituée trois cases : deux cases bleues et une case blanche (qui correspondent a deux entiers impairs suivi d'un entier pair) répétée une infinité des fois. Cette séquence constitue l'unité constructive répétitive de l'arrangement infini considérée. Elle est représentée ci-dessous :

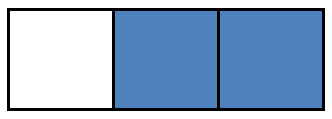

Séquence constitutive de longueur 3 d'une chromoligne cyclique infini

La longueur de la séquence constitutive est appelée la période structurelle de la suite cyclique elle est notée $\Delta_s$. Pour cet exemple :

$$\Delta_s = 3$$

Les cases des rangs $(i + k\Delta_s)$ ont la même couleur avec $i \in \{1,2,3\}$ et k un entier naturel quelconque

Une chromoligne est dite glissante ou bien une suite binaire infinie est dite glissante si leurs générateurs minimaux tend vers l'infini lorsque la longueur n des leurs séquences tend vers l'infini.

**Exemple 10.8**

Prenons le cas d'une chromoligne cyclique régulière constituée par une même séquence comme suit :

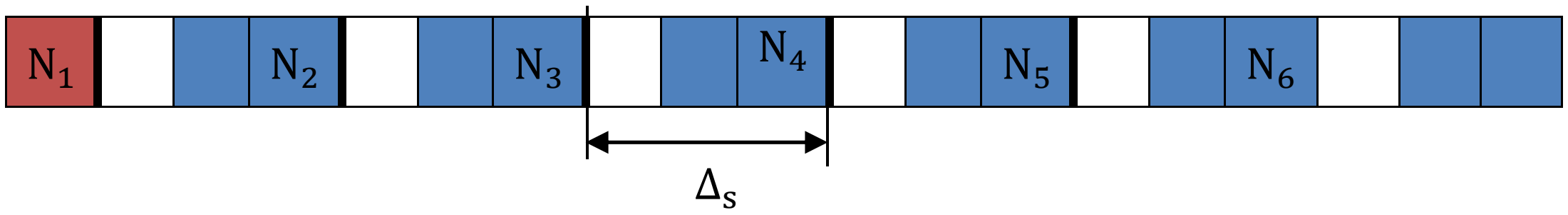

On considère une suite des termes extraits de la chromoligne considérée. Les termes extraits sont séparés par le même nombre des cases. Ce nombre correspond à la période structurelle du chromoforme cyclique. L'expression de chaque terme de la suite extraite en fonction du terme qui lui précède est comme suit :

$$N_2 = \frac{3^2}{2^3} N_1 + \frac{7}{8} \; ; \; N_3 = \frac{3^2}{2^3} N_2 + \frac{7}{8} ; \ldots , N_{k+1} = \frac{3^2}{2^3} N_k + \frac{7}{8}$$

Noter que le nombre k correspond au rang de la $k^{eme}$ séquence de la chromoligne considérée.

L'expression du terme du rang (k+1) de la suite extraite en fonction du premier terme est comme suit :

$$N_{k+1} = \frac{3^2}{2^3}(\frac{3^2}{2^3}N_{k-1} + \frac{7}{8}) + \frac{7}{8} = \frac{3^2}{2^3}(\frac{3^2}{2^3}(\frac{3^2}{2^3}N_{k-2} + \frac{7}{8}) + \frac{7}{8}) + \frac{7}{8}$$

$$= \left(\frac{3^2}{2^3}\right)^k N_1 + \frac{7}{8}\left(1 + \frac{3^2}{2^3} + \left(\frac{3^2}{2^3}\right)^2 + \cdots + \left(\frac{3^2}{2^3}\right)^k\right)$$

$$= \frac{3^{2k}}{2^{3k}} N_1 + \frac{7}{8}\frac{\left(\frac{3^2}{2^3}\right)^k - 1}{\frac{3^2}{2^3} - 1}$$

On peut déduire la relation suivante entre les deux termes de la suite extraite :

$$N_{k+1} = \left(\frac{3^2}{2^3}\right)^k \left(N_1 + 7\left(1 - \left(\frac{2^3}{3^2}\right)^k\right)\right) \qquad \text{(e 10.15)}$$

On représente ci-dessous les quatre premières suites cycliques de différentes longueurs :

| | | | |
|---|---|---|---|
| 9 | 14 | 7 | 11 |

| | | | | | | |
|---|---|---|---|---|---|---|
| 121 | 182 | 91 | 137 | 206 | 103 | 155 |

| | | | | | | | | | |
|---|---|---|---|---|---|---|---|---|---|
| 1017 | 1526 | 763 | 1145 | 1418 | 859 | 1289 | 1934 | 967 | 1451 |

| | | | | | | | | | | | | |
|---|---|---|---|---|---|---|---|---|---|---|---|---|
| 8185 | 12278 | 6139 | 9209 | 13814 | 6907 | 10361 | 15542 | 7771 | 11657 | 17486 | 8743 | 13115 |

Lorsque le nombre des séquences structurelles augmente le générateur aussi augmente et si ce nombre est infini, le générateur tend aussi vers l'infini ce qui signifie que la chromoligne cyclique constituée par un nombre infini de la séquence structurelle cyclique est non réalisable dans $\mathbb{N}^*$ autrement, il n'existe aucun entier naturel non nul P qui peut générer une suite de longueur infinie dont la distribution structurelle correspond à la chromoligne étudiée.

**Remarque 10.1**

En réalité, cette propriété est générale et elle est vérifiée pour toutes les distributions structurelles linéaires cycliques de dimensions infinies qui sont non convertibles et elles n'admettent pas des générateurs finies dans $\mathbb{N}^*$.

**Remarque 10.2**

Les seules suites cycliques qui peuvent admettre des générateurs finis sont celles qui présentent simultanément un comportement cyclique structurel et un comportement quantitatif (ou numérique) cyclique. La seule suite connue est la suite $S^R(1)$ puisque on n'a pas encore montré s'il existe d'autre cycles stables ou non.

| | | | | | | | | | | | |
|---|---|---|---|---|---|---|---|---|---|---|---|
| 1 | 2 | 1 | 2 | 1 | 2 | 1 | 2 | 1 | | | |

Dans le cas général on peut distinguer deux types de distributions structurelles cycliques

Cycle structurel dynamique : Il est caractérisé par une stabilité structurelle uniquement alors que les valeurs des différents termes de la suite varient continuellement.

**Exemple 10.9**

La séquence finie suivante ayant un cycle structurel dynamique correspond à la suite finie générée par 1017.

| | | | | | | | | | |
|---|---|---|---|---|---|---|---|---|---|
| 1017 | 1526 | 763 | 1145 | 1418 | 859 | 1289 | 1934 | 967 | 1451 |

Cycle stable : il est caractérisé par une périodicité numérique alors il est évident qu'il présente une périodicité structurelle aussi.

C'est l'exemple du cyclique stable correspond à la suite de Collatz générée par 1

**Théorème 10.1**

Il n'existe pas aucun entier naturel non nul qui peut générée une suite de Collatz de longueur infini

**Démonstration**

On montre par absurde que toutes les suites cycliques (comme dans le cas de l'exemple précédent) de longueurs infinies sont non réalisables (non convertible en suites de Collatz) dans $\mathbb{N}^*$.

On suppose qu'elle existe une suite infinie qui possède une distribution structurelle périodique de générateur fixe.

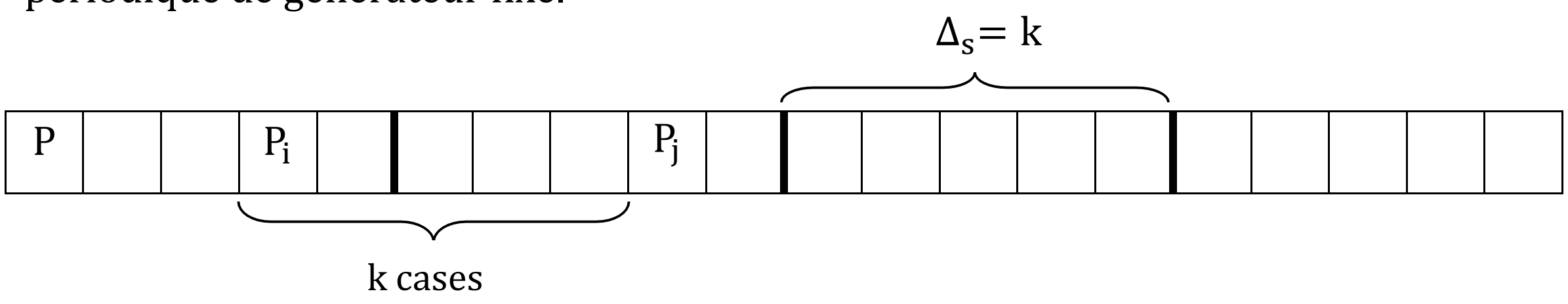

De plus on suppose que la séquence constitutive de cette suite est de longueur k (constituée de k cases) qui représente la période structurelle de la suite considérée.
Prenons deux termes quelconques de cette suite tel que et de dernier terme donc les deux termes ont la même parité. Les deux suites générées par sont isochromatiques elles possèdent la même distribution structurelle.
Si on considère la chromatrice d'ordre dimensionnel n tel que les deux sont deux termes de la suite génératrice de la chromatrice.
La matrice considérée contient deux suites générées isochromatiques absurde parce qu'on sait qu'une chromatrice parfaite de Collatz ne peut pas contenir deux suites isochromatiques.

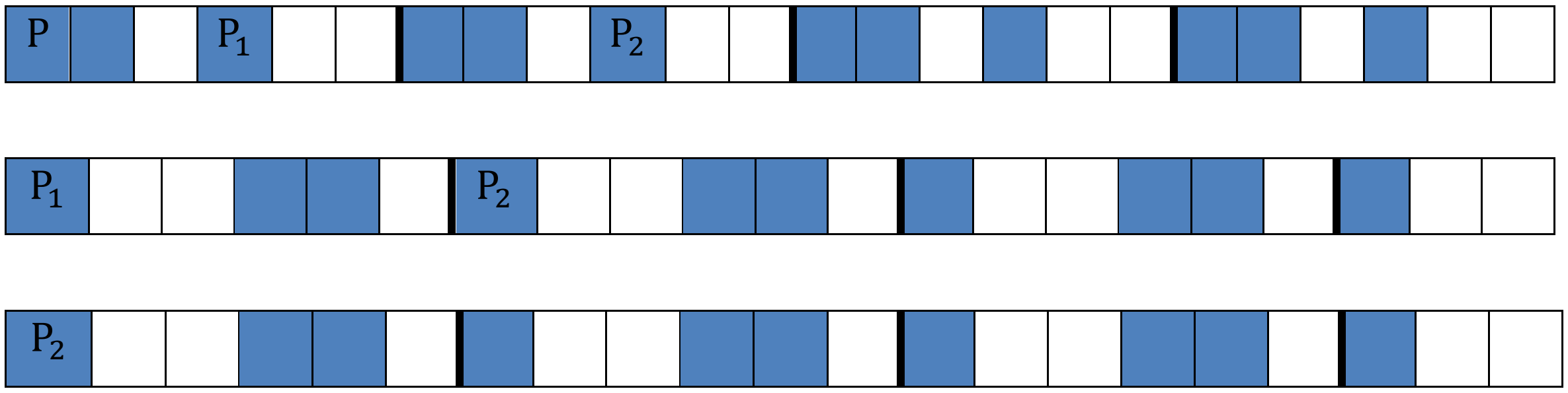

Comme P est un entier naturel fixe donc les deux termes sont constants et ils génèrent deux suites isochromatiques dans la même matrice parfaite à condition de choisir un ordre dimensionnel n de cette matrice qui vérifie la condition suivante :

$$(P_2 - P_1) < 2^n \quad \text{(e 10.16)}$$

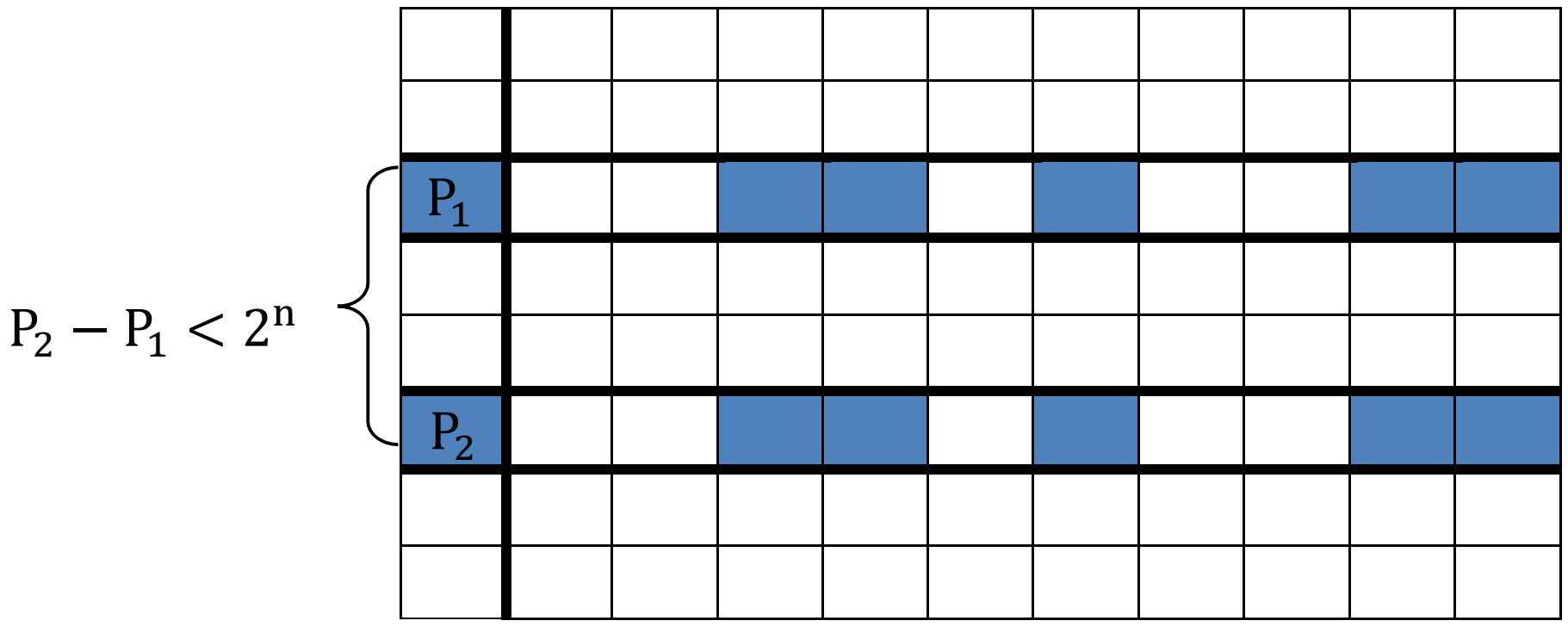

*Figure 30:* Impossibilité de la génération d'une suite cyclique de longueur infinie

Si cette condition est vérifiée alors les deux entiers appartenant a la suite génératrice de cette matrice. La possibilité de l'existence d'une suite cyclique de longueur infinie équivaut à l'existence de deux suites isoformes dans une même chromatrice parfaite du Collatz c'est qui est impossible.

**Remarque importante 10.3**

En réalité les chromolignes régulières ou cycliques ne sont pas seules qui sont non convertibles en suites numériques de Collatz. En effet on peut montrer que la quasi-totalité des chromolignes infinies (ou des polyformes) sont non convertibles en des suites de Collatz dans $\mathbb{N}^*$ c'est à dire que toutes les distributions linéaires structurelles relatives aux toutes les suites générées par tous les entiers naturels non nuls ne constituent qu'une proportion quasi-nulle de toutes les arrangements de l'ensemble $\mathbb{D}^{\mathrm{b}}$. La démonstration de cette propriété fait l'objet de la partie suivante.

**Théorème 10.2**

La quasi-totalité des arrangements binaires avec répétition de longueurs infinies de l'ensemble $\mathbb{D}^{\mathrm{b}}$ ne possèdent pas des conversions entières en suites de Collatz c'est-à-dire qu'ils ne possèdent pas des générateurs dans $\mathbb{N}^*$. En contre partie la proportion des arrangements convertibles en suites infinies de Collatz est quasi-nulle c'est à dire que toutes les distributions structurelles qui correspondent aux toutes les suites de Collatz des longueurs infinies et qui sont générées par tous les entiers naturels non nuls ne constituent qu'une proportion quasi-nulle de toutes les arrangements avec répétition de l'ensemble $\mathbb{D}^{\mathrm{b}}$.

**Démonstration**

On considère la matrice parfaite d'ordre n de Collatz $\mathbb{M}^{\mathrm{T}}(1,\mathrm{n})$ et on choisi un réel strictement positif quelconque r tel que :

$$0 < r < 1 \qquad (e\ 10.17)$$

On pose :

$$\mathrm{r_n} = 2^{\mathrm{n}}\mathrm{r} \qquad (\mathrm{e}\ 10.18)$$

On sait que la suite génératrice de cette matrice correspond à la suite parfaite $\mathrm{Y^S}(1,\mathrm{n})$

$$\mathrm{Y^S}(1,\mathrm{n}) = (1,3,5,7,\dots,2^{\mathrm{n}+1}-1)$$

On désigne par $P_{i_0}(n)$ le plus petit entier naturel strictement supérieur à $r_n$ alors cet entier vérifie l'équation suivante :

$$P_0(n) = E(r_n) + 1 \qquad \text{(e 10.19)}$$

On subdivise la matrice génératrice parfaite en deux parties comme suit :

$$(\underbrace{(1,3,5,7,\dots,P_0(n)-1)}_{Y^-},\underbrace{(P_0(n),P_0(n)+1,,\dots,2^{n+1}-1)}_{Y^+})$$

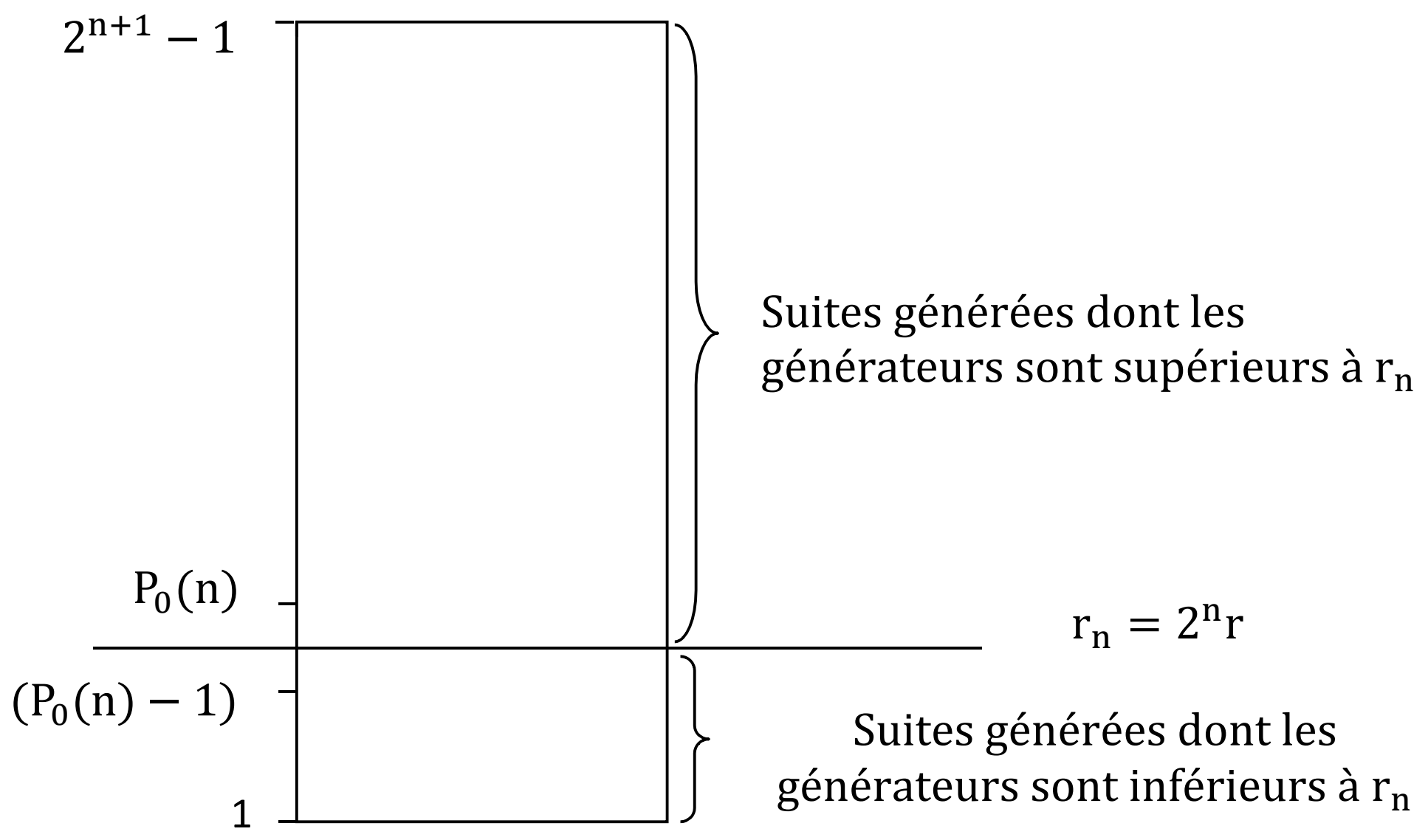

On subdivise la matrice parfaite de Collatz en deux zones comme suit:

Les suites qui ont des générateurs strictement supérieurs à $P_{i_0}(n)$ c'est à dire l'ensemble des suites générées par $Y^+$ dont les générateurs vérifient l'équation ci-dessous:

$$P_i \geq P_{i_0}(n) > r_n \quad \text{avec} \quad 2^n \geq i \geq i_0$$

Le nombre de ces suites est comme suit :

$$z^+(n) = 2^n - \frac{P_{i_0}(n) - 1}{2} \qquad \text{(e 10.20)}$$

La proportion des ces suites dans la matrice parfaite considérée est donnée par :

$$\rho^+(n) = \frac{z^+(n)}{2^n} = 1 - \frac{P_{i_0}(n) - 1}{2^{n+1}} \qquad \text{(e 10.21)}$$

Comme $(P_{i_0}(n) - 1) \leq r_n$ alors on peut écrire :

$$\left(1 - \frac{P_{i_0}(n) - 1}{2^{n+1}}\right) \geq (1 - \frac{r_n}{2^{n+1}})$$

Ceci équivaut à :

$$\rho^+(n) \geq 1 - \frac{r_n}{2^{n+1}}$$

Comme on a :

$$1 - \frac{r_n}{2^{n+1}} = 1 - \frac{r}{2}$$

On peut déduire que :

$$\rho^+(n) \geq 1 - \frac{r}{2} \quad \text{(e 10.22)}$$

Pour r constant, lorsque n tend vers l'infini, tous les générateurs des suites générées qui sont supérieurs ou égaux à $r_n$ tendent aussi vers l'infini ce qui signifie que toutes les distributions structurelles linéaires de longueurs infinies qui correspondent aux suites appartenant à cette zone sont non réalisables. Cette limite est comme suit :

$$\begin{aligned} \lim_{n\to+\infty} (P_i)_{i\geq i_0} &= \lim_{n\to+\infty} P_{i_0}(n) \\ &= \lim_{n\to+\infty} r_n \\ &= \lim_{n\to+\infty} 2^n r \\ &= +\infty \end{aligned} \quad \text{(e 10.23)}$$

Puisque cette propriété est vraie pour n'importe quelles valeurs de r strictement positif, on peut tendre r vers 0 dans ce cas la proportion de la zone supérieur tend vers 1. En effet, on peut écrire :

$$\lim_{r\to 0} \rho^+(n) \geq \lim_{r\to 0} \left(1 - \frac{r}{2}\right)$$

Comme on a :

$$\lim_{r\to 0}(1 - \frac{r}{2}) = 1$$

On peut conclure que :

$$\lim_{r\to 0} \rho^+(n) = 1 \quad \text{(e 10.24)}$$

Cette dernière limite correspond à la proportion des chromolignes de longueurs infinies qui sont non convertibles en suites de Collatz.

## Références